\documentclass[reqno]{amsart}
\usepackage{amsfonts}
\usepackage{amsmath}
\usepackage{amsthm, amscd}
\usepackage{mathtools}
\usepackage{amssymb}
\usepackage{mathrsfs}
\usepackage{graphicx}
\usepackage{caption}
\usepackage{subcaption}
\usepackage{float}
\usepackage{xypic}
\usepackage[abs]{overpic}
\usepackage[alphabetic]{amsrefs}
\usepackage{etex}
\usepackage{hyperref}
\hypersetup{
  colorlinks = true, 
  linkcolor = black,
  urlcolor=black
}

\usepackage{listings}
\allowdisplaybreaks

\theoremstyle{plain}\newtheorem{Theorem}{Theorem}[section]
\theoremstyle{plain}\newtheorem{Corollary}[Theorem]{Corollary}
\theoremstyle{plain}\newtheorem{Lemma}[Theorem]{Lemma}
\theoremstyle{plain}\newtheorem{Definition}[Theorem]{Definition}
\theoremstyle{plain}\newtheorem{Proposition}[Theorem]{Proposition}
\theoremstyle{plain}
\theoremstyle{plain}\newtheorem{Conjecture}[Theorem]{Conjecture}
\theoremstyle{plain}\newtheorem*{Theorem*}{Theorem}
\theoremstyle{plain}\newtheorem{Question}[Theorem]{Question}

\theoremstyle{remark}\newtheorem{remark}[Theorem]{Remark}
\theoremstyle{remark}\newtheorem{Example}[Theorem]{Example}
\theoremstyle{remark}

\theoremstyle{plain}
\makeatletter
\newtheorem*{rep@theorem}{\rep@title}
\newcommand{\newreptheorem}[2]{
\newenvironment{rep#1}[1]{
 \def\rep@title{#2 \ref{##1}}
 \begin{rep@theorem}}
 {\end{rep@theorem}}}
\makeatother
\newreptheorem{Theorem}{Theorem}
\newreptheorem{Proposition}{Proposition}
\newreptheorem{Corollary}{Corollary}

\numberwithin{equation}{section}

\usepackage{enumerate}
\usepackage{lipsum}

\usepackage{mathtools}
\usepackage{tikz-cd}

\DeclareMathOperator{\Imm}{Im}

\DeclareMathOperator{\id}{id}

\DeclareMathOperator{\Vol}{Vol}

\DeclareMathOperator{\Homeo}{Homeo}
\DeclareMathOperator{\Diff}{Diff}
\DeclareMathOperator{\BHomeo}{BHomeo}
\DeclareMathOperator{\BDiff}{BDiff}
\DeclareMathOperator{\BSO}{BSO(4)}
\DeclareMathOperator{\ESO}{ESO(4)}
\DeclareMathOperator{\BO}{BO(4)}

\DeclareMathOperator{\BTop}{BTop}
\DeclareMathOperator{\ETop}{ETop}
\DeclareMathOperator{\Top}{Top}
\DeclareMathOperator{\Fr}{Fr}
\DeclareMathOperator{\Sm}{Sm}
\DeclareMathOperator{\Bl}{Bl}
\DeclareMathOperator{\Dim}{Dim}
\DeclareMathOperator{\cech}{Cech}
\DeclareMathOperator{\even}{even}
\DeclareMathOperator{\FrT}{Fr^{Top}}
\DeclareMathOperator{\dom}{dom}
\DeclareMathOperator{\Inte}{Int}
\DeclareMathOperator{\pj}{pj}
\DeclareMathOperator{\Fra}{Fr}
\DeclareMathOperator{\ks}{ks}

\newcommand{\bQ}{\mathbb{Q}}
\newcommand{\bR}{\mathbb{R}}

\newcommand{\bZ}{\mathbb{Z}}
\newcommand{\bfn}{\mathbf{n}}
\newcommand{\bfp}{\mathbf{p}}

\newcommand{\bfk}{\mathbf{k}}

\newcommand{\cA}{\mathcal{A}}
\newcommand{\cG}{\mathcal{G}}
\newcommand{\cS}{\mathcal{S}}
\newcommand{\cN}{\mathcal{N}}
\newcommand{\cM}{\mathcal{M}}
\newcommand{\cP}{\mathcal{P}}
\newcommand{\cT}{\mathcal{T}}

\author{Jianfeng Lin and Yi Xie}
\title{Configuration space integrals and formal smooth structures}
\date{}

\begin{document}

\maketitle
\begin{abstract}
 Watanabe disproved the 4-dimensional Smale conjecture by constructing topologically trivial $D^{4}$-bundles over spheres and showing that they are smoothly nontrivial using configuration space integrals. In this paper, we define a new version of configuration space integrals that only relies on a formal smooth structure on the $D^{4}$-bundle (i.e., a vector bundle structure on the vertical tangent microbundle). It coincides with Watanabe's definition when the $D^{4}$-bundle is smooth. We obtain several applications. First, we give a lower bound (in terms of the graph homology) on the dimension of the rational homotopy and homology groups of $\Top(4)$ and $\Homeo(S^4)$ (the homeomorphism group of $\mathbb{R}^4$ and $S^4$). In particular, this implies that $\Top(4)$ and $\Homeo(S^4)$ are not rationally equivalent to any finite-dimensional CW complexes. Second, we discover a generalized Miller-Morita-Mumford class $\kappa_{\theta}(\pi)\in H^{3}(B;\mathbf{Q})$, which is defined for any topological 4-manifold bundle $X\to E\to B$. This class obstructs the existence of a formal smooth structure on the bundle. Third, we show that for any compact, orientable, smooth 4-manifold $X$ (possibly with boundary), the inclusion map from its diffeomorphism group to its homeomorphism group is not rationally $2$-connected (hence not a weak homotopy equivalence). This implies that the space of smooth structures on $X$ has a nontrivial rational homotopy group in dimension 2.
\end{abstract}

\section{Introduction}

The diffeomorphism group $\Diff(S^{n})$ plays an essential  role in many important problems in topology. For example, by the work of Smale \cite{Smale62} and Cerf \cite{Cerf1970}, smooth structures on $S^{n+1}$ one-one correspond to path components of $\Diff(S^n)$ when $n\geq 5$. Therefore, there have been numerous studies about the homotopy types of $\Diff(S^{n})$ and the closely related object $\Diff_{\partial}(D^{n})$ (the group of diffeomorphisms on $D^n$ relative to the boundary) \cite{FarrellHsiang78,Igusa88,Kupers-Randel-Williams20,KupersRandel-Williams23,Krannich22,Krannich-Randal-Williams,Wa18,Wa21,Weiss21}. (See \cite{Randel-Williams22} for a survey about recent progresses.) There is the following central question:
\begin{Question}[Generalized Smale Problem]\label{ques: generalized smale} Given $n\geq 1$, is the inclusion 
$O(n+1)\hookrightarrow \Diff(S^n)$ a homotopy equivalence?
\end{Question}
By the homotopy equivalence $\Diff(S^{n})\simeq \Diff_{\partial}(D^n)\times O(n+1)$, the answer to Question \ref{ques: generalized smale} is affirmative if and only if $\Diff_{\partial}(D^n)$ is contractible. This is true when $n\leq 3$: $n=1$ case  is straightforward. Smale \cite{Smale59} proved the case $n=2$. And Hatcher \cite{Hatcher83} proved the case $n=3$. (See \cite{Bamler23} for an alternative using Ricci flow, which applies to many other 3-manifolds.) For all $n\geq 5$, Question \ref{ques: generalized smale} has a negative answer. This follows from the works of Novikov \cite{Novikov65}, Antonelli-Burghelea-Khan \cite{Antonelli70}, Smale \cite{Smale61},
Cerf \cite{Cerf1970} and Kervarie-Milnor \cite{Kervarie63}. And Crowley-Schick \cite{Crowley13} further proved that for all $n\geq 7$, $\pi_{k}(\Diff_{\partial}(D^n))\neq 0$ for infinitely many $k$. The last case $n=4$ remains open until 2018, when Watanabe \cite{Wa18} proved that $\pi_{k}(\Diff_{\partial}(D^4))\otimes \bQ\neq 0$ for many $k$ including $1, 4, 8$.

The current paper reflects our attempt to understand and generalize Watanabe's result. We first briefly recall Watanabe's proof. Let $\mathcal{A}^{\even}_{k}$ be the graph homology, defined as the $\mathbb{Q}$-vector space spanned by trivalent graphs with $2k$ vertices modulo the IHX relation (see Definition \ref{defi: A_k}). Given an integer class $\eta\in \mathcal{A}^{\even}_{k}$, Watanabe constructed a boundary trivialized, framed smooth disk bundle 
\begin{equation}\label{eq：Watanabe family}
D^4\to E^{\eta}_{D}\xrightarrow{\pi^{\eta}}S^{k}.
\end{equation} 
This is done by generalizing Goussarov-Habiro's clasper surgery \cite{Goussarov99,Habiro2000} to higher dimensions. The bundle gives an element \[\alpha_{\eta}\in \pi_{k}(\BDiff^{\Fra}_{\partial}(D^4))\otimes \bQ\cong \pi_{k}(\BDiff_{\partial}(D^4))\otimes \bQ\cong \pi_{k-1}(\Diff_{\partial}(D^4))\otimes  \bQ.\]
Here $\BDiff^{\Fra}_{\partial}(D^4)$ denotes the classifying space of framed disk bundles (see (\ref{eq: framed moduli space})). By an higher-dimensional analogue of Kuperberg–Thurston and Lescop's computation of configuration space for homology 3-spheres \cite{Kuperberg99,Lescop04}, Watanabe computed that $I(\alpha_{\eta})=\eta$. Hence $\pi^{\eta}$ is a nontrivial bundle. Here \[I: \pi_{k}(\BDiff^{\Fra}_{\partial}(D^4))\otimes \bQ\to \cA^{\even}_{k}\]
is the map defined by evaluating the configuration space integral, an powerful invariant first proposed by Kontsevich \cite{Kontsevich94}) and further studied by many people \cite{Axelrod-Singer,Bar-Natan1995,Bott1998,BottTaubes94,Kuperberg99,Lescop04}. This argument actually shows that $I$ is surjective. Hence $\pi_{k-1}(\Diff_{\partial}(D^4))\otimes \bQ\neq 0$ whenever $\cA^{\even}_{k} \neq 0$. This holds for many values of $k$ including $2,5,9$. Note that Watanabe also proved similar results in higher dimensions \cite{Wa09,Wa09II,Wa21}.

By the Alexander trick, the homeomorphism group $\Homeo_{\partial}(D^4)$ is contractible. So $E^{\eta}_{D}$ is an ``exotic disk bundle'' (i.e., a topologically trivial bundle that is smoothly nontrivial). This is surprising because other than gauge theoretic invariants (Donaldson invariants \cite{Donaldson90}, Seiberg-Witten invariants \cite{Witten94}, Ozsv\'ath-Szab\'o mixed invariants \cite{Ozsvath06}), the only known invariants that can detect exotic phenomena on (families of) orientable 4-manifolds are Khovanov homology \cite{Piccirillo20,Hayden22} and configuration space integrals \cite{Wa18}.  Given this fact, it is natural to ask how  $I$ depends on smooth structures, and what's the difference between exotic phenomena detected by gauge theory and exotic phenomena detected by configuration space integrals. This problem has been studied in several papers. In particular, Chen \cite{Chen23} proved that a smooth structure affects the configuration space integral by affecting the stratification of the compactified 2-point configuration space bundle. And Knudsen-Kupers \cite{Knudsen-Kupers} raised the following two questions:

\begin{Question}\label{ques: IW factors through T-infty}
Does the configuration space integral factor through the limit of the embedding calculus Taylor tower \cite{Boavida18}?
\end{Question}

\begin{Question}\label{ques: IW factors through Top}
Does the configuration space integral only depend on formal smooth structures?
\end{Question}
Recall that any topological $n$-manifold $X$ has a \emph{topological tangent bundle} \[\mathbb{R}^n\to \cT X\to X,\]
defined as a suitable neighborhood of the diagonal $\Delta$ in $X\times X$. A formal smooth structure  on $X$ is a vector bundle structure on $\cT X$ and a framing on $X$ is a trivialization of $\cT X$. More generally, any topological bundle $X\to E\to B$ has a vertical topological tangent bundle 
 \[
 \mathbb{R}^n\to \cT^{v} E\to E.
 \] 
 We assume the bundle $E$ has a boundary trivialization $\nu(\partial M)\times B\hookrightarrow E$, where $\nu(\partial M)$ is a collar neighborhood of $\partial M$ in $M$. Then a formal smooth structure (resp. a framing) on the bundle $E\to B$ is a vector bundle structure (resp. a trivialization) on $\cT^{v} E$ that extends a fixed formal smooth structure (resp. trivialization) on $\nu(\partial M)$. 
 A topological bundle equipped with a formal smooth structure is called a formally smooth bundle. And a   topological bundle equipped with a framing is called a framed topological bundle. 
 
 The following result is the main theorem of our paper. It gives an affirmative answer to Question \ref{ques: IW factors through Top} for $I$ and provides evidence for Question \ref{ques: IW factors through T-infty}.
 \begin{Theorem}\label{thm: main}
 There is a topological version of configuration space integral $I^{\Top}$ that is defined for all framed topological disk bundles over general bases, as well as formally smooth bundles over spheres. It coincides with $I$ when the bundle is smooth and it induces a surjection 
 \begin{equation}\label{eq: I-top to Ak}
I^{\Top}: \pi_{k+4}(\Top(4))\otimes \bQ\to \cA^{\even}_{k}     
 \end{equation}
 that fits in to the commutative diagram 
\begin{equation}\label{eq: diagram I-top=I-W}
\begin{tikzcd}
 \pi_{k}(\BDiff^{\Fra}_{\partial}(D^4))\otimes \bQ\arrow[r,"m_*"] \arrow[rd,"I"]  & \pi_{k+4}(\Top(4))\otimes \bQ \arrow[d,"I^{\Top}"] & \\
   ~  & \cA^{\even}_{k}
\end{tikzcd}
\end{equation}
Here $\Top(4)$ is the group of homeomorphisms on $\mathbb{R}^4$ that fix the origin, and $m_*$ is induced by  the ``scanning map'' $m: \BDiff^{\Fra}_{\partial}(D^4)\to\Omega^{4} \Top(4)$ (see (\ref{eq: scanning map})). 
 \end{Theorem}
 By choosing a Riemannian metric, any smooth structure determines a formal smooth structure. In dimension $\neq 4$, the smoothing theory of Burghelea-Lashof \cite{Burghelea-Lashof} and Kirby-Siebenmann \cite[Essay V]{Kirby-Siebenmann} implies that any formal smooth structure determines an essentially unique smooth structure. So a result analogous to Theorem \ref{thm: main} has no content in dimension $\neq 4$. However, it's known that smoothing theory fails in dimension $4$ (e.g. any two smooth structures on a simply-connected 4-manifold are isomorphic as formally smooth structures). And Theorem \ref{thm: main} has concrete topological applications, which we discuss now. 
 
 First, as pointed out in \cite[Proposition 4.7]{Knudsen-Kupers}, a result like Theorem \ref{thm: 
 main} immediately leads to a better understanding of the topology of $\Top(4)$ and $\Homeo(S^4)$.

\begin{Corollary} \label{cor: Top and Homeo} The following result holds.
\begin{enumerate}
    \item $\dim(\pi_{k}(\Top(4))\otimes \bQ)\geq \dim \cA^{\even}_{k-4}$ for any $k\geq 0$.
\item $\dim(\pi_{k}(\Homeo(S^4))\otimes \bQ)\geq \dim\cA^{\even}_{k-4}$ for any $k\geq 0$.
\item The groups $H^{k}(\Top(4);\bQ)$ and $H^{k}(\Homeo(S^4);\bQ)$ are nonvanishing for infinitely many $k$. More precisely, consider the graded vector space \[V=\oplus_{k\in \mathbb{N}} V_{k},\text{ where } V_{k}:=(\cA^{\even}_{k-4})^*.\] Let $\operatorname{Sym}(V)$ be the symmetric algebra generated by $V$. Then both $H^{*}(\Top(4);\bQ)$ and $H^{*}(\Homeo(4);\bQ)$ contain  $\operatorname{Sym}(V)$ as a subring. 

\item $\Top(4)$ and $\Homeo(S^4)$ are are not rationally homotopy equivalent to any finite dimensional CW complexes. In particular $\Top(4)$ is not rational homotopy equivalent to $O(4)$ and $\Homeo(S^4)$ is not rational homotopy equivalent to $O(5)$. 
 \end{enumerate}
\end{Corollary}
\begin{remark} It is proved in \cite{Randall94} that $\Homeo(S^4)$ is not homotopy equivalent to a finite CW complex. The same argument applies to $\Top(4)$.
\end{remark}

The second application concerns  characteristic classes of topological bundles. Given any homology class $c\in H^{k}(\BTop(4))$, one can associate a characteristic class $\kappa_{c}$ of topological bundles \cite{Ebert2014}. This class $\kappa_{c}$, called the generalized MMM (Miller-Morita-Mumford) class, is constructed as follows.  Let $X\to E\to B$ be an oriented topological bundle whose fiber $X$ a closed 4-manifold. The topological vertical tangent bundle $\cT^{v}E$ is classified by a map 
\[
f_{\pi}: E\to \BTop(4)
\]
So we define \[\kappa_{c}(\pi):=\pi_{!}\circ f^{*}_{\pi}(c)\in H^{k-4}(B).\]
Here $\pi_{!}$ denotes integration along fibers. Taking $c$ to be the generator of \[\ker(H^{4}(\BTop(4);\mathbb{F}_2)\to H^{4}(\BO;\mathbb{F}_2))\cong \mathbb{Z}/2,\]
then $\kappa_{c}$ recovers the Kirby-Siebenmann invariant $\ks(X)\in\mathbb{Z}/2$ of the fiber. Using Theorem \ref{thm: main}, we are able to find a new nontrivial MMM class, which can be thought of as a ``higher Kirby-Siebenmann invariant'' for topological bundles. 

\begin{Theorem}\label{thm: MMM class} The configuration space integral for the complete graph with 4 vertices gives a canonical nonzero element $\theta\in H^{7}(\BTop(4);\bQ)$. The corresponding MMM class $\kappa_{\theta}(\pi)\in H^3(B;\bQ)$ is defined for all oriented topological bundle $X\to E\xrightarrow{\pi}B$ whose fiber is a closed, topological 4-manifold $X$. This class is nontrivial for a general bundle. It is natural under pullbacks and it vanishes if the bundle has a formal smooth structure.  
\end{Theorem}
\begin{remark}
In the proof Theorem \ref{thm: MMM class}, we will construct a topological bundle \[
X\# N(S^2\times S^2)\to E\xrightarrow{\pi} Y
\]
such that $\kappa_{\theta}(\pi)\neq 0$. Here $N$ is a fixed number and $Y$ is a fixed 3-manifold, and $X$ can be any smooth closed 4-manifold. The bundle is constructed by combining Watanabe's construction with Galatius-Randal-Williams' calculation of the stable homology of moduli spaces of 4-manifolds \cite{Galatius2017}. In particular, this implies the class 
\[
\kappa_{\theta}\in
H^{3}(\lim_{\substack{\rightarrow\\n}}\cM^{t}(\mathring{X}\# n(S^2\times S^2));\bQ) \]
of the stable moduli space is nonvanishing for any smooth $X$. Here $\cM^{t}$ denotes the moduli space in the topological category (see (\ref{eq: M^t})) and  $\mathring{X}=X\setminus \Inte(D^4)$.
\end{remark}

\begin{remark}
Weiss \cite{Weiss21} proved a surprising result that the topological Pontryagin classes do not satisfy the elementary relations of smooth  Pontryagin classes. Galatius and Randal-Williams  \cite{Galatius2023} further proved that the topological Pontryagin classes are algebraically independent in dimension $\geq 4$. This implies that the map
\[
\mathbb{Q}[e,p_{1},p_{2},\cdots]\to H^{*}(\BTop(4);\mathbb{Z})\otimes \mathbb{Q}
\]
is injective. It's unknown whether the same result holds for $H^{*}(\BTop(4);\mathbb{Q})$.
\end{remark}

The third application concerns the comparison of the diffeomorphism group and the homeomorphism group of a general manifold. There is the following question.
\begin{Question}\label{ques: Diff vs Homeo}
Given a smooth, compact manifold $X$, is the inclusion map 
$i^{X}: \Diff_{\partial}(X)\to \Homeo_{\partial}(X)$ a weak homotopy equivalence?
\end{Question}
We first discuss the motivation of Question \ref{ques: Diff vs Homeo}. Consider the classifying spaces 
\begin{equation}\label{eq: M^t}
\cM^{t}(X):=\BHomeo_{\partial}(X)    
\end{equation}
and 
\[
\cM^{s}(X):=\BDiff_{\partial}(X).
\]
$\cM^{t}(X)$ and $\cM^{s}(X)$ are called the moduli spaces of $X$ in the topological category and the smooth category respectively, because they classify families of topological/smooth manifolds (i.e. bundles) over a given base. The map $i^{X}$ induces the forgetful map 
\begin{equation}\label{eq: forgetful map top to smooth}
\rho^{X}:\cM^{t}(X)\to \cM^{s}(X).    
\end{equation}
The homotopy fiber of this map is denoted by $\Sm(X)$. This space is called \emph{the space of smooth structures} diffeomorphic to $X$ and it measures the difference between the two moduli spaces. Question \ref{ques: Diff vs Homeo} is asking about the uniqueness of smooth structures in a family setting, i.e., whether $\rho^{X}$ is a weak homotopy equivalence (equivalently, whether $\Sm(X)$ is weakly contractible).

In dimension $\leq 3$, $i^{X}$ is always a weak homotopy equivalence \cite{Earle,Bing83}. This is the family version of the statement that any manifold of dimension $\leq 3$ admits an essentially unique smooth structure. 
In dimension $4$, the study of Question \ref{ques: Diff vs Homeo} has a long history since the work of Donaldson \cite{Donaldson90} and Ruberman \cite{Ruberman98}. Both gauge theoretic and non-gauge theoretic techniques have been used to establish many examples for which $i^{X}$ is not a homotopy equivalence \cite{BaragliaKonno22,Kato2021,Baraglia21,KronheimerMrowkaDehn,Lin2023, Wa18,RubermanInpreparation,RubermanInpreparation2}.

The following result is an application of Theorem \ref{thm: main}. It answers  Question \ref{ques: Diff vs Homeo} negatively for all orientable 4-manifolds.
\begin{Theorem}\label{thm: Sm noncontractible}
Let $X$ be any compact, orientable 4-manifold with (possibly empty) boundary. Then $\pi_{2}(\Sm(X))\otimes \mathbb{Q}\neq 0$. If we further assume that $X$ has nonempty boundary or vanishing signature, then the inequality
\[
\dim(\pi_{k}(\Sm(X))\otimes \bQ)\geq \dim \cA^{\even}_{k} 
\]
holds for arbitrary $k$.
\end{Theorem}

\begin{Corollary}\label{cor: Diff vs Homeo} For any compact, orientable 4-manifold $X$, the map $i^{X}$ is \text{not} rationally 2-connected (so not a weak homotopy equivalence). 
\end{Corollary}

Theorem \ref{thm: Sm noncontractible} is proved by studying the smooth bundle $X\to E^{\eta}_{X}\to S^{k}$ obtained by gluing the $D^4$-bundle (\ref{eq：Watanabe family}) with $S^{k}\times \mathring{X}$. It suggests the following conjecture:
\begin{Conjecture}\label{conj: E-X nontrivial}
For any oriented smooth 4-manifold $X$ and any nonzero class $\eta\in \cA^{\even}_{k}$, the bundle $E^{\eta}_{X}\to S^k$ is smoothly nontrivial. In particular, it gives a nontorsion element \[\alpha^{X}_{\eta}\in \ker(\pi_{k}(\cM^{s}(X))\to \pi_{k}(\cM^{t}(X))).\]
\end{Conjecture}
\begin{remark}
One may try to approach Conjecture \ref{conj: E-X nontrivial} by developing  configuration space integrals for general manifold bundles. See \cite{Wa20} for relevant works in the case of homology $S^1\times S^3$s, using twisted coefficients. However, for simply connected $X$ with $b_2\neq 0$ (e.g. $S^2\times S^2$ or $\mathbb{CP}^2$), there is an obstruction for existence of a propagator
because the link $S^{3}$ of the diagonal $\Delta$ in $\mathring{X}\times \mathring{X}$ is null-homologous in $\mathring{X}\times \mathring{X}-\Delta$. Another approach is to study rational homotopy groups of the embedding calculus Taylor tower, which could be thought of as a generalization of configuration space integrals that works for any $X$. We also note that when $X=n(S^2\times S^2)$ with $n$ large enough, the Hurewicz image of $\alpha^{X}_{\eta}$ in $H_{k}(\BDiff_{\partial}(X);\bQ)$ is zero. (See Lemma \ref{lem: beta-vanishing}). See \cite{Ebert2023} for similar phenomena in higher dimensions. \end{remark}

\begin{remark}\label{rmk: two types of exotic phenomena}
The map $\rho^X$ has a factorization
\[
\cM^{s}(X)\xrightarrow{\rho^{X}_{1}} \cM^{fs}(X)\xrightarrow{\rho^{X}_{2}} \cM^{t}(X).
\]
Here $\cM^{fs}(X)$ is the moduli space of $X$ in the formally smooth category. (See \cite{Knudsen-Kupers} for definition.) 
The maps $\rho^{X}_{1}$ and $\rho^{X}_{2}$  compare the formally smooth category with the smooth category and the topological category respectively. By smoothing theory,  $\rho^{X}_{1}$ is a weak equivalence when $\dim X\neq 4$. However, in dimension $4$, both  $\rho^{X}_{1}$ and $\rho^{X}_{2}$ can fail to be a weak equivalence. As suggested in \cite{Knudsen-Kupers}, it would be fruitful to study them separately. We say Type  I (resp. Type II) exotic phenomena appears if $\rho^{X}_{1}$ (resp. $\rho^{X}_{2}$) is not a weak equivalence. Theorem \ref{thm: main} implies that configuration space integrals can only detect Type II phenomena. Therefore, those exotic disk bundles established by Watanabe are Type II, while many exotic bundles established via gauge theory are Type I (see \cite{KonnoLin2022} for relevant discussions).
\end{remark}


\subsection*{Organization of the paper} In Section \ref{section: road map}, we give an outline of our proof. This 3-page section summarises the main ideas of the paper.
In Section \ref{section: Preliminaries}, we recall preliminaries on microbundles, smoothing theory and graph complexes. In Section \ref{section: FM compactification}, we construct the Fulton-Macpherson compactifications for the configuration spaces of formally smooth manifolds and formally smooth bundles. And we prove that these compactifications are stratified topological manifolds. In Section \ref{section: propagators}, we define the propagator system as a compatible collection of maps from the Fulton-Macpherson compactifications of various configuration spaces to $S^3$. The propagator system is used in Section \ref{section: Kontsevich integral} to define the configuration space integral $I^{\Top}$ for formally smooth bundles. In Section \ref{section: equivalence of definitions}, we prove that $I^{\Top}$ and $I$ are equal for smooth bundles. The main result (Theorem \ref{thm: main}) and its topological applications (Corollary \ref{cor: Top and Homeo}, Theorem \ref{thm: MMM class}, Theorem \ref{thm: Sm noncontractible}) are proved in Section \ref{section: applications}.
\subsection*{
Acknowledgements} The authors would like to thank Xujia Chen, Ruizhi Huang, Manuel Krannich, Mark Powell, Daniel Ruberman and Tadayuki Watanabe for several enlightening discussions. The motivating question of the current paper (relation between configuration space integrals and (formal) smooth structures) was first mentioned to us independently by Michael Freedman and Sander Kupers. We are very grateful to them for their willingness to share their insight.

\section{Sketch of the argument}\label{section: road map}
In this short section, we sketch the proofs of our main results. Let $D^{4}\to E_{D}\xrightarrow{\pi} B$ be a topological bundle,  equipped with a boundary trivialization and a framing $\theta$. We assume the base $B$ is also a closed topological manifold. We attach a product bundle $B\times (S^4-D^{4})$ and form the sphere bundle $S^4\to E^{+}\to B$ with a canonical section at infinity. For $n\geq 1$, we denote $\{0,\cdots, n\}$ by $\mathbf{n}$. Consider the space
\[
E^{n}:=\{(x_{0},\cdots, x_{n})\in E^{+}\times_{B}\cdots \times _{B}E^{+}\mid x_0=\infty\}.
\]
For each $S\subset \mathbf{n}$ with $|S|\geq 2$, we have the diagonal 
\[
\Delta_{S}:=\{(x_0,\cdots, x_n)\in E^n\mid x_{i}=x_{j},\ \forall i,j\in S  \}.\]
The configuration space $C_{\mathbf{n}}(E^+)$ is defined to be $E^{n}-\bigcup_{S\subset \mathbf{n}, |S|\geq 2}\Delta_{S}$. A main step in our proof is a construction of the compactification $\bar{C}^{G}_{\mathbf{n}}(E^+)$, which is a version of the Fulton-Macpherson compactification  for formally smooth bundles. To construct this space, we extend $\theta$ to a formal smooth structure on $E^{+}$, represented by an embedding $\tau:\mathcal{N}(\Delta)\hookrightarrow T^{v}E^{+}$ such that $\tau(x_1,x_2)\in T^{v}_{x_1}E^{+}$. Here $\mathcal{N}(\Delta)$ is a tubular neighborhood of the diagonal $\Delta$ in $E^{+}\times_{B}E^{+}$. And $T^{v}E^{+}$ is a vector bundle over $E^+$, called the vertical tangent bundle. Next, we choose an observer system $G$, i.e.  a collection of maps $G=\{g_{S}: \Delta_{S}\to X\}_{S\subset \mathbf{n}}$, which are  small perturbations of the maps \[f_{S}:\Delta_{S}\to X,\quad (x_0,\cdots, x_n)\mapsto x_{\min S}\]
that satisfies the following compatibility condition:
For any $\Delta_{S}\subset \Delta_{T}$, the following diagram must commute
\begin{equation}\label{eq: g_S compatible}
\begin{tikzcd}
 \cN(\Delta_S)\cap \Delta_T \arrow[r,"g_T"] \arrow[d,"\rho_{S}"] & X \\
 \Delta_S \arrow[ru,"g_S",swap]    &  ~
\end{tikzcd}.
\end{equation}
Here $\cN(\Delta_S)$ is a tubular neighborhood of $\Delta_{S}$. And $\rho_{S}:E^{n}\to \Delta_{S}$ sends $(x_0,\cdots,x_n)$ to $(x'_0,\cdots, x'_{n})$ where $x'_{i}=x_{i}$ if $i\notin S$ and $x'_{i}=x_{\min S}$ if $i\in S$. 
Given any $i\in S-\min S$, we have a commutative diagram 
\begin{equation*}
\begin{tikzcd}
 \cN(\Delta_S) \arrow[r,"h_{i}"] \arrow[d,"\rho_{S}"] & T^{v}E^{+}\ar[d] \\
 \Delta_{S}\arrow[r,"g_{S}"] & E^{+}
\end{tikzcd},
\end{equation*}
with the map $h_{i}$ defined by 
\[
h:\vec{x}=(x_{0},\cdots, x_{n})\mapsto \tau(g_{S}\circ \rho_{S}(\vec{x}),x_{i})-\tau(g_{S}\circ \rho_{S}(\vec{x}),x_{\min S}).
\]
Here we use the linear structure on $T^{v}E$ to subtract two points on the same fiber. The maps $\{h_{i}\}_{i\in S-\min S}$ together give an open embedding 
\[
\mathcal{N}(\Delta_{S})\to \operatorname{\oplus}\limits_{i\in S-\min S} g_{S}^*(T^vE^{+}).
\] 
This gives a linear structure $\theta_S$ on the normal microbundle of $\Delta_{S}$ in $E^n$. By (\ref{eq: g_S compatible}), $\theta_{S}$ and $\theta_{T}$ are compatible when $\Delta_{S}\subset \Delta_{T}$. Following \cite{Axelrod-Singer}, we define the space $\bar{C}^{G}_{\mathbf{n}}(E^{+})$ as the closure of the canonical embedding
\[
\iota: C_\mathbf{n}(E) \to E^n\times \prod_{S\subset\mathbf{n}, |S|\ge 2} \Bl_{\Delta_S}^{\theta(g_S)} E^n.
\]
Here $\Bl_{\Delta_S}^{\theta(g_S)} E^n$ denotes the blow up of $E^n$ along $\Delta_{S}$, using $\theta(g_S)$. The compatibility condition allows us to perform a local coordinate change near each point of $\cup_{S} \Delta_{S}\subset E^{n}$ and locally turn $\iota$ into the embedding of smooth configuration spaces. As a consequence, we can show that $\bar{C}^{G}_{\mathbf{n}}(E^{+})$ satisfies many important properties of the classical Fulton-Macpherson compatification. Most importantly, $\bar{C}^{G}_{\mathbf{n}}(E^{+})$ is a topological manifold stratified by various configuration spaces, associated with quotient sets and subsets of $\mathbf{n}$. In particular, we have a canonical decomposition 
\begin{equation}\label{eq: codim-1 bounary decomposition}
\partial \bar{C}^{G}_{\mathbf{n}}(E^{+})\cong \bigcup\limits_{S\subset \mathbf{n}, |S|\geq 2} \bar{C}^{G_{S}}_{\mathbf{n}/S}(E^{+})\times \bar{C}^{*}_{S}(\mathbb{R}^4).
\end{equation}
of the boundary into the  codimension-1 faces. Note that $G$ is unique up to homotopy. So $\bar{C}^{G}_{\mathbf{n}}(E^{+})$ is well-defined up to concordance. It is worth mentioning that a formal smooth structure is essential in this construction. Actually, Kupers \cite{Kupers2020} proved that there exists no natural definition of Fulton-Macpherson compactification for topological bundles (without formal smooth structures).

The next step is to define the so-called ``propagator'', which measures the relative direction of different coordinates of a point in the configuration space. In the classical setting, such a propagator is defined as a differential form on $\bar{C}_2(E^+)$ and then pulled back via forgetful maps $\bar{C}_{\mathbf{n}}(E^+)\to \bar{C}_{\mathbf{2}}(E^+)$ induced by various inclusions $\mathbf{2}\to \mathbf{n}$. In our setting, there are two essential difficulties: First, an inclusion $\mathbf{2}\to \mathbf{n}$ no longer induces a forgetful map between compactified configuration spaces. Second, $\bar{C}^{G}_{\mathbf{n}}(E^+)$ is just a topological manifold, so differential forms are not defined. To overcome the first difficulty, we define the so-called ``propagator system'', which is a compatible collection of maps
\[
\cP:=\{P^{\cS}_{[i],[j]}: \bar{C}^{G_\cS}_{\mathbf{n}/\cS}(E^{+})\to S^3\mid \cS\subset \cM^{+}, [i]\neq [j]\in \mathbf{n}/\cS-\{[0]\}\}.
\] Here $\cS$ is a collection of disjoint subsets of $\mathbf{n}$, each containing at least two elements. And $\cM^{+}$ denotes the set of all possible $\cS$ (including $\cS=\emptyset$).  The function $P^{\cS}_{[i],[j]}$ is called a propagator in our context. They are constructed inductively on the boundary of various faces of $\bar{C}^{G}_{\mathbf{n}}(E^{+})$ and then extended to their interiors via the boundary decomposition like (\ref{eq: codim-1 bounary decomposition}). A crucial point is that $S^3$ is a rational Eilenberg-Maclane space, so as long as the homological obstructions vanish, the extension problem can be solved  after composing with a nonzero-degree map on 
$g:S^3\to S^3$. To overcome the second difficulty, we use Čech cochain as a replacement for differential forms. 
Notice that any Čech cochain on a topological manifold will automatically vanish if its degree is greater than the dimension of the manifold. This is analogous to the fact that a differential form vanishes if its degree exceeds the dimension. This similarity allows us to integrate a Čech cocyle on a topological manifold with boundary and repeat many arguments in \cite{Wa21}. Note that this approach using Čech cohomology also appears in \cite{Chen23}. 

In general, the action of the symmetry group $\Sigma_{n}$ on $C_{\bfn}(E^{+})$ does not extend to $\bar{C}^{G}_{\mathbf{n}}(E^{+})$. So we need a variation of the classical graph complex $\cG$, which we call the \emph{decorated graph complex} and denote by $\widetilde{\cG}$. A generator of $\widetilde{\cG}$ is given by an admissible graph $\Gamma$ with an orientation $o$ and a decoration $\rho$. Admissibility means no multiple edges and self-loops, and each vertex is at least 3-valent. An orientation means an order on the set of edges $E(\Gamma)$, modulo even permutations. And a decoration means a surjection 
 from  $\{1,\cdots,n\}$ to the set of vertices, such that each $k$-valent vertex has exactly $k-2$ preimages. We prove that $H^{*}(\widetilde{\cG})$ contains $H^{*}(\cG)$ as a direct summand. (It might be possible that they are actually be isomorphic.)

Given $(\Gamma,\rho,o)$, we consider the map 
\[
P(\Gamma,\rho):=\prod_{e\in E(\Gamma)} P^{\cS}_{[i(e)],[j(e)]}  :\bar{C}^{G_\cS}_{\mathbf{n}/\cS}(E^{+})\to X(\Gamma):=\prod_{E(\Gamma)} S^3
\]
Here $i(e), j(e)$ denote the initial and terminal vertex of an edge $e$. We pick a Čech cocycle $\xi$ on $X(\Gamma)$ that is equivariant under the permutation action of $\Sigma_{|E(\Gamma)|}$ and represents the generator of the homology $H^{3|E(\Gamma)|}(X(\Gamma);\mathbb{Q})$ given by $o$. Then we define the configuration integral \[I^{\Top}_{(\Gamma,\rho,o)}(\pi):=\frac{1}{|\operatorname{degree}(g)|^{|E(\Gamma)|}}\langle (P^{\cS}_{[i(e)],[j(e)]})^*(\xi), [\bar{C}^{G_\cS}_{\mathbf{n}/\cS}(E^{+})]\rangle\in \mathbb{Q}. \] 
A class $[\eta]\in H^{*}(\widetilde{G})$ is represented by a linear combination of oriented, decorated graphs, hence we can define $I^{\Top}_{[\eta]}(\pi)$ as the corresponding linear combination. By adapting the argument in \cite{Wa21} to the Čech setting, we prove that $I^{\Top}_{[\eta]}(\pi)$ is well-defined and invariant under any framed cobordism. 

Suppose the bundle $\pi:E_{D}\to B$ is smooth. Then $\bar{C}^{G_\cS}_{\mathbf{n}/\cS}(E^+)$ is homeomorphic to the classical Fulton-Macpherson compactification $\bar{C}_{\mathbf{n}/\cS}(E^+)$. By taking a Čech cycle supported near a generic point $x\in X(\Gamma)$, we show that $I^{\Top}_{(\Gamma,\rho,o)}(\pi)$ equals the local degree of $P(\Gamma,\rho)$ at $P(\Gamma,\rho)^{-1}(x)$. Similarly, we can define $I_{(\Gamma,o)}(-)$ using a differential form supported near $x$ and prove that it also equals this local degree. Hence $I^{\Top}_{[\eta]}(\pi)=I_{[\eta]}(\pi)$ when $\pi$ is smooth. 

Now assume the base $B=S^{k}$. We let $F_{b}$ be the fiber over a chosen base point $b\in B$. Since $\Homeo_{\partial}(D^4)$ is contractible, we have a canonical trivialization $E_{D}$ as a topological bundle, which induces a trivialization of   $\cT^{v} E\cong E\times \mathbb{R}^4$. Under this trivialization, we have a $1-1$ correspondence between the sets
\[
\{\text{framings on $E_{D}\xrightarrow{\pi} B$ that is standard on $\partial E\cup F_{b}$}\}/\text{concordance}
\]
and 
\[ \{\text{maps $(E_{D},\partial E_{D}\cup F_{b})\to (\Top(4),e)$}\}/\text{homotopy}.
\]
The later set can be identified with $\pi_{k+4}(\Top(4))$. Therefore, the evaluation of $I^{\Top}_{-}(-)$ gives a map \[\pi_{k+4}(\Top(4))\otimes H^{*}(\widetilde{\cG})\to \mathbb{Q}.\] By taking the dual, we obtain the map $I^{\Top}: \pi_{k+4}(\Top(4))\otimes \mathbb{Q}\to \mathcal{A}_{k}^{\even}$ in (\ref{eq: I-top to Ak}).
The surjectivity of (\ref{eq: I-top to Ak}) follows from the surjectivity of $I$.
To this end, we have proved Theorem \ref{thm: main}. Corollary \ref{cor: Top and Homeo} and Theorem \ref{thm: MMM class} are proved by further calculations in rational homotopy theory, using the topological group structure on $\Top(4)$ and $\Homeo(S^4)$. Let $X$ be a smooth 4-manifold. We glue the bundle $E^{\eta}_{D}\to S^{k}$ defined (\ref{eq：Watanabe family}) with the product bundle $S^{k}\times (X-\Inte(D^4))\to S^{k}$ and get the smooth bundle 
\[
X\to E^{\eta}_{X}\to S^{k}.
\]
The canonical trivialization  $E^{\eta}_{D}\cong_{\Top} S^{k}\times D^{4}$ extends to a trivialization $E^{\eta}_{X}\cong_{\Top}S^{k}\times X$. Hence the smooth structure $E^{\eta}_{X}\to S^{k}$ induces a (nonstandard) smooth structure on the product bundle 
$X\to S^{k}\times X\to X$, denoted by $s$. Note that a formal smooth structure corresponds to reduction of the structure group from $\Top(4)$ to $O(4)$. By studying the Postnikov tower of $(\Top(4)/O(4))_{\bQ}$, we show that $s$ does not even extend to a formal smooth structure on the product bundle $X\to D^{k+1}\times X\to D^{k+1}$.
This implies Theorem \ref{thm: Sm noncontractible} and Corollary \ref{cor: Diff vs Homeo}.

\section{Preliminaries}\label{section: Preliminaries}
In this section, we recall basic facts of microbundle theory \cite{Milnor64} and graph complexes. We also briefly introduce the smoothing theory for higher dimensional manifolds. In particular,  we will state the classification theorem of smooth structures due to Burghelea-Lashof \cite{Burghelea-Lashof} and Kirby-Siebenmann \cite[Essay V]{Kirby-Siebenmann}. Although we will not use this result in later sections, we will compare it with the 4-dimensional situation.

\subsection{Microbundle theory}
\begin{Definition}
A microbundle is a diagram of topological spaces
$$
B\xrightarrow{s} E \xrightarrow{p} B
$$
such that
\begin{itemize}
\item $p\circ s=\id $;

\item For any $p\in B$, there exists a neighborhood $U$ of $p$, a neighborhood $V$ of $s(U)$, a homeomorphism 
$V\cong U\times \mathbb{R}^n$ so that the following diagram commutes
\begin{equation*}
\begin{tikzcd}
 U\arrow[dr,"s_0"] \arrow[r,"s"]  & V \arrow[d,"\cong"] \arrow[r,"p"] & U \\
   ~  & U\times \mathbb{R}^n  \arrow[ur,"p_2"] & ~
\end{tikzcd}
\end{equation*}
Here $s_0(x)=(x,0)$ and $p_2$ is the projection to the second coordinate.
\end{itemize}
In this microbundle $B, E, s$ are called the base, the total space and the zero section respectively.
\end{Definition}
We may  use the total space $E$ or the map $p:E\to B$ to denote 
a microbundle for simplicity when it does not cause any confusion.

\begin{Definition}
Given two microbundles $p:E\to B$ and $p':E'\to B$, an isomorphism between $E$ and $E$ is an homeomorphism
between neighborhoods of the images of the zero sections of $E$ and $E'$ which is compatible with the maps $p,p'$ and the zero sections.
\end{Definition}

\begin{Definition}
Suppose $X\to E\to B$ is a family of topological manifold $X$. The vertical tangent microbundle $\mathcal{T}(E/B)$ is defined as 
$$
E\to E\times_B E\to E
$$
$$
x\mapsto (x,x),(x,y)\mapsto x
$$
\end{Definition}

\begin{Definition}
Suppose $M\subset X$ is an inclusion of topological manifolds. 
If there exists a neighborhood $U$ of $M$ and a retraction $r: U\to M$ such that the diagram
$$
M\xhookrightarrow{} U \xrightarrow{r} M
$$
constitutes a microbundle, then we call this microbundle a \emph{normal microbundle} of $M$ in $X$.
\end{Definition}

\begin{Example}
Suppose $M$ is a topological manifold. Then tangent microbundle $M\times M\to M$ can also be viewed as the normal microbundle 
of the diagonal $\Delta\subset M\times M$.
\end{Example}

\begin{Definition}
Suppose $E\to B$ is a microbundle. We call a microbundle isomorphism $\theta: E\to V$
a \emph{linear structure} of $E$ where $V \to B$ is a vector bundle. 
\end{Definition}

\begin{remark}\label{remark_Top_redection_O} According to a theorem of Kister and Mazur \cite{Kister64}, any microbundle $E$ over  a locally finite simplicial complex $B$
is equivalent to a $\mathbb{R}^n$ fiber bundle with structure group $\Top(n)$, associated to a $\Top(n)$-principal bundle $P$. Then
a linear structure on $E$ can be viewed as 
a reduction of the structure group from $\Top(n)$ to $GL(n)$. Such reductions one-one correspond to points in the section space 
\[
\Gamma(B,P\times_{\Top(n)}\Top(n)/GL(n))
\]
Since $GL(n)\simeq O(n)$, this space is homotopy equivalent to \[\Gamma(B,P\times_{\Top(n)}\Top(n)/O(n)).
\]
\end{remark}

Suppose $Y$ is a topological manifold and $X\subset Y$ is a submanifold. Moreover we assume $X$ has a normal microbundle
$U$ and $U$  admits a linear structure $\theta: U\to V$. After shrinking $U$ if necessary, we may use $\theta$ to identify 
$U\to X$ is with a convex neighborhood 
the zero section of $V$. Given $p,q\in U-X$, we define $p\sim q$ if $p,q$ lies in the same fiber and
there exists $r\in \mathbb{R}^+$ such that $rp=q$. It is clear that $(U-X)/\sim$ is homeomorphic
to the unit disk bundle of $V$. We define the blowup of $U$ along $X$
as 
$$
\Bl_X^\theta U:= \{(p,l)\in U\times (U-X)/\sim | p\in l~\text{or}~p\in X \}
$$
We have a projection $\pi_U:\Bl_X^\theta Y\to U$. 

Given $(p,l)\in \Bl_X^\theta U$, $p$ determines $l$ uniquely whenever $p\notin X$. Therefore we could identify
$\Bl_X^\theta U-\pi^{-1}(X)$ with $U-X$. Now we define the blowup of $Y$ along $X$ as 
\begin{equation}\label{eq_blow_X_Y}
\Bl_X^\theta Y:=(Y-X) \bigcup_{U-X} \Bl_X U.
\end{equation}
We have an obvious projection map $\pi:\Bl_X^\theta Y\to  Y$ and we call it the \emph{blowup map}.
It is clear
from the definition that the exception divisor $\pi^{-1}(Y)$ can be identified with
$$
\mathbb{P}^+(V):=(V-\text{zero section})/\text{positive scaling}.
$$
Given $y\in Y$, the fiber $\pi^{-1}(y)$ can be identified with 
$$
\mathbb{P}^+(E_y):=(V_y-\{0\})/\text{positive scaling}.
$$

\subsection{Smoothing theory}
Now we give a very brief review of the smoothing theory for higher dimensional manifolds. Let $M$ is a $n$ dimensional topological manifold with (possibly empty) boundary. Let $\mathcal{T}M$ be the topological tangent bundle. Then $\mathcal{T}M$ is isomorphic to a fiber bundle associated to a $\Top(n)$-bundle $\FrT(M)$, called the topological frame bundle. We fix a collar neighborhood $\nu(\partial M)\cong \partial M\times [0,1)$ and a smooth structure $\Sigma_{\partial}$  on $\partial M$. Then $\Sigma_{\partial}$ induces a linear structure on $\mathcal{T}M|_{\partial M}$ so gives a section $s_{\partial}$ of the bundle $(\FrT(M)\otimes \Top(n)/O(n))|_{\partial M}$. We use 
\[
\Gamma_{\partial} (\FrT(M)\otimes \Top(n)/O(n))
\]
to denote the space of sections that equals $s_\partial$ on $\partial M$. 

Given a smooth structure $\Sigma$ on $M$ that equals $\Sigma_{\partial}$ on $\partial M$, we denote the smooth manifold by $M_{\Sigma}$. We  define the space 
of smooth structures $\Sm_{\Sigma}(M)$ over $M$ which are diffeomorphic
to $\Sigma$ as the homotopy fiber of 
$$
\BDiff_{\partial} (M_\Sigma)\to \BHomeo_{\partial}(M).
$$ 
Or equivalently, $\Sm_{\Sigma}(M)$ can be defined as the homotopy quotient 
\[\text{Homeo}_{\partial}(M)//\text{Diff}_{\partial}(M_\Sigma).\]

\begin{Theorem}[Burghela-Lashof \cite{Burghelea-Lashof}, Kirby-Siebenmann \cite{Kirby-Siebenmann}]\label{thm: Kirby-Siebenmann}
Let $M$ be a topological manifold of dimension $n\neq 4$. Let $\Sigma_{\partial}$ be a smooth structure on $\partial M$ and let $S$ be the set of diffeomorphism classes of smooth structures on $M$ that equals $\Sigma_{\partial}$ on $\partial M$. Then  the spaces 
$$
\bigsqcup_{[\Sigma]\in S}\Sm_{\Sigma}(M) 
$$
and 
$$
\Gamma_{\partial}(M,\FrT(M)\times_{\Top(n)} \Top(n)/O(n)).
$$
weakly homotopy equivalent.
\end{Theorem}

Now we take $M=D^n$ with the standard smooth structure. Since $\Homeo_{\partial}(M)\simeq \ast$, we have 
$\Sm(D^n)\simeq \BDiff_{\partial}(D^n)$. In this case, Theorem \ref{thm: Kirby-Siebenmann} gives the Morlet correspondence 
\begin{equation}\label{eq: morlet}
\BDiff_{\partial}(D^n)\simeq \Omega^{n}_0(\Top(n)/O(n))\quad \text{when }n\neq 4.   \end{equation}
Here $\Omega^{n}_0$ denotes the unit component of the loop space. There is a framed version of this weak equivalence. Let $\mathcal{F}$ be the space of trivializations of $TD^n$ that are standard near the boundary. Then $\mathcal{F}$ has a natural action by $\Diff_{\partial}(D^n)$. We let 
\begin{equation}\label{eq: framed moduli space}
\BDiff^{\Fr}_{\partial}(D^n):=\operatorname{EDiff}_{\partial}(D^n)
\times_{\Diff_{\partial}(D^n)}
\mathcal{F}
\end{equation}
By pulling back the universal disk bundle over $\BDiff_{\partial}(D^n)$, we obtain the boundary smooth bundle 
\[
D^n\to E\to \BDiff^{\Fr}_{\partial}(D^n).
\]
which carries a canonical boundary trivialization and a canonical framing $\operatorname{fr}_0$ on the vertical tangent bundle $T^{v}E$. After choosing  fiberwise metrics, $\operatorname{fr}_0$ induces a framing   $\operatorname{fr}_1$ on the vertical tangent microbundle $\mathcal{T}^{v}E$. On the other hand, since $\Homeo_{\partial}(D^n)$ is contractible, 
$E$ has a trivialization 
\[
E\cong D^n\times \BDiff^{\Fr}_{\partial}(D^n).
\]
as a topological bundle. This trivialization induces another framing $\operatorname{fr}_2$ on $\mathcal{T}^{v}E$. By taking the difference of $\operatorname{fr}_1$ and $\operatorname{fr}_2$, we obtain a map 
\[
s:D^n\times \BDiff^{\Fr}_{\partial}(D^n)\cong E\to \Top(n).
\]
The two farmings coincide in a neighborhood of $\partial D^n\times \BDiff^{\Fr}_{\partial}(D^n)$ so $s$ sends this neighborhood to the identity element $e$. By taking the adjunction, we obtain a map 
\begin{equation}\label{eq: scanning map}
m:    \BDiff^{\Fr}_{\partial}(D^n)\to \Omega^{n}\Top(n) 
\end{equation}
called the ``scanning map''. When $n\neq 4$, the weak equivalence (\ref{eq: morlet}) implies that $m$ is a weak equivalence to a union of path components of $\Omega^{n}\Top(n)$. When $n=4$, $m$ is still well-defined but unknown to be homotopy equivalence. As we mentioned in Remark \ref{rmk: two types of exotic phenomena},  a failure of $m$ being a homotopy equivalence would correspond to a Type I exotic phenomenon on $D^4$.  

\subsection{The (decorated) graph complex}
\begin{Definition}
Let $\Gamma$ be a graph. We make the following definitions.

\begin{enumerate}
    \item We say $\Gamma$ is admissible if it has no self-loop and multiple edges. And all vertices are at least 3-valent.
    \item We use  $V(\Gamma)$ to denote the set of vertices  and use $E(\Gamma)$ to denote the set of edges of $\Gamma$.
    \item We let 
    \[
 n(\Gamma):=2|E(\Gamma)|-2|V(\Gamma)|
    \]
     and let 
     \[m(\Gamma):=2|E(\Gamma)|-3|V(\Gamma)|.
     \]
     We call $\frac{n(\Gamma)}{2}$ the degree of $\Gamma$ and call $m(\Gamma)$ the excess of $\Gamma$.
      \item An orientation $o$ on $\Gamma$ is an orientation on $\mathbb{R}^{E(\Gamma)}$. Equivalently, it is a total order on $E(\Gamma)$ modulo even permutations. We use $\bar{o}$ to denote the orientation reversal of $o$.
    \item A decoration on $\Gamma$ is a surjective map 
    \[
    \rho: \{1,\cdots ,n(\Gamma)\}\rightarrow V(\Gamma)
    \]
    such that a $k$-valent vertex has exactly $k-2$ preimages.
    \item Given an edge $e\in E(\Gamma)$, we use $\Gamma/e$ to denote the graph obtained by collapsing the edge $e$. Note that $V(\Gamma/e)$ is a quotient of $V(\Gamma)$. So a decoration $\rho$ on $\Gamma$ induces a decoration $\rho/e$ on $\Gamma/e$. Given an orientation $o$ on $\Gamma$ represented an order on $E(\Gamma)$ such that $e$ is the minimal element. By restricting this order to $E(\Gamma/e)=E(\Gamma)-\{e\}$, we get an induced orientation on $\Gamma/e$, denoted by $o/e$.
\end{enumerate}
\end{Definition}

\begin{Definition}\label{defi: graph homology} (1) For positive integers $p,q$, we let $V^{p,q}$ be the $\mathbb{Q}$-vector space spanned by isomorphism classes of oriented, admissible graphs with $p$ vertices and $q$ edges, modulo the relation  $(\Gamma,o)=-(\Gamma,\bar{o})$.
 We define $V^{p,q}=0$ if $p$ or $q$ is non positive. For an integer $m$ and an even integer $n$, we let $\cG^{m,n}=V^{n-m,\frac{3n}{2}-m}$. In other words, $\cG^{m,n}$ is spanned by admissible, oriented graphs with degree $\frac{n}{2}$ and excess $m$. Then we define the linear map 
$\delta: \cG^{m,n}\to \cG^{m+1,n}$
by 
\[\delta(\Gamma,o):=\sum_{e\in E'(\Gamma)} (\Gamma/e,o/e).\]
Here $E'(\Gamma)$ denotes the set of edge $e$ such that $\Gamma/e$ is admissible. The \textbf{graph cohomology} is defined as
\[
H^{m,n}(\cG):=\frac{\ker(\delta^{m,n}: \cG^{m,n} \to \cG^{m+1,n})}{\operatorname{image}(\delta^{m-1,n}: \cG^{m-1,n} \to \cG^{m,n}))}.
\]
(2) Similarly, we let  $\widetilde{V}^{p,q}$
 be the vector space spanned by isomorphism classes of oriented, decorated, admissible graphs with $p$ vertices and $q$ edges, modulo the relation $(\Gamma,\rho,o)=-(\Gamma,\rho,\bar{o})$. We let $\widetilde{\cG}^{m,n}=\widetilde{V}^{n-m,\frac{3n}{2}-m}$ and define the linear map 
$\widetilde{\delta}: \widetilde{\cG}^{m,n}\to \widetilde{\cG}^{m+1,n}$
by 
\[\widetilde{\delta}(\Gamma,\rho,o):=\sum_{e\in E'(\Gamma)} (\Gamma/e,\rho/e,o/e).\]
The \textbf{decorated graph cohomology} is defined as
\[
H^{m,n}(\widetilde{\cG}):=\frac{\ker(\widetilde{\delta}^{m,n}: \widetilde{\cG}^{m,n} \to \widetilde{\cG}^{m+1,n})}{\operatorname{image}(\widetilde{\delta}^{m-1,n}: \widetilde{\cG}^{m-1,n} \to \widetilde{\cG}^{m,n})}.
\]
\end{Definition}
\begin{Definition}\label{defi: A_k} Given integer $k>0$, we use $\mathcal{A}^{\even}_{k}$ to denote graph homology $H_{0,2k}(\mathcal{G}^*)$. Concretely, $\mathcal{A}^{\even}_{k}$ is the vector space spanned by isomorphism classes of oriented, admissible trivalent graphs $(\Gamma,o)$ with $2k$-vertices, modulo the following two 
relations
\begin{enumerate}[(i)]
    \item $(\Gamma,o)+(\Gamma,\bar{o})=0$
    \item $(\Gamma_1,o_1)+(\Gamma_2,o_2)+(\Gamma_3,o_3)=0$ if $\Gamma_{i}$ satisfies the IHX relation. Here the orientation $o_i$ must be compatible with each other under the natural identification between edges of $\Gamma_i$.
\begin{figure}[h]    
\includegraphics[scale=0.6]{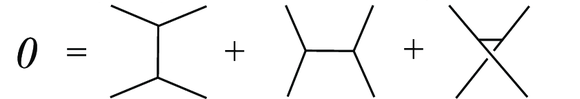}    \centering
\end{figure}
\end{enumerate}
\end{Definition}

We have a map $f:\widetilde{\cG}^{m,n}\to \cG^{m,n}$ defined by $f(\Gamma,\rho,o)=(\Gamma,o)$. It
induces the forgetful map \begin{equation}\label{eq: forgetful}
f_{*}:H^{m,n}(\widetilde{\cG})\to H^{m,n}(\cG).    
\end{equation} One can define a left inverse 
\begin{equation}\label{eq: average}
i_{*}:H^{m,n}(\cG)\to H^{m,n}(\widetilde{\cG})    
\end{equation}
as follows: Given a $(\Gamma,o)$ be an oriented admissible graph, let $A(\Gamma)$ be the set of decorations on $\Gamma$. Then we define $i: \cG^{m,n}\to \widetilde{\cG}^{m,n}$ by 
\[
i(\Gamma,o):=\frac{1}{|A(\Gamma)|}\sum_{\rho\in A(\Gamma)}(\Gamma,\rho,o).
\]
\begin{Lemma}
$i$ is a chain map. And the induced map $i_{*}:  H^{m,n}(\cG)\to H^{m,n}(\widetilde{\cG})$ satisfies $f_{*}\circ i_{*}=\id$. Therefore, $f_{*}$ is a split surjection.   
\end{Lemma}
\begin{proof}
Let $(\Gamma,o)$ be an admissible graph with $p$ vertices and $q$ edges. We denote the vertices of $\Gamma$ by $v_1,\cdots, v_p$ and denote the degree of $v_{i}$ by $k_{i}$. Then we have 
\[
|A(\Gamma)|=\frac{(2q-2p)!}{(k_1-2)!\cdots (k_{p}-2)!}.
\]
Consider an edge $e\in E'(\Gamma)$ that connects the vertices $v_{i},v_{j}$. 
In the graph $\Gamma/e$, these two vertices are merged into a single vertex $v_{ij}$ with degree $k_{i}+k_{j}-2$. Consider the map
\[
A(\Gamma)\to A(\Gamma/e),\quad  \rho\to \rho/e.
\] 
Then each element in $A(\Gamma/e)$ has exactly $\frac{(k_{i}+k_{j}-4)!}{(k_{i}-2)!(k_{j}-2)!}$ preimages in $A(\Gamma)$. Notice that this number also equals $\frac{|A(\Gamma)|}{|A(\Gamma/e)|}$. From this observation, we see that 
\[
\begin{split}
\delta\circ i(\Gamma,o)&=\frac{1}{|A(\Gamma)|}\sum_{\rho\in A(\Gamma)}\sum_{e\in E'(\Gamma)}(\Gamma/e,\rho/e,o/e)\\&=\frac{1}{|A(\Gamma)|}\sum_{e\in E'(\Gamma)}\sum_{\rho\in A(\Gamma/e)}\frac{|A(\Gamma)|}{|A(\Gamma/e)|}(\Gamma/e,\rho/e,o/e)\\
&=\sum_{e\in E(\Gamma')}\sum_{\rho\in A(\Gamma/e)}\frac{1}{|A(\Gamma/e)|}(\Gamma/e,\rho/e,o/e)\\&=\sum_{e\in E(\Gamma')}i(\Gamma/e,o/e)
\\
&=i\circ \delta(\Gamma,o).
\end{split}
\]
This shows that $i$ is a chain map. The rest assertions are straightforward.
\end{proof}


\section{Construction of the compactification}\label{section: FM compactification}

In this section, we define the Fulton-Macpherson compactification of the configuration space of formally smooth manifold and more generally,  a formally smooth bundle. In the smooth setting, such compactification can be obtained by inductively blowing up various diagonals in the Cartesian product. Similar strategy works in the formally smoothing setting. But extra care have to be taken to ensure that the compactification has a nice boundary stratification structure.

Throughout the section, $S,R,T,A,D$ will always denote subsets of at least two elements.

\subsection{The Fulton-Macpherson compactification for a formally smooth manifold}
Suppose $X$ is a closed, formally smooth  $4$-manifold, i.e. a  topological manifold whose tangent microbundle carries a linear structure $\theta$. 
According to the definitions in the previous section,
this means that
the diagonal $\Delta\subset X\times X$ has a neighborhood $\cN(\Delta)$ equipped with
an open embedding $\theta:\cN(\Delta)\to TX$
where $TX$ is a vector bundle over $X$. We call $TX$ the \emph{tangent bundle} of $X$.
Moreover, $\theta$ preserves the zero section and respects the projection map
$$
\rho : X\times X\to X, ~(x,y)\mapsto x.
$$
We also denote $\cN(\Delta)$ by $\mathcal{T}X$ and denote a fiber $\rho|_{\cN(\Delta)}^{-1}(x)$ by $\mathcal{T}_x X$.
From now on we will use $\theta$ to identify $\cN(\Delta)=\mathcal{T}X$ with an open subset of $TX$.
Therefore we also have an inclusion $\mathcal{T}_xX\subset T_xX$ for any $x\in X$.

We call a neighborhood $\cN$ of $\Delta$ \emph{symmetric}
if $(x,y)\in \cN$ implies $(y,x)\in \cN$. 
Let $\phi:\ X\times X\to X\times X$ be the map swapping the two components.
By taking
$\phi(\cN(\Delta))\cap \cN(\Delta)$, we will assume $\cN(\Delta)$ is symmetric.

We fix an integer $n\geq 1$. Given any nonempty $S\subset \{1,\cdots, n\}$ with $i_0=\min S$, we define the (pure) diagonal
$$
\Delta_S:=\{(x_1,\cdots,x_n)\in X^n|x_i=x_j~\text{whenever}~i,j\in S\}
$$
It has a normal microbundle defined by
\begin{equation}\label{eq: rho_S}
\begin{split}
\rho_S :X^n&\to \Delta_S    \\
\rho_S(x_1,\cdots,x_n)&=(x_1',\cdots,x_n')
\end{split}
\end{equation}
where $x'_i=x_i$ if $i\notin S$ and $x'_i=x_{i_0}$ if $i\in S$. Notice that given two subsets $S\subset S'$,
we have
\begin{equation}\label{eq_composition_rho_S}
\rho_{S'}\circ \rho_S=\rho_{S'}
\end{equation}
Consider the map 
\begin{equation}\label{eq: f_S}
   f_S:\Delta_S\to X 
\end{equation}
defined by $f(x_1,\cdots, x_n)=x_{i_0}$. The we have an isomorphism 
$$
X^n\cong \bigoplus_{i\in S-\{i_0\}} f_S^\ast X^2
$$
$$
(x_1,\cdots, x_n) \mapsto (x_{i_0},x_i)_{i\in S-\{i_0\} }.
$$
compatible with the microbundle structures $\rho_S$ and $\rho$. 
In particular, we have
an isomorphism
\begin{equation}\label{eq_f_S_pullback_NDelta_decomposition}
\cN(\Delta_S)\cong \bigoplus_{i\in S-\{i_0\}} f_S^\ast \cN(\Delta)
=\bigoplus_{i\in S-\{i_0\}} f_S^\ast \mathcal{T}X
\end{equation}

An intersection of pure diagonals is called a \emph{mixed diagonal}. Precisely, let $\cS$ be a collection of disjoint subsets of $\{1,\cdots,n\}$ where all the subsets contain at least 2 elements. We also denote the set consisting all such $\cS$ by $\cM$.
Then we define the mixed diagonal as
$$
\Delta_{\cS}=\bigcap_{S\in \cS} \Delta_S.
$$

Suppose $\cS=\{S_1,\cdots,S_k\}$, then a normal microbundle of $\cS$ can be defined by
$$
\rho_\cS:X^n\to \Delta_\cS
$$
$$
\rho_\cS(\vec{x})= \rho_{S_1}\circ \rho_{S_2}\circ \cdots \circ\rho_{S_k} (\vec{x})
$$
It does not matter which order of $S_1,\cdots, S_k$ is used in the definition of $\rho_\cS$ because they all commute with each other. 
Similar to \eqref{eq_f_S_pullback_NDelta_decomposition},
we have an isomorphism
\begin{equation}\label{eq_f_S_pullback_NDelta_decomposition_mixed_diagonal}
\cN(\Delta_\cS)\cong \bigoplus_{S\in\cS}\bigoplus_{i\in S-\{i_0\}} f_S|_{\Delta_{\mathcal S}}^\ast \mathcal{T}X.
\end{equation}

To define the Fulton-Macpherson compactification, we will need to blow up $X^n$ along various pure diagonals $\{\Delta_{S}\}_{S}$. As we explained in Section \ref{section: Preliminaries}, a linear structure on $\cN(\Delta_S)$ is needed. It seems tempting  to use the pull back of $\theta$ via (\ref{eq_f_S_pullback_NDelta_decomposition}). However, linear structures on various diagonals constructed in this way are not compatible with each other. So we will need a modification. 

We start with the following technical lemma.

\begin{Lemma}\label{lemma_good_nbh_Delta} Let $\cN$ be any symmetric neighborhood of $\Delta$ contained in  $\cN(\Delta)$. There exists a symmetric neighborhood $\cN'\subset \cN$ such that for any $x,y,z\in X$ that satisfies 
\[(x,y),(x,z)\in \cN',\] 
we have
\begin{itemize}
\item $(y,z)\in \cN$;

\item For all $t\in[0,1]$, we have $(y,t\cdot z)\in  \cN$ and $(x,t\cdot z)\in \cN$. Here $t\cdot z$ is defined using the
linear structure on $\cN(\Delta)|_y$.
\end{itemize}
\end{Lemma}
\begin{proof}
Since $\rho:\cN\to X$ can be identified with a neighborhood of the zero section of 
a vector bundle $E\to X$, we may equip $E$ with a Riemannian metric on the fibers and assume
the open unit disk bundle $\mathring{D}(E)$ is included in $\cN$. We use 
$\mathring{D}_\epsilon(E)$ to denote the open disk bundle of radius $\epsilon$. Consider the map
$$
F:X^3=X\times X\times X \to X\times X
$$
$$
(x,y,z)\mapsto (y,z)
$$
Since $F^{-1}\cN$ is an open neighborhood of the diagonal $\Delta_{1,2,3}(X^3)$,
the compactness of $X$ implies that when $\epsilon$ is small enough, 
$$
\pi_{12}^{-1}\mathring{D}_\epsilon(E) \cap \pi_{13}^{-1}\mathring{D}_\epsilon(E)\subset 
F^{-1}\cN
$$
where $\pi_{ij}: X^3\to X\times X$ is the projection map onto the $i$th and $j$th components.
Therefore $\mathring{D}_\epsilon(E)$ satisfies the requirement of the first item.

Define 
$$
U=\{(x,y,z,t)\in 
\pi_{12}^{-1}\mathring{D}_\epsilon(E) \cap \pi_{13}^{-1}\mathring{D}_\epsilon(E)\times [0,1]\mid 
(y,tz)\in \cN, (x,tz)\in \cN\}
$$
Notice that $U$ is open and $\Delta_{1,2,3}(X^3)\times [0,1]\subset U$. Therefore 
the compactness of $X$ and $[0,1]$ implies there is an open neighborhood $V$ of $\Delta_{1,2,3}(X^3)$
in $X^3$ such that $V\times [0,1]\subset U$. Again the compactness of $X$ implies that when 
$\epsilon$ is small enough, we have
$$
\pi_{12}^{-1}\mathring{D}_\epsilon(E) \cap \pi_{13}^{-1}\mathring{D}_\epsilon(E)\subset V.
$$
By taking 
$$
\cN'=\phi(\mathring{D}_\epsilon(E))\cap \mathring{D}_\epsilon(E)
$$
all the requirements in the lemma are satisfied.
\end{proof} 

We apply Lemma \ref{lemma_good_nbh_Delta} to $\cN=\cN(\Delta)$ and obtain a smaller neighborhood  $\cN'(\Delta)$.

\begin{Definition}
(1) Given $x,y\in X$, if $(x,y)\in \cN'(\Delta)$, then we say $x$ \emph{observes} $y$.\\
(2) A function $g_S:\Delta_S\to  X$ is called an \emph{observer function} on $\Delta_S$ if $g_S(\vec{x})$ observes $f_S(\vec{x})$ for any $\vec{x}\in \Delta_S$. In this case, we call $g_S({\vec{x}})$ an \emph{observer} for 
$f_S(\vec{x})$.
\end{Definition}

Given an observer function $g_S:\Delta_S\to X$, we shrink 
$\mathcal{T}X$ to obtain an neighborhood $\mathcal{T}'X$ of the zero section so that
$$
\mathcal{T}'_{f_S(\vec{x})}X\subset \mathcal{T}_{g_S(\vec{x})}X ~\text{for all}~\vec{x}\in \Delta_S.
$$

We can define an equivalence of microbundles 
\begin{align}
h(g_S): f_S^\ast\mathcal{T}'X&\to g_S^\ast {T}X \label{eq_def_of_h_map}\\
     (\vec{x},y) &\mapsto (\vec{x}, y-f_S(\vec{x})) \nonumber
\end{align}
where $\vec{x}\in \Delta_S,y\in \mathcal{T}_{f_S(\vec{x})}'X\subset \mathcal{T}_{g_S(\vec{x})}X$ 
and the subtraction 
$y-f_S(\vec{x})$ is defined using the linear structure on  
${T}_{g_S(\vec{x})}X$.  
The diagram
\begin{equation}\label{diagram_normal_bundle_def_component}
\begin{tikzcd}
 f_S^\ast\mathcal{T}'X \arrow[r,"h(g_S)"] \arrow[d] & g_S^\ast TX \arrow[d] \arrow[r,"\tilde{g}_S"]& TX \arrow[d,"\rho"] \\
 \Delta_S \arrow[r,equal]  & \Delta_S \arrow[r,"g_S"]  &  X
\end{tikzcd}
\end{equation}
commutes. 
We could use $ \tilde{g}_S \circ h(g_S)$ 
to define a new linear structure $\theta(g_S)$ on $f_S^\ast \mathcal{T}'X$.
In particular, the linear 
structure on
$$
f_S^\ast \mathcal{T}'X|_{\vec{x}}= \mathcal{T}_{f_S(\vec{x})}'X \subset \mathcal{T}_{g_S(\vec{x})}X
\subset {T}_{g_S(\vec{x})}X
$$
is defined using $T_{g_S(\vec{x})}X$ together with a translation of the origin 
(from $g_S(\vec{x})$ to $f_S(\vec{x})$). We denote this origin-translated new linear space by
$T_{f_S(\vec{x})}^{\theta(g_S(\vec{x}))}X$. We call $\theta(g_{S}(\vec{x}))$ the observer of for this linear space.


We could use $\theta(g_S)$ to equip $\cN(\Delta_S)$ with a linear structure via the isomorphism (\ref{eq_f_S_pullback_NDelta_decomposition}). By abuse of notation, we still denote this linear structure on $\cN(\Delta_S)$
by $\theta(g_S)$. We denote the corresponding vector bundle by $N^{\theta(g_S)}(\Delta_S)$ and call it
the \emph{normal bundle} of $\Delta_S$. According to the above construction, we have
\begin{equation}\label{eq_fiber_of_N_Delta_S}
N^{\theta(g_S)}_{\vec{x}}(\Delta_S)=N^{\theta(g_S)}(\Delta_S)|_{\vec{x}}
=\bigoplus_{i\in S-\{i_0\}}T_{f_S(\vec{x})}^{\theta(g_S(\vec{x}))}X
\end{equation}
for $\vec{x}\in \Delta_S$.

Next, we generalize the above constructions to mixed diagonals. 

\begin{Definition} An observer function $g^{\cS}$ on a mixed diagonal $\Delta_{\cS}$ is a collection of function 
\[\{g_S^{\cS}:\Delta_\cS\to X)\mid S\in \cS\}\]
such that $g_S^\cS(\vec{x})$ observes $f_S(\vec{x})$ for all $\vec{x}\in \Delta_{\cS}$.    
\end{Definition}

\begin{Definition}\label{Definition: compatibility} Consider two general collections of functions 
\[g^{\cS}=\{g_S^{\cS}:\Delta_\cS\to X\mid S\in \cS\}\]
and 
\[g^{\cT}=\{g_T^{\cT}:\Delta_\cT\to X\mid T\in \cT\}\]
defined on mixed diagonals $\Delta_{\cS}$ and $\Delta_{\cT}$ respectively. Suppose $\Delta_{\cS}\subset \Delta_{\cT}$. We say $g^{\cS}$ is compatible with $g^{\cT}$ if the following condition is satisfied: For any $S\in \cS$ and any $T\in \cT$ such that $T\subset S$, there exits a neighborhood $\cN(\Delta_\cS)$ of $\Delta_\cS$
 so that
 the following diagram commutes:
\begin{equation}\label{diagram_compatible_g_S}
\begin{tikzcd}
 \cN(\Delta_\cS)\cap \Delta_\mathcal{T} \arrow[r,"g_T^\mathcal{T}"] \arrow[d,"\rho_{\mathcal S}"] & X \\
 \Delta_\cS \arrow[ru,"g_S^\cS",swap]    &  ~
\end{tikzcd}
\end{equation}
\end{Definition}

Note that $\Delta_\cS\subset \Delta_\mathcal{T}$
if and only if 
for any $T\in \mathcal{T}$, there exits $S\in \cS$ such that $T\subset S$.
\begin{Definition}
An observer system a collection of observer functions 
$$
G=\{g^{\cS}\mid \cS\in\cM\}
$$
such that $g^{\cS}$ is compatible with $g^{\cT}$ whenever $\Delta_\cS\subset \Delta_\mathcal{T}$.
\end{Definition}

\begin{Lemma}\label{lemma_g_S_existence}
There exits an observer system $G$.
\end{Lemma}
\begin{proof}
$G$ is given by a collection of maps 
\[
\{g^{\cS}_{S}: \Delta_{\cS}\to X\mid S\in \cS, \cS\in\cM\}.
\]
We will construct $G$ by induction. 
Suppose 
$$
\cN'(\Delta)=\cN_0(\Delta)\supset \cN_1(\Delta)\supset \cdots 
\supset \cN_{2n-4}(\Delta)
$$
is a sequence of symmetric neighborhood of $\Delta$ so that any adjacent pair 
$(\cN_k(\Delta),\cN_{k+1}(\Delta))$ satisfies the requirements of 
Lemma \ref{lemma_good_nbh_Delta}.

We define 
$g^{\{\{1,\cdots, n\}\}}_{\{1,\cdots, n\}}$ to be $f_{\{1,\cdots, n\}}$. Suppose 
$g_S^{\cS}$ has been defined for  all diagonals of dimension $\le 4(l-1)$ such that
\begin{itemize}
\item $(g_S^\cS(\vec{x}),f_S(\vec{x}))\in \cN_{2n-2l}(\Delta)$ for all $\vec{x}\in \Delta_\cS$;
\item \eqref{diagram_compatible_g_S} holds for these diagonals.
\end{itemize}

Let $\Delta_\mathcal{T}$ be a diagonal of dimension $4l$. We construct a regular neighborhood $\cN(\Delta_{\cS})$ for each  $\Delta_{\cS}\subsetneq \Delta_{\cT}$ as follows: Any such $\Delta_{\cS}$ has dimension $\le 4(l-1)$. By induction, for all $\vec{x}\in \Delta_\cS$, we have
\begin{equation}\label{eq_gs_fs}
(g_S^\cS(\vec{x}),f_S(\vec{x}))\in \cN_{2n-2l}(\Delta).
\end{equation}
Picking a small enough neighborhood $\cN(\Delta_{\cS})$ so that for all $\vec{y}\in \cN(\Delta_\cS)\cap \Delta_\mathcal{T}$ and all $T\in \cT,~S\in \cS$ with $T\subset S$, we have
\begin{equation}\label{eq_fs_rhos_ft}
(f_{{S}}\circ \rho_\cS(\vec{y}), f_T(\vec{y}))\in \cN_{2n-2l}(\Delta)
\end{equation}
Taking $\vec{x}=\rho_\cS(\vec{y})$ in \eqref{eq_gs_fs}, Lemma \ref{lemma_good_nbh_Delta} implies
\begin{equation}\label{eq_gs_rhos_ft}
(g_S^\cS\circ  \rho_\cS(\vec{y}),f_T(\vec{y}))\in  \cN_{2n-2l-1}(\Delta)~
\text{for all}~ 
\vec{y}\in \cN(\Delta_\cS)\cap \Delta_\mathcal{T}.
\end{equation}

By inductively shrinking $\cN(\Delta_{S})$, we also assume that for any two diagonal $\Delta_{\cS}$ and $\Delta_{\mathcal{R}}$ contained in $\Delta_{\cT}$, their normal neighborhoods satisfies
\begin{equation}\label{eq_NS_NR_intersection}
\cN(\Delta_{\cS})\cap \cN(\Delta_{\mathcal{R}})\subset 
\cN(\Delta_{\cS}\cap \Delta_{\mathcal{R}}).
\end{equation}
Let 
$$
U=\Delta_{\mathcal{T}}\cap \bigcup_{ \Delta_{\cS}\subsetneq \Delta_{\cT}} \cN(\Delta_{\cS}).
$$
We define a map
$$
g: \cN(\Delta_\cS)\cap \Delta_\mathcal{T}\to X
$$
using \eqref{diagram_compatible_g_S}, i.e. $g=g_S^\cS\circ  \rho_\cS$. 
The induction assuption and \eqref{eq_NS_NR_intersection} guarantee that all these maps can patch together to give us
a continuous map $g:U\to X$. According to \eqref{eq_gs_rhos_ft}, we have
\begin{equation}\label{eq_U_g_ft}
(g(\vec{y}), f_T(\vec{y}))\in \cN_{2n-2l-1}(\Delta) ~\text{for all}~ \vec{y}\in U
\end{equation}
Let $\cN'(\Delta_\cS)$ be a neighborhood such that 
$\overline{\cN'(\Delta_\cS)}\subset \cN(\Delta_\cS)$. We define
$$
V=\Delta_{\mathcal{T}}\cap \bigcup_{\Delta_{\cS}\subsetneq \Delta_{\cT}} \cN'(\Delta_{\cS}).
$$
Then we have $\bar{V}\subset U$. We pick a cut-off function $h:\Delta_{\mathcal{T}}\to [0,1]$ which
equals $1$ on $\bar{V}$ and equals $0$ on $U^c$. According to \eqref{eq_U_g_ft}, we have
$$
g(\vec{y})\in   \cN_{2n-2l-1}(\Delta)|_{f_T(\vec{y})}~\text{for all}~ \vec{y}\in U
$$
We define
$
g':U\to X
$
by
$$
g'(\vec{y})=h(\vec{y})\cdot g(\vec{y}) \in \cN(\Delta)_{2n-2l-2}|_{f_T(\vec{y})} \subset X
$$
where the scalar product is defined using the linear structure on $\cN(\Delta)_{2n-2l-2}|_{f_T(\vec{y})}$.
The scalar product is well-defined according to Lemma \ref{lemma_good_nbh_Delta} 
(taking $x=y=f_T(\vec{y}), z= g(\vec{y})$).
It is clear from the definition that $g'$ extends to a continuous function on $\Delta_{\mathcal{T}}$
which is equal to $f_T$ on $U^c$. Now we set $g_T^{\mathcal{T}}=g'$. Similarly we could define 
$g_T^{\mathcal{T}}$ for all $T\in \mathcal{T}$ and $\dim \Delta_\mathcal{T}=4l$. In summary we obtain
a collection of maps $\{g_S^\cS\}$
for all diagonals $\Delta_\mathcal{T}$ with dimension $\le 4l$ such that
\begin{itemize}
\item $(g_S^\cS(\vec{x}),f_S(\vec{x}))\in \cN_{2n-2l-2}(\Delta)$ for all $\vec{x}\in \Delta_\cS$;
\item \eqref{diagram_compatible_g_S} holds for these diagonals.
\end{itemize}
Now this induction process will finally give us a collections of maps satisfying the requirements of the lemma.
\end{proof}

\begin{Definition} Let $G_0$ and $G_1$ be two observer systems. A homotopy from $G_0$ to $G_1$ is a continuous family of collections of maps 
$$
G_t=\{g_{S,t}^{\cS}:\Delta_\cS\to X)| S\in\cS, ~\cS\in\cM\},~t\in [0,1]
$$
such that 
\begin{itemize}
    \item $(g_S^\cS(\vec{x}),f_S(\vec{x}))\in \cN(\Delta)$ for all $\vec{x}\in \Delta_\cS$;
    \item For any $\Delta_{\cS}\subset \Delta_{\cT}$ and any $t\in[0,1]$, the collection $\{g^{\cS}_{S,t}\mid S\in \cS\}$ is compatible with the collection $\{g^{\cT}_{T,t}\mid T\in \cT\}$ (see Definition \ref{Definition: compatibility}). 
\end{itemize}
\end{Definition}

\begin{Proposition}\label{prop_homotopy_G0_G1}
Any two observer systems $G_0$ and $G_1$ are homotopic. 
\end{Proposition}
\begin{proof}
By assumption we have
$$
(g_{S,0}^{\cS}(\vec{x}),f_S(\vec{x})),
(g_{S,1}^{\cS}(\vec{x}),f_S(\vec{x}))
 \in \cN'(\Delta) ~\text{for all}~\vec{x}\in \Delta_\cS.
$$
Lemma \ref{lemma_good_nbh_Delta} implies that 
$$
g_{S,1}^{\cS}(\vec{x})\in \cN(\Delta)|_{g_{S,0}^{\cS}(\vec{x})} 
$$
and 
$$
t\cdot g_{S,1}^{\cS}(\vec{x})\in \cN(\Delta)|_{g_{S,0}^{\cS}(\vec{x})}
\cap \cN(\Delta)|_{f_{S}(\vec{x})}, 
t\in [0,1]
$$
where the scalar product is defined using the linear structure on $\cN(\Delta)|_{g_{S,0}^{\cS}(\vec{x})}$.
Therefore we could define
$$
g_{S,t}^{\cS}(\vec{x}):=t\cdot g_{S,1}^{\cS}(\vec{x})
$$
to obtain the required family $G_t$.
\end{proof}

Let $G$ be an observer system. We still denote $g^{\{S\}}_S$ by $g_S$. We define
$$
C_n(X):=\{(x_1,\cdots,x_n)\in X^n|x_i\neq x_i~\forall i,j\}
$$
We have inclusion maps $\iota_0: C_n(X)\to X^n$ and $\iota_S :C_n(X)\to \Bl_{\Delta_S}^{\theta(g_S)} X^n$.
We define
\begin{equation}\label{eq_def_iota}
\iota=(\iota_0, (\iota_S)): C_n(X)\to  X^n\times \prod_{S\subset \{1,\cdots,n\}} \Bl_{\Delta_S}^{\theta(g_S)} X^n
\end{equation}
\begin{Definition}\label{def_FM_cpt}
Given an observer system $G$, we define the Fulton-Macpherson compactification of $C_n(X)$ as
$
\bar{C}_n^G(X):= \overline{\Imm\iota}
$
where $\iota$ is the map in \eqref{eq_def_iota}.
\end{Definition}
\begin{remark}
Even though the mixed diagonal $\Delta_{\cS}$ and map $g^{\cS}_S$ in Lemma \ref{lemma_g_S_existence}
do not appear in the definition of $\bar{C}_n^G(X)$ when $|\mathcal S|\ge 2$, they will be essential
to prove that $\bar{C}_n^G(X)$ is a manifold and has the desired stratification in Proposition \ref{prop_Cn_bar_manifold}. 
\end{remark}

\begin{remark}
If $X$ is a smooth manifold and the blowups in
\eqref{eq_def_iota} are defined using the smooth structure, then Definition \ref{def_FM_cpt} is exactly the definition of 
Fulton-Macpherson conpactification $\bar{C}_n(X)$
for the smooth case used in
\cite{Axelrod-Singer}. If $X$ is a non-compact
smooth manifold, then  $\bar{C}_n(X)$ can still be defined in the same way but $\bar{C}_n(X)$ is not compact any more.
\end{remark}

The boundary of $\bar{C}_n^G(X)$ can be described in a way similar to that in \cite[Section 5]{Axelrod-Singer}.
The exceptional divisor of $\Bl_{\Delta_S}^{\theta(g_S)} X^n$ can be identified with
$\mathbb{P}^+(N(\Delta_S,\theta(g_S))$ and the fiber of the exceptional divisor over a point
$\vec{x}\in \Delta_S$ can be identified with
$$
\mathbb{P}^+(
N^{\theta(g_S)}_{\vec{x}}(\Delta_S))
=\mathbb{P}^+( (T_{f_S(\vec{x})}^{\theta(g_S(\vec{x}))}X)^{S-\min S} )
$$
according to \eqref{eq_fiber_of_N_Delta_S}. 

Suppose $(\vec{x},(y_S)_{S\subset \{1,\cdots,n\}})\in X^n\times \prod_{S} \Bl_{\Delta_S}^{\theta(g_S)} X^n$ is a point in
$\bar{C}_n^G(X)$. Let $\pi_S:\Bl_{\Delta_S}^{\theta(g_S)} X^n\to X^n$ be the blowup map. It is clear
from the definition that $\pi_S(y_S)=\vec{x}$.
If $\vec{x}\notin \Delta_S$, then $y_S$ is away from the exceptional divisor and
determined by $\pi_S(y_S)=\vec{x}$. If $\vec{x}\in \Delta_S$, then $y_S$ lies in the exceptional divisor and the data of $y_S$ consist of  
$$
[u_S]=[(u_i^S)_{i\in S-\min S}]\in \mathbb{P}^+((T_{f_S(\vec{x})}^{\theta(g_S(\vec{x}))}X)^{S-\min S} )
$$
Or equivalently,
$$
[u_S]=[(u_i^S)_{i\in S}]\in ((T_{f_S(\vec{x})}^{\theta(g_S(\vec{x}))}X)^{S} -\text{total diagonal})/\sim
$$
where the equivalence relation ``$\sim$'' is defined by positive scaling and diagonal translations.
Instead of $(\vec{x},(y_S)_{S\subset \{1,\cdots,n\}})$, we will use 
$(\vec{x},\{[u_S]|\vec{x}\in \Delta_S\})$ to represent a point in $\bar{C}_n^G(X)$. 

We will need the following definition. 
\begin{Definition}\label{def_nested_set}
We say a collection $\mathcal{F}$ of subsets of $\{1,\cdots,n\}$ is \emph{nested} if
\begin{itemize}
\item If $A\in \mathscr{F}$, then $|A|\ge 2$.

\item If $R,T\in \mathscr{F}$, then they are either disjoint or one contains the other.
\end{itemize}
\end{Definition}
\begin{Proposition}\label{prop_point_in_Cn_bar}
Suppose $\alpha=(\vec{x},(y_S)_{S\subset \{1,\cdots,n\}})\in X^n\times \prod_{S} \Bl_{\Delta_S}^{\theta(g_S)} X^n$.
 The data of $y_S$ consist of  
$$
[u_S]=[(u_i^S)_{i\in S}]\in ((T_{f_S(\vec{x})}^{\theta(g_S(\vec{x}))}X)^{S} -\text{total diagonal})/\sim
$$
when $y_S$ lies in the exceptional divisor. Then $\alpha\in \bar{C}_n^G(X)$ if and only if
\begin{itemize}
\item[(i)] $\pi_S(y_S)=\vec{x}$;

\item[(ii)] If $\vec{x}\in \Delta_S\subset \Delta_T$ (hence
$T\subset S$) and $u_S|_T=(u_i^S)_{i\in T}$ are not all equal, then 
$$
[(u^T_i)_{i\in T}]=[(u^S_i)_{i\in T}] \in  
((T_{f_T(\vec{x})}^{\theta(g_T(\vec{x}))}X)^{T} -\text{total diagonal})/\sim.
$$
Notice that we have
$g_T(\vec{x})=g_S(\vec{x})$ and $f_T(\vec{x})=f_S(\vec{x})$ according to \eqref{diagram_compatible_g_S}.
\end{itemize}
\end{Proposition}
\begin{proof}
We first prove the ``only if'' part. (i) is clear from the definition. For (ii), we assume
$$
\vec{x}(l)=(x_1(l),\cdots,x_n(l))
$$
 is a sequence in
$C_n(X)$ converging to $\alpha$, $\iota_S(\vec{x}(l))$ converges to $[u_S]$ and 
$\iota_T(\vec{x}(l))$ converges to $[u_T]$.
Suppose $i_0=\min S$ and $j_0=\min T$. 
Since $\vec{x}(l)$ converges to $\vec{x}\in \Delta_S\subset \Delta_T$ in $X^n$,  we have
$$
x_i(l)\in \mathcal{T}_{x_{i_0}(l)}X \subset T_{x_{i_0}(l)}^{\theta(g_S\circ \rho_S(\vec{x}(l)))}X ~\text{for all}~
i\in S-\{i_0\}
$$
and
$$
x_i(l)\in \mathcal{T}_{x_{j_0}(l)}X \subset T_{x_{j_0}(l)}^{\theta(g_T\circ \rho_T(\vec{x}(l)))}X ~\text{for all}~
i\in T-\{j_0\}
$$
when $l$ is large enough. When $l$ is large enough, we also have 
$\rho_T(\vec{x}(l))\in \cN(\Delta_S)\cap \Delta_T$ where $\cN(\Delta_S)$ the neighborhood 
in  \eqref{diagram_compatible_g_S}. Therefore \eqref{diagram_compatible_g_S} and \eqref{eq_composition_rho_S}
imply 
$$
g_T\circ \rho_T(\vec{x}(l))=g_S\circ \rho_S\circ  \rho_T(\vec{x}(l))=
g_S\circ \rho_S(\vec{x}(l)).
$$
Thus
$$
x_i(l)\in \mathcal{T}_{x_{j_0}(l)}X \subset T_{x_{j_0}(l)}^{\theta(g_S\circ \rho_S(\vec{x}(l)))}X ~\text{for all}~
i\in T-\{j_0\}.
$$
Notice that the two linear spaces $T_{x_{i_0}(l)}^{\theta(g_S\circ \rho_S(\vec{x}(l)))}X$ and  $T_{x_{j_0}(l)}^{\theta(g_S\circ \rho_S(\vec{x}(l)))}X$ are defined using the same observer $g_S\circ \rho_S(\vec{x}(l))$
but have origins $x_{i_0}(l)$ and $x_{j_0}(l)$ respectively. Both spaces converge to 
$T_{x_{i_0}}^{\theta(g_S (\vec{x}))}X$ when $l$ tends to infinity. Now we equip $TX$ with a Riemannian metric. Then
$[u_S]$ can be obtained by taking
$
u^S_{i_0}=0
$
and
$$ 
u^S_i=\lim_{l\to \infty} \frac{x_i(l)-x_{i_0}(l)}{(\sum_{k\in S-\{i_0\}}|x_k(l)-x_{i_0}(l))|^2)^{\frac12}}\in \lim_{l\to \infty} T_{x_{i_0}(l)}^{\theta(g_S\circ \rho_S(\vec{x}(l)))}X =T_{x_{i_0}}^{\theta(g_S (\vec{x}))}X
$$
for $i\in S-\{i_0\}$.
Similarly, 
$[u_T]$ can be obtained by taking
$
u^T_{j_0}=0
$
and
$$ 
u^T_i=\lim_{l\to \infty} \frac{x_i(l)-x_{j_0}(l)}{(\sum_{k\in T-\{j_0\}}|x_k(l)-x_{j_0}(l))|^2)^{\frac12}}\in \lim_{l\to \infty} T_{x_{j_0}(l)}^{\theta(g_S\circ \rho_S(\vec{x}(l)))}X =T_{x_{i_0}}^{\theta(g_S (\vec{x}))}X
$$
for $i\in T-\{j_0\}$.  We may scale and re-define
\begin{align*}
u^T_i&=\lim_{l\to \infty} \frac{x_i(l)-x_{j_0}(l)}{(\sum_{k\in S-\{j_0\}}|x_k(l)-x_{i_0}(l))|^2)^{\frac12}}\\
     &=\lim_{l\to \infty} \frac{(x_i(l)-x_{i_0}(l))-(x_{j_0}(l)-x_{i_0}(l))
     }{(\sum_{k\in S-\{j_0\}}|x_k(l)-x_{i_0}(l))|^2)^{\frac12}}\\
     &=u^S_i-u^S_{j_0}
\end{align*}
Since $u_S|_T=(u_i^S)_{i\in T}$ are not all equal, $[(u^S_i-u^S_{j_0})_{i\in T-\min T}]$ is a well-defined 
element in $\mathbb{P}^+((T_{f_T(\vec{x})}^{\theta(g_T(\vec{x}))}X)^{T-\min T} )$. Now it is clear that
$[(u^T_i)_{i\in T}]=[(u^S_i)_{i\in T}]$.

Now we prove the ``if'' part. Let $\Delta_\cS$ be the minimal diagonal containing 
$\vec{x}$. For simplicity we assume $\cS=\{S\}$ and hence $\Delta_\cS=\Delta_S$. The general case is  analogous.
Let $\min S=i_0$.
We  define $\mathscr{C}$  to be the minimal collection  of subsets of $S$  such that
\begin{itemize}
  \item $S\in \mathscr{C}$;
\item $\forall T\subset S (u^T_i=u^T_j \Leftrightarrow \exists R\in \mathscr{C}(R\subsetneqq T,i,j\in R))$.
\end{itemize}
Alternatively, the collection $\mathscr{C}$ can be defined by induction: we call $S$ a set at level 1. 
If entries of $u_S=(u^S_i)$ are all distinct, then we stop and define $\mathscr{C}=\{S\}$. Otherwise
the minimal
diagonal $\Delta_{\mathcal{L}_2}((T_{f_S(\vec{x})}^{\theta(g_S(\vec{x}))}X)^{S})$ containing $u_S$
determines a collection of disjoint subsets of $S$. If  $T\in \mathcal{L}_2$, then we call $T$ a set at level 2. 
Given a set $T$ at level $k$, if entries of $u_T$ are all distinct then we stop. Otherwise
we define a collection of sets at level $k+1$ by taking the index of the
minimal diagonal in  $(T_{f_S(\vec{x})}^{\theta(g_S(\vec{x}))}X)^{T}$  containing $u_T$. Finally we define
$\mathscr{C}$ to be the collection of subsets at all levels. 
It is direct to check that $\mathscr{C}$ is nested.

We take a representative $u_T$ for each $[u_T]$
such that $u^T_{\min T}=0$ and
define a sequence 
$\vec{x}(l)=(x_1(l),\cdots,x_n(l))\in X^n$ by taking
$$
x_i(l)=x_i+\sum_{\substack{T\in \mathscr{C}\\i\in T}} l^{-\text{level}(T)} u^T_i
$$
where the additions and scalar products are computed in 
$\mathcal{T}_{x_{i_0}}X\subset T_{x_{i_0}}^{\theta(g_S(\vec{x}))}X$. Notice that 
$x_i(l)\equiv x_i$ when $i\notin S$ or $i=i_0$.
It is straightforward to check that $\vec{x}(l)\in C_n(X)$ when $l$ is large enough and
$\lim_{l\to \infty} \vec{x}(l)=\vec{x}$. Simlary, we define a sequence 
$\vec{x}(l,k)=(x_1(l,k),\cdots,x_n(l,k))\in X^n$ by
$$
x_i(l,k)=x_i+\sum_{\substack{T\in \mathscr{C},i\in T\\ \text{level}(T)\le k} }l^{-\text{level}(T)} u^T_i
$$
When $k$ is no less than largest level of sets in $\mathscr{C}$, $\vec{x}(l,k)=\vec{x}(l)$. We have
$\rho_S(\vec{x}(l,k))=\vec{x}$ for all $l,k$. When $l$ is large, we have
$\vec{x}(l,k)$ lies in the neighborhood $\cN(\Delta_S)$ in \eqref{diagram_compatible_g_S}
for all $k$ and all $T\subset S$
 (we take $\cS=\{S\}$ and $\mathcal{T}=\{T\}$ in \eqref{diagram_compatible_g_S}).
Then \eqref{diagram_compatible_g_S} and \eqref{eq_composition_rho_S}
imply 
\begin{equation}\label{eq_observer_xlk=x}
g_T\circ \rho_T(\vec{x}(l,k))=g_S\circ \rho_S\circ  \rho_T(\vec{x}(l,k))=
g_S\circ \rho_S(\vec{x}(l,k))=g_S(\vec{x})
\end{equation}
for any $T\subset S$ when $l$ is large. 

Now pick any $T\in \mathscr{C}$ at level $k+1$. Suppose $j_0=\min T$.
It is clear from the definition 
that $\vec{x}(l,k)\in \Delta_T$ and $x_{j_0}(l,k')=x_{j_0}(l,k)$ for all $k' \ge k$. 
The equality \eqref{eq_observer_xlk=x} implies that
linear structures on 
$$
\mathcal{T}_{x_{i_0}}X\subset T_{x_{i_0}}^{\theta(g_S(\vec{x}))}X
$$
and 
$$
\mathcal{T}_{x_{j_0}(l,k)}X\subset T_{x_{j_0}(l,k)}^{\theta(g_T(\vec{x}(l,k)))}X
$$
are defined using the same observer. Therefore for any $i\in R\subset T$, we could also view 
$u^R_i\in T^{\theta(g_S(\vec{x}))}_{x_{i_0}}X$ as a vector in 
 $T_{x_{j_0}(l,k)}^{\theta(g_T(\vec{x}(l,k)))}X$.
For any $i\in T$, we have
\begin{align*}\label{eq_xl=xlk+higher_order_terms}
{x}_i(l)&={x}_i(l,k)+
\sum_{\substack{R\in \mathscr{C},i\in R\\ \text{level}(R)> k} }l^{-\text{level}(R)} u^R_i \\
&={x}_i(l,k)+\sum_{\substack{R\in \mathscr{C}\\i\in R\subset T}}l^{-\text{level}(R)} u^R_i \\
&=x_i(l,k)+l^{-k-1}u^T_i + O(l^{-k-2}) \\
&=x_{j_0}(l,k)+l^{-k-1}u^T_i + O(l^{-k-2})
\end{align*}
where the second equal sign follows from the fact $\mathscr{C}$ is nested. The additions in the above
equalities can also be view as additions in 
$\mathcal{T}_{x_{j_0}(l,k)}X\subset T_{x_{j_0}(l,k)}^{\theta(g_T(\vec{x}(l,k)))}X$. Now it is clear
from the above equalities that 
$$
\lim_{l\to \infty} \iota_T(\vec{x}(l))=u_T=(u^T_i)_{i\in T}.
$$
For any $\Delta_R$ with $R\notin \mathscr{C}$, either $\vec{x}\notin \Delta_R$ or there exists some 
$D\in\mathscr{C}$  such that $R\subset D$ and $(u^D_i)_{i\in R}$ are not all equal.  
In both cases $y_R$ or $[u_R]$
are determined by (i) and (ii) in the proposition. Since $\lim \vec{x}(l)$ satisfies (i)
 and (ii) from the proof of the ``only if'' part, we must have $\lim \iota_R(\vec{x}(l))=[u_R]$.
 In summary we obtain
 $\alpha= \lim_{l\to \infty} \iota(\vec{x}(l))\in \bar{C}_n^G(X)$.
 \end{proof}

We define
\begin{equation*}
\mathring{\partial}_S:=\{(\vec{x},([u_T])_{T\subset \{1,\cdots,n\}})\in \bar{C}_n^G(X)|
\vec{x}\in \mathring{\Delta}_S, u^S_i\neq u^S_j~\text{whenever}~i\neq j\}
\end{equation*}
where 
$$
\mathring{\Delta}_S:=\{(z_1,\cdots,z_n)\in X^n| z_i=z_j~\text{if and only if}~i,j\in S   \}.
$$
According to Proposition \ref{prop_point_in_Cn_bar} we see that a point in $\mathring{\partial}_S$ is completely determined by
$\vec{x}$ and $[u_S]$. So it is clear that $\mathring{\partial}_S$ is a (non-compact) manifold
of dimension $4n-1$.

Given a nested collection $\mathscr{F}$ of subsets of $\{1,\cdots,n\}$, 
we define
\begin{align}
\mathring{\partial}_\mathscr{F}:=&
\{(\vec{x},([u_T]))\in \bar{C}_n^G(X)|
x_i=x_j \Leftrightarrow \exists T\in \mathscr{F}(i,j\in T), \nonumber\\
&\forall S\in \mathscr{F}
(u^S_i=u^S_j \Leftrightarrow \exists R\in \mathscr{F}(R\subsetneqq S\text{ and }i,j\in R)) \label{def_partial_F}
\}
\end{align}
A point in $\mathring{\partial}_\mathscr{F}$ is determined by $\vec{x}$ and $[u_T]$ ($T\in \mathscr{F}$).
Therefore it straightforward to check that $\mathring{\partial}_\mathscr{F}$ is a (non-compact) manifold
of dimension $4n-|\mathscr{F}|$. We define $\partial_S:=\bar{\mathring{\partial}}_S$ and 
${\partial}_\mathscr{F}=\bar{\mathring{\partial}}_\mathscr{F}$. We have
$$
{\partial}_\mathscr{F}=\bigsqcup_{\mathscr{G}\supset \mathscr{F}} {\mathring{\partial}}_\mathscr{G}
$$
In particular, this implies 
$$
\partial_S \cap \partial_R=\partial_{\{S,R\}}
$$
if $\{S,R\}$ is nested.
We also obtain a stratification 
$$
\bar{C}_n^G(X)=C_n(X)\sqcup \bigsqcup_{\mathscr{F}} {\mathring{\partial}}_\mathscr{F}.
$$
\begin{Proposition}\label{prop_Cn_bar_manifold}
The space $\bar{C}_n^G(X)$ is a compact manifold with boundary.
\end{Proposition} 
\begin{proof}

Given any point $\alpha=(x_1,\cdots,x_n)\in X^n$,
 there exists $\cS=\{S_1,\cdots,S_k\}\in \cM$ such that $x_i=x_j$ if and only $i,j\in S_l$
 for some $1\le l \le k$. This implies that $\alpha\in \Delta_S$ if and only if $S\subset S_l$ for some $l$.
To simplify the notation we assume $k=2$. Let $S_0$ be the complement of $S_1\cup S_2$ in $\{1,\cdots,n\}$.
We could identify $X^n$ with $X^{S_1}\times X^{S_2}\times X^{S_0}$. Pick a small neighborhood $W$ of $\alpha$
so that  $TX|_{g_{S_j}(W\cap \Delta_{S_j})}$ ($j=1,2$) is a trivial vector bundle.
A trivialization gives us a continuous map 
$$
d_j:TX|_{g_{S_j}(W\cap \Delta_{S_j})}\to \mathbb{R}^4 ~(j=1,2)
$$
which is a linear isomorphism restricted to each fiber $T_{g_{S_j}(\vec{x})}X$. 
By \eqref{eq_f_S_pullback_NDelta_decomposition} we have
\begin{align*}
e_j:\cN(\Delta_{{S_j}})|_{W\cap \Delta_{{S_j} }}\cong 
&( f_{S_j}|_{W\cap \Delta_{S_j}}^\ast \mathcal{T}X)^{S_j-\{\min S_j\}}.
\end{align*}
Using the maps in \eqref{diagram_normal_bundle_def_component}, we define
$$
r_j:f_{S_j}|_{W\cap \Delta_{S_j}}^\ast \mathcal{T}X \xrightarrow{\tilde{g}_{S_j}\circ h(g_{S_j})} TX|_{g_{S_j}(W\cap \Delta_{S_j})}\xrightarrow{d_j}\mathbb{R}^4
$$
We have an open embedding
$$
(\rho_{S_j}, H_j=r_j^{S_1-\{\min S_1\}}\circ e_j):
\cN(\Delta_{{S_j}})|_{W\cap \Delta_{{S_j} }} \to (W\cap \Delta_{{S_j} }) \times (\mathbb{R}^4)^{S_j-\{\min S_j\}}
$$
which maps $W\cap \Delta_{S_j}$ to $(W\cap\Delta_{S_j})\times \{0\}$ and identifies $N^{\theta(g_{S_{j}})}(\Delta_{S_j})$
with  the obvious projection map together with the standard linear structure on $\mathbb{R}^4$.
We define the map
$$
H=(\rho_{\mathcal S}, H_1,H_2): \cN(\Delta_{\mathcal S})|_{W\cap \Delta_{\mathcal S}}
\to 
(W\cap \Delta_{\mathcal S}) \times (\mathbb{R}^4)^{S_1-\{\min S_1\}} \times (\mathbb{R}^4)^{S_2-\{\min S_2\}}
$$

Let $U,V\subset X$ and $Q\subset X^{S_0}$ be open subsets such that 
$$
\alpha\in U^{S_1}\times V^{S_2}\times Q \subset \cN(\Delta_{\cS})|_{W\cap \Delta_{\mathcal S}}
\subset \cN(\Delta_{{S_1}})|_{W\cap \Delta_{{S_1} }} \cap \cN(\Delta_{{S_2}})|_{W\cap \Delta_{{S_2} }} 
\subset X^{S_1}\times X^{S_2}\times X^{S_0}.
$$
The restriction of $H$ gives us the following open embedding after a reordering of the sets in the Cartesian product:
\begin{align}
F:U^{S_1}\times V^{S_2}\times Q&\to U\times (\mathbb{R}^4)^{S_1-\{\min S_1\}} 
\times V\times (\mathbb{R}^4)^{S_2-\{\min S_2\}} \times Q \label{eq_embed_nbh_of_alpha} \\
\gamma=(y_1,\alpha_1,y_2,\alpha_2,\beta)&\mapsto 
(y_1, H_1(\gamma), 
y_2, H_2(\gamma),\beta) \nonumber
\end{align}
where $y_1\in U,y_2\in V,\alpha_1\in U^{S_1-\{\min S_1\}}$, $\alpha_2\in V^{S_2-\{\min S_2\}}$.
Taking $T=S_1, S_2$ in \eqref{diagram_compatible_g_S}, we obtain
$$
g_{S_j}\circ \rho_{S_j}=g^{\cS}_{S_j}\circ \rho_{\cS}\circ \rho_{S_j}=
g^{\cS}_{S_j}\circ \rho_{\cS}~(j=1,2).
$$
Therefore the definition of $H_j$ implies that it
does not depend on $\alpha_{2-j}$. Hence we could rewrite the map $F$ as
$$
F(y_1,\alpha_1,y_2,\alpha_2,\beta)=
(y_1, H_1(y_1,\alpha_1,y_2,\beta), 
y_2, H_2(y_1,y_2,\alpha_2,\beta),\beta).
$$
For each fixed $y_1,y_2,\beta$, the restriction of $F$ gives an open embedding 
$$
U^{S_1-\{\min S_1\}} 
\times V^{S_2-\{\min S_2\}} \to 
(\mathbb{R}^4)^{S_1-\{\min S_1\}} 
\times (\mathbb{R}^4)^{S_2-\{\min S_2\}}.
$$
Hence $F$ is also an open embedding. Under the map $F$, the image of $\Delta_{S_1}$ is an open subset of 
$$
U\times \{\vec{0}\}\times  V\times (\mathbb{R}^4)^{S_2-\{\min S_2\}} \times Q
\subset U\times (\mathbb{R}^4)^{S_1-\{\min S_1\}}
\times V\times (\mathbb{R}^4)^{S_2-\{\min S_2\}} \times Q
$$
and the image of  $\Delta_{S_2}$ is an open subset of 
$$
U\times (\mathbb{R}^4)^{S_1-\{\min S_1\}}
\times V\times \{\vec{0}\} \times Q\subset
U\times (\mathbb{R}^4)^{S_1-\{\min S_1\}}
\times V\times (\mathbb{R}^4)^{S_2-\{\min S_2\}} \times Q
$$
Moreover, their normal microbundles and normal bundles induced by 
$F$, $\cN(\Delta_{S_j})$ and $N^{\theta(g_{S_j})}(\Delta_{S_j})$
are given by the obvious projection maps and  the standard linear structure on $\mathbb{R}^4$.

If $S\subset S_1$ and $\min S=\min S_1$, then $F(\Delta_S)$ is an open subset of 
$$
U\times \{\vec{0}\}^{S-\min S_1}\times (\mathbb{R}^4)^{S_1-S}\times V\times (\mathbb{R}^4)^{S_2-\{\min S_2\}} \times Q
\subset
U\times (\mathbb{R}^4)^{S_1-\{\min S_1\}}\times V\times (\mathbb{R}^4)^{S_2-\{\min S_2\}} \times Q 
$$ 
Taking $T=S$ and $\cS=\{S_1\}$ in \eqref{diagram_compatible_g_S}, we see that the normal microbundle
and normal bundle
of $F(\Delta_S)$ induced by $F$, $\cN(\Delta_S)$ and  $N^{\theta(g_{S})}(\Delta_{S})$ is given by the obvious projection map
and the standard linear structure on $\mathbb{R}^4$.

If $S\subset S_1-\{\min S_1\}$, then $F(\Delta_S)$ is an open subset of 
$$
U\times \Delta_S ((\mathbb{R}^4)^{S_1-\min S_1})\times V\times (\mathbb{R}^4)^{S_2-\{\min S_2\}} \times Q
\subset
U\times (\mathbb{R}^4)^{S_1-\{\min S_1\}}\times V\times (\mathbb{R}^4)^{S_2-\{\min S_2\}} \times Q 
$$ 
Taking $T=S$ and $\cS=\{S_1\}$ in \eqref{diagram_compatible_g_S}, we see that the normal microbundle
and normal bundle
of $F(\Delta_S)$ induced by $F$, $\cN(\Delta_S)$ and  $N^{\theta(g_{S})}(\Delta_{S})$ is given by a projection map
similar to the map $\rho_S$ defined as before and the standard linear structure on $\mathbb{R}^4$.

Similar description also works for the case $S\subset S_2$. We may also assume $U$, $V$ and $Q$ are homeomorphic
to open subsets of Euclidean spaces. 
In summary, all the diagonals near the point $\alpha$ becomes smooth submanifolds of a Euclidean space after the
continuous map $F$. Moreover, the image of the normal microbundles $\cN(\Delta_S)$ and normal 
bundles $N^{\theta(g_{S})}(\Delta_{S})$
of the diagonals under the map $F$ coincides
with the standard normal bundles for smooth submanifolds. Let $\dom(F)$ be the domain of $F$. Then the open subset
$$
Z=\dom(F)\times \prod_S \Bl_{\Delta_S\cap \dom(F)}^{\theta(g_S)}\dom(F)\subset
X^n\times \prod_{S\subset \{1,\cdots,n\}} \Bl_{\Delta_S}^{\theta(g_S)} X^n
$$
can be identified with the definition in the smooth case in \cite[Section 5]{Axelrod-Singer} and 
after this (local) identification
$
\overline{\Imm \iota} \cap Z
$
also coincides with the definition in the smooth case. Hence $\overline{\Imm \iota} \cap W$
is a manifold with boundary. Since $\alpha$ is arbitrary, we conclude that $\bar{C}_n^G(X)$ is a compact manifold with boundary.
\end{proof}

\subsection{The Fulton-Macpherson compactification for a family}
Now we deal with the family case. Let $B$ be a compact topological manifold (possibly with boundary). Let $D^{4}\to E_{D}\to B$ be a boundary-trivialized topological disk bundle, i.e., a fiber bundle whose structure group lies in $\Homeo_{\partial}(D^4)$.  We attach a trivial bundle $(\mathbb{R}^{4}-D^4)\times B$ to $E_{D}$ and obtain a fiber bundle $\mathbb{R}^{4}\to E\to B$. We let $E^{+}$ be the $S^{4}$-bundle obtained by a fiberwise one-point compactification on $E$. (We  treat $S^{4}$ as the unit sphere in $\bR^5$ and identify it with $\mathbb{R}^{4}\cup \{\infty\}$ via the stereographic projection.) Then $E^{+}\setminus E_{D}$ is a trivial bundle with fiber $S^{4}\setminus D^{4}$. We call $E^{+}\setminus E_{D}$ (or $E\setminus E_{D}$) the ``product region'' and use 
\[
q: E^{+}\setminus E_{D}\rightarrow (\mathbb{R}^{4}\cup \{\infty\})\setminus D^{4}
\]
to denote the projection map to the fiber. 
We also fix a length function 
\[
|\cdot |: E^+\rightarrow [0,\infty]
\]
which equals standard length $|q(x)|$ for any $x\notin E_{D}$ 
 and satisfies $|x|<1\iff x\in \operatorname{Int}(E_{D})$ .

We assume $E_{D}$ has a framing, i.e.,  a trivialization  $\theta$ of the vertical tangent microbundle $\cT^{v}E$ that is standard in a neighborhood of $E\setminus E_{D}$. Concretely, $\theta$ is given by a continuous map 
\begin{equation}\label{eq: difference function}
d: V_0\rightarrow \mathbb{R}^{4}    
\end{equation}
defined on a small neighborhood $V_0$ of the diagonal $\Delta_0$ in $E\times_{B}E$, which satisfies the following conditions:
\begin{itemize}
    \item For any $(x_{1},x_{2})\in V_0$ with $x_{1},x_{2}\notin E_{D}$, one has 
    \begin{equation}\label{eq: d is standard outside ED}
    d(x_{1},x_{2})=q(x_{2})-q(x_1).    
    \end{equation}
    
    \item  For any fixed $y\in E$, the restriction of $d$ gives a homeomorphism
    \[
    \{(x_1,x_{2})\in V_0\mid x_{1}=y\}\xrightarrow{\cong} \mathring{D}^{4}
    \]
    that takes $(y,y)$ to $0$. 

    \item There exits  $M_0>1$ such that if $y\in E, |y|\ge M_0$, then
    \begin{equation}\label{eq: V_0 outside E_D}
    \{x_2|(y,x_{2})\in V_0\}\subset E-E_D    
    \end{equation}
\end{itemize}

The linear structure $\theta$ on $\mathcal{T}^vE=\cT(E^{+}/B)$ extends to a linear structure $\theta^{+}$ on 
$\mathcal{T}^v E^+=\cT(E^{+}/B)$. 
Let $\Delta\subset E^+\times_B E^+$ be the diagonal and $\cN(\Delta)$ is a small neighborhood of $\Delta$.
Then $\theta^+$ is an open embedding $\theta^+: \cN(\Delta)\to T^v E^+$ which is also a microbundle 
equivalence. Here 
$T^v E^{+}$ is a vector bundle over $E^+$ and we call it the \emph{vertical tangent bundle} of $E$.
We also denote $\cN(\Delta)$ by $\mathcal{T}^v E^+$ and denote a fiber $\cN(\Delta)|_x$ by 
$\mathcal{T}^v_x(E^+)$. We use $\theta^+$ to identify $\cN(\Delta)= \mathcal{T}^v E^+$ with an open subset
of $T^v E^+$.
Then we have an inclusion $\mathcal{T}^v_x(E^+)\subset {T}^v_x(E^+)$ for all $x\in X$. We also assume $\cN(\Delta)$
satisfies the following condition: 
\begin{itemize}
\item There exists a large real number $M_1> M_0$ such that 
      \begin{equation}\label{eq: N_delta outside ED}
      \{(x_1,x_{2})\in \cN(\Delta)\mid |x_{1}|\le M_1\}=\{(x_1,x_{2})\in V_0\mid |x_{1}|\le M_1\}.     \end{equation}
      \end{itemize}

This condition implies that when $|y|\le M_1$, the linear structure $\theta^+$ on $\cN(\Delta)|_y$ is given by \eqref{eq: difference function}
and when $|y|> M_0$, the linear structure $\theta^{+}$ on $\cN(\Delta)|_y$ is consistent with the smooth structure 
on $S^4\backslash D^4$.

With a slight abuse of notation, we use $\infty$ to denote the infinity point in any fiber of $E^{+}$.  
For $n\geq 1$, we let 
\[E^{n}:=\{(x_{0},\cdots, x_{n})\in E^{+}\times _{B}\cdots\times_{B}E^{+} \mid x_{0}=\infty\}.\] 
Throughout the paper, we will use  $\mathbf{n}$ to denote the set $\{0,1,\cdots,n\}$.
Given any subset $S\subset \mathbf{n}$ with $|S|\geq 2$ and $\min S=i_0$, we have the diagonal 
\[
\Delta_{S}:=\{(x_{0},\cdots ,x_{n})\in E^{n}\mid x_{i}=x_{j},\ \forall i,j\in S\}.
\]
More generally, given any subset $A\subset \mathbf{n}$ that contains $0$, we let
$$
E^A=\{(x_a)_{a\in A}\in E^{+}\times _{B}\cdots\times_{B}E^{+} \mid 
x_{0}=\infty\}.
$$
For any subset of $S$ of $A$, we define the diagonal 
$$
\Delta^{A}_{S}=\{(x_a)_{a\in A}\in E^A \mid 
x_{i}=x_j, \forall i,j\in S \}.
$$

We have maps
$$
\rho_S : E^n\to \Delta_S, ~f_S :\Delta_S \to E^+
$$
which coincides with the previous definitions on the fiber over any point in $B$ (see (\ref{eq: rho_S}) and (\ref{eq: f_S})). Notice that if $i_0=\min S=0$, then
$f_S(\vec{x})=\infty$ for all $\vec{x}$. Similar to \eqref{eq_f_S_pullback_NDelta_decomposition}, we have
\begin{equation}\label{eq_NDelta_decomposition_family}
\cN(\Delta_S)\cong \bigoplus_{i\in S-\{i_0\}} f_S^{\ast} \cN(\Delta)=\bigoplus_{i\in S-\{i_0\}} 
 f_S^{\ast} \mathcal{T}^v E^+ \subset \bigoplus_{i\in S-\{i_0\}} 
 f_S^{\ast}{T}^v E^+
\end{equation}
Therefore 
$\rho_S :\cN(\Delta_S) \to \Delta_S$ inherits a linear structure $\theta_S$ using the above isomorphism
and the linear structure for $\cN(\Delta)$.

Suppose $0\notin S$ and $g_S:\Delta_S\to E^+$ is a fiber-preserving map (as fiber bundles over $B$) which satisfies
$$
(g_S(\vec{x}),f_S(\vec{x}))\in \cN(\Delta)~\text{for all}~\vec{x}\in \Delta_S,
$$
then we call $g_S$ an \emph{observer function}. Similar to \eqref{eq_def_of_h_map}, 
we can define an equivalence of microbundles 
\begin{align}
h(g_S): f_S^\ast\mathcal{T}^v E^+&\to g_S^\ast {T}^v E^+ \label{eq_def_of_h_map_family}\\
     (\vec{x},y) &\mapsto (\vec{x}, y-f_S(\vec{x})) \nonumber
\end{align}
where $\vec{x}\in \Delta_S,y\in \mathcal{T}_{f_S(\vec{x})}^{v}E^+\subset \mathcal{T}_{g_S(\vec{x})}^vE^+$ 
(we may shrink $\mathcal{T}_{f_S(\vec{x})}^{v}E^+$ when necessary)
and the subtraction 
$y-f_S(\vec{x})$ is defined using the linear structure on  
${T}_{g_S(\vec{x})}^vE^+$.  
We could use $h(g_S)$ to equip $f_S^\ast\mathcal{T}^v E^+$ with a new linear structure $\theta(g_S)$. 
More concretely,
the linear structure on 
$$
f_S^\ast\mathcal{T}^v E^+|_{\vec{x}}=\mathcal{T}^v _{f_S(\vec{x}))} E^+ \subset 
\mathcal{T}_{g_S(\vec{x})}^vE^+ \subset {T}_{g_S(\vec{x})}^vE^+
$$
is defined using ${T}_{g_S(\vec{x})}^vE^+$ with a translation of the origin. We denote 
the origin-translated new linear space (or equivalently, the tangent space of ${T}_{g_S(\vec{x})}^vE^+$
at $f_S(\vec{x})$) by $T^v E^+|_{f_S(\vec{x})}^{\theta(g_S(\vec{x}))}$. Roughly speaking, we use $g_S(\vec{x})$ to observe $f_S(\vec{x})$ and its nearby points. 
We could use $\theta(g_S)$ to equip 
\begin{equation*}
 \cN(\Delta_S) \cong \bigoplus_{i\in S-\{i_0\}}f_S^\ast \mathcal{T}^v(X)
\end{equation*}
with a new linear structure (still denoted by $\theta(g_S)$). In short, this linear structure
is defined by the following equivalence of microbundles:
\begin{equation}\label{eq_linear_on_N_Delta_S}
\cN(\Delta_S)\cong  \bigoplus_{i\in S-\{i_0\}} 
 f_S^{\ast} \mathcal{T}^v E^+ \xrightarrow{h(g_S)^{\oplus S-\{i_0\}} } \bigoplus_{i\in S-\{i_0\}} 
 g_S^{\ast}{T}^v E^+.
\end{equation}
We denote the corresponding vector bundle of this linear structure by $N^{\theta(g_S)}(\Delta_S)$, then we have 
\begin{equation}\label{eq_normal_bundle_of_Delta_S}
\cN(\Delta_S)|_{\vec{x}}\subset N^{\theta(g_S)}(\Delta_S)|_{\vec{x}}=\bigoplus_{i\in S-\{i_0\}} 
T^v E^+|_{f_S(\vec{x})}^{\theta(g_S(\vec{x}))}
\end{equation}

Let $\cS$ be a collection of disjoint subsets of $\mathbf{n}$ where all the subsets contain at least 2 elements. We also denote the set consisting all such $\cS$ by $\cM_{\mathbf{n}}$.
Then we define
$$
\Delta_{\cS}=\bigcap_{S\in \cS} \Delta_S
$$
Suppose $\cS=\{S_1,\cdots,S_k\}$, then a normal microbundle of $\cS$ can be defined by
$$
\rho_\cS:E^n\to \Delta_\cS
$$
$$
\rho_\cS(\vec{x})= \rho_{S_1}\circ \rho_{S_2}\circ \cdots \circ\rho_{S_k} (\vec{x})
$$
Similar to \eqref{eq_NDelta_decomposition_family},
we have an isomorphism
\begin{equation}\label{eq_f_S_pullback_NDelta_decomposition_mixed_diagonal_family}
\cN(\Delta_\cS)\cong \bigoplus_{S\in\cS}\bigoplus_{i\in S-\{i_0\}} f_S|_{\Delta_{\mathcal S}}^\ast \mathcal{T}^vE^+.
\end{equation}

The following lemma is an analogue of Lemma \ref{lemma_g_S_existence} and Proposition \ref{prop_homotopy_G0_G1}
in the family case.
We skip its proof since there is no essential change in the proof.
\begin{Lemma}\label{lemma_g_S_existence_family}
There exit a neighborhood $\cN'(\Delta)\subset \cN(\Delta)$ of $\Delta$ and
 a collection of observer functions 
$$
G=\{g_S^{\cS}:\Delta_\cS\to E^+)| S\in\cS, ~\cS\in\cM\}
$$
 such that
\begin{itemize}
\item[(a)]  If $0\in S$, then $g_S^{\cS}=f_S|_{\Delta_\cS}$;

\item[(b)] $(g_S^\cS(\vec{x}),f_S(\vec{x}))\in \cN'(\Delta)$ for all $\vec{x}\in \Delta_\cS$;

\item[(c)]  If $\Delta_\cS\subset\Delta_{\mathcal{T}}$, $S\in \cS, T\in \mathcal{T}, T\subset S$, 
 then there exits a normal neighborhood $\cN(\Delta_\cS)$ of $\Delta_\cS$
 so that
 the following diagram commutes:
\begin{equation}\label{diagram_compatible_g_S_family}
\begin{tikzcd}
 \cN(\Delta_\cS)\cap \Delta_\mathcal{T} \arrow[r,"g_T^\mathcal{T}"] \arrow[d,"\rho_{\mathcal S}"] & E^+\\
 \Delta_\cS \arrow[ru,"g_S^\cS",swap]    &  ~
\end{tikzcd}
\end{equation}
\end{itemize}

Moreover, given two such collections $G_0$ and $G_1$ satisfying \emph{(a) (b) (c)}, there exists a 
homotopy $G_t$ ($0\le t \le 1$)
such that for any $t\in [0,1]$, $G_t$ satisfies \emph{(a) (b) (c)} after replacing $\cN'(\Delta)$ in (b)
by $\cN(\Delta)$.
\end{Lemma}

\begin{Definition}
We call a collection $G$ of maps satisfying the requirements in Lemma \ref{lemma_g_S_existence_family} 
an \emph{observer system} for $E^n$.
\end{Definition}
Now we assume $G$ is an observer system
and still denote $g_S^{\{S\}}\in G$ by $g_S$. Using $G$, we can 
equip any diagonal $\Delta_S$ with a normal bundle structure $\theta(g_S)$. 
 We define
$$
C_\mathbf{n}(E):=\{(x_0,\cdots,x_n)\in E^n|x_i\neq x_i~\forall i,j\}
$$
Similar to \eqref{eq_def_iota}, we have an obvious map
\begin{equation}\label{eq_def_iota_family}
\iota: C_\mathbf{n}(E) \to E^n\times \prod_{S\subset\mathbf{n}, |S|\ge 2} \Bl_{\Delta_S}^{\theta(g_S)} E^n
\end{equation}
\begin{Definition}\label{Def_FM_cpt_family}
Given an observer system $G$, we define the Fulton-Macpherson compactification of $C_\mathbf{n}(E)$ as
$
\bar{C}_{\mathbf{n}}^G(E^+):= \overline{\Imm\iota}
$
where $\iota$ is the map in \eqref{eq_def_iota_family}.
\end{Definition}
\begin{remark}
In Definition \ref{Def_FM_cpt_family} and all the previous constructions in this subsection, 
we may replace $\mathbf{n}$ by any finite ordered set. For example, any quotient set of $\mathbf{n}$
is still an ordered set whose order is defined by comparing the minimal representatives in 
equivalence classes.
\end{remark}
\begin{Proposition}\label{prop_Cn_bar_manifold_family}
The space $\bar{C}_{\mathbf{n}}^G(E^+)$ is a compact manifold with boundary.
\end{Proposition}
\begin{proof}[Sketch of the proof]
The proof can be carried over almost verbatim from the proof of Proposition \ref{prop_Cn_bar_manifold}. Given a point $\alpha\in E^n$ lying over a point $b\in B$, we could choose a small neighborhood  
$N(b)$ of $b$ so that $E_D|_{N(b)}\cong N(b)\times D^4$. Then we have
$$
E^n|_{N(b)}\cong N(b)\times (\{\infty\}\times X^n )
$$
Similar to \eqref{eq_embed_nbh_of_alpha}, we could construct an open embedding of a neighborhood $W$ 
of $\alpha$
into $N(b)\times \mathbb{R}^m$. Moreover, this embedding maps any diagonal containing $\alpha$ into linear 
subspace of $N(b)\times \mathbb{R}^m$ whose normal bundle coincide with its normal bundle as a smooth submanifold. Therefore $\pi^{-1}(W)$ coincides with the definition in the smooth case and  
hence is a manifold where 
$\pi: \bar{C}_{\mathbf{n}}^G(E^+)\to E^n$ is the projection map. Since $\alpha$ is arbitrary, we conclude that
$\bar{C}_{\mathbf{n}}^G(E^+)$ is a compact manifold with boundary. 
\end{proof}

Proposition \ref{prop_point_in_Cn_bar}, as well as
the definitions and discussions below it, can all be adapted to 
$\bar{C}_{\mathbf{n}}^G(E^+)$. For example, a point in $\bar{C}_{\mathbf{n}}^G(E^+)$ is determined by 
$$
(\vec{x},([u_S])_{S\subset \mathbf{n}, |S|\geq 2, \vec{x}\in \Delta_S})
$$
where $\vec{x}\in E^n$ and \[[u_S]=[(u^S_i)_{i\in S}]\in 
((T^vE^+|_{f_S(\vec{x})}^{\theta(g_S(\vec{x})} )^{S} -\text{total diagonal})/\sim.\]
Here $\sim$ again means translation and rescaling.
We have
\begin{equation}\label{eq_partial_S_inner_family}
\mathring{\partial}_S=\{
(\vec{x},([u_T]) \in \bar{C}_{\mathbf{n}}^G(E^+)|\vec{x}\in \mathring{\Delta}_S,
[u_S]\in 
C^\ast_S( T^vE^+|_{f_S(\vec{x})}^{\theta(g_S(\vec{x}))}
\}
\end{equation}
where we define
\begin{align*}
\mathring{\Delta}_S &=\{\vec{x}\in E^n| x_i=x_j\Leftrightarrow i,j\in S\}\\
C_S^\ast(V)&=\{(u_i)_{i\in S}\in V^S| u_i\neq u_j~\text{whenever}~i\neq j\}/ \sim
\end{align*}
for any vector space $V$. We define
the Fulton-Macpherson compactification
$\bar{C}_S^\ast(V)$ as  $\bar{C}_S(V)$
modulo diagonal translations and positive scalings
(see Definition 4.11 and Lemma 4.12 in \cite{Sinha} for more details). Notice that $\bar{C}_S^\ast(V)$
is compact even though $\bar{C}_S(V)$ is non-compact. More generally, given a (continuous)
vector bundle $W\to B$, we could define
$\bar{C}_S^\ast(W)$ by taking the compactification
$\bar{C}_S^\ast(W_b)$ ($b\in B$)
fiberwise.

It is clear that $\mathring{\partial}_S$ is a fiber bundle over $\mathring{\Delta}_S$ 
with fiber $C^\ast_S( T^vE^+|_{f_S(\vec{x})}^{\theta(g_S(\vec{x}))})$. We define
$\partial_S$ as the closure of $\mathring{\partial}_S$ in $\bar{C}_{\mathbf{n}}^G(E^+)$. Given a nested collection $\mathscr{F}$  of subsets of $\mathbf{n}$ (see Definition \ref{def_nested_set}), we define
$\mathring{\partial}_{\mathscr{F}}$ as in \eqref{def_partial_F} let
$\partial_{\mathscr{F}}$ be its closure in $\bar{C}_{\mathbf{n}}^G(E^+)$.


If $T\not\subset S$, then  $\mathring{\Delta}_S\subset E^n-\Delta_T$ and we have an inclusion map
$$
\iota_{S,T}:\mathring{\Delta}_S \to \Bl_{\Delta_T}^{\theta(g_T)} E^n.
$$
The above inclusions give us an embedding map
\begin{equation*}
\iota_S:\mathring{\Delta}_S\hookrightarrow E^n\times \prod_{T\not\subset S} \Bl_{\Delta_T}^{\theta(g_T)} E^n
\end{equation*}
We have a commutative diagram
\begin{equation}\label{diagram_partial_S_inner_Delta_S_inner}
\begin{tikzcd}
 \mathring{\partial}_S \ar[d] \arrow[r,hookrightarrow] &  E^n\times \prod_{T} \Bl_{\Delta_T}^{\theta(g_T)} E^n
 \ar[d] \\
  \mathring{\Delta}_S \arrow[r,hookrightarrow,"\iota_S"] &  E^n\times \prod_{T\not\subset S} \Bl_{\Delta_T}^{\theta(g_T)} E^n
\end{tikzcd}
\end{equation}

Given any diagonal $\Delta_S$, we have $\Delta_{\mathcal{T}}\subset \Delta_S$ if and only if there exists
$T\in \mathcal T$ such that $S\subset T$. We have an 1-1 correspondence 
\begin{equation}\label{eq_quotient_S_correspondence}
\cdot/S:\{\mathcal{T}\in \cM_{\mathbf{n}}|\Delta_{\mathcal{T}}\subsetneq \Delta_S\}\cong \cM_{\mathbf{n}/S}
\end{equation}
$$
\mathcal{T} \mapsto \{T/S|T\in \mathcal{T}, T\not\subset S\}
$$
where  $T/S$ is the image of $T$ under the quotient map $\mathbf{n}\to \mathbf{n}/S$. Notice that we also have  obvious identifications
\begin{equation}\label{eq_diagonal_n_and_ES}
\Delta_S\cong E^{\mathbf{n}/S}=\{(x_a)_{a\in \mathbf{n}/S} \in E\times_B \cdots \times_B\times E| 
x_0=\infty  \}    
\end{equation}
where we still use $0$ to represent the equivalence class of $0$. Moreover, for any $\Delta_\mathcal{T}\subset \Delta_S$, (\ref{eq_diagonal_n_and_ES}) restricts to the identification
\begin{equation}\label{eq_diagonal_n_and_S}
\Delta_\mathcal{T}\cong \Delta^{\mathbf{n}/S}_{\mathcal{T}/S}.
\end{equation}
Therefore \eqref{eq_quotient_S_correspondence} and \eqref{eq_diagonal_n_and_S}
give us an observer system $G_S$ for $E^{\mathbf{n}/S}$ by taking the subset
\begin{equation*}
\{g_T^\mathcal{T}\in G|  \Delta_{\mathcal{T}}\subset \Delta_S     \} \subset G
\end{equation*}

\begin{Proposition}\label{prop_closure_Delta_S_Cn/S}
We have a homeomorphism
$$
\overline{\Imm\iota_S}\cong \bar{C}_{\mathbf{n}/S}^{G_S}(E^+)
$$
\end{Proposition}
\begin{proof}
The closure of $\Imm\iota_{S,T}$ in $\Bl_{\Delta_T}^{\theta(g_T)} E^n$ is just the proper transform
$\tilde{\Delta}_S(T)$ of $\Delta_S$ under the blowup, i.e. the closure of $\Delta_S-\Delta_T$ in
$\Bl_{\Delta_T}^{\theta(g_T)} E^n$. 
By \eqref{eq_NDelta_decomposition_family}, we have an identification
$$
\cN(\Delta_T)= \bigoplus_{i\in T-\min T} 
 f_T^{\ast} \mathcal{T}^v E^+ .
$$

If $S\cap T=\emptyset$, then using the above identification we have
$$
(\Delta_S-\Delta_T)\cap \cN(\Delta_T)=\Delta_S\cap \cN(\Delta_T)-\Delta_T
=\bigoplus_{i\in T-\min T} 
 f_T^{\ast} \mathcal{T}^v E^+|_{\Delta_S\cap \Delta_T}-\text{zero section}
$$
where the microbundle $\bigoplus_{i\in T-\min T} 
 f_S^{\ast} \mathcal{T}^v E^+|_{\Delta_S\cap \Delta_T}$ can also be viewed as the normal microbundle
 of $\Delta_{\{S,T\}}=\Delta_S\cap \Delta_T$ in $\Delta_S$
 (cf. \eqref{eq_f_S_pullback_NDelta_decomposition_mixed_diagonal_family}).
Thus the intersection of $\tilde{\Delta}_S(T)$ and the exceptional divisor in $\Bl_{\Delta_T}^{\theta(g_T)} E^n$ 
is 
$$
\tilde{\Delta}_S(T)\cap \mathbb{P}^+ (N^{\theta(g_T)}(\Delta_T))
=\mathbb{P}^+ (N^{\theta(g_T)}(\Delta_T)|_{\Delta_S\cap \Delta_T})
$$
Therefore we have the identification 
$$
\tilde{\Delta}_S(T)\cong \Bl_{\Delta_{\{S,T\}}}^{\theta(g_T|_{\Delta_{\{S,T\}}})} \Delta_S
= \Bl_{\Delta_{\{S,T\}}}^{\theta(g_T^{\{S,T\}})} \Delta_S
$$
where the second equality follows from the compatibility condition \eqref{diagram_compatible_g_S_family}.

If $S\cap T\neq \emptyset$, then $\Delta_S\cap \Delta_T=\Delta_{S\cup T}$. 
The intersection of $\tilde{\Delta}_S(T)$ and the exceptional divisor in $\Bl_{\Delta_T}^{\theta(g_T)} E^n$ 
is a fiber bundle over $\Delta_{S\cup T}$ whose fiber over $\vec{x}\in \Delta_{S\cup T}$ is
\begin{align}
&\mathbb{P}^+ \{ (u_i)_{i\in T}\in (T^v E^+|_{f_T(\vec{x})}^{\theta(g_T(\vec{x}))})^{T}|
u_{\min T}=0,u_i=u_j~\text{if}~i,j\in S\cap T \} \nonumber \\
&=\{ (u_i)_{i\in T}\in (T^v E^+|_{f_T(\vec{x})}^{\theta(g_T(\vec{x}))})^{T}|
u_i=u_j~\text{if}~i,j\in S\cap T \} /\sim      \label{eq_AAAA}
\end{align}
On the other hand, the exceptional divisor of $\Bl_{\Delta_{S\cup T}}^{\theta(g_{S\cup T})}\Delta_S$
is a fiber bundle over $\Delta_{S\cup T}$ whose fiber over $\vec{x}\in \Delta_{S\cup T}$ is
\begin{align}
&\mathbb{P}^+ \{ (u_i)_{i\in S\cup T}\in (T^v E^+|_{f_{S\cup T}(\vec{x})}^{\theta(g_{S\cup T}(\vec{x}))})^{S\cup T}|
u_{\min {S\cup T}}=0,u_i=u_i~\text{if}~i,j\in S \} \nonumber 
\\
&=\{ (u_i)_{i\in S\cup T}\in (T^v E^+|_{f_{S\cup T}(\vec{x})}^{\theta(g_{S\cup T}(\vec{x}))})^{S\cup T}|
u_i=u_j~\text{if}~i,j\in S \}/\sim             \label{eq_BBBB}
\end{align}
Since $f_{S\cup T}(\vec{x})=f_{T}(\vec{x})$ and  $g_{S\cup T}(\vec{x})=g_{T}(\vec{x})$ according to 
\eqref{diagram_compatible_g_S_family}, it is clear that \eqref{eq_AAAA} and \eqref{eq_BBBB} can be canonically identified with each other. Therefore we have an identification
$$
\tilde{\Delta}_S(T)
\cong \Bl_{\Delta_{S\cup T}}^{\theta(g_{S\cup T})} \Delta_S.
$$

Now $\iota_S$ factors through
$$
\Delta_S\times \prod_{T\not\subset S} \tilde{\Delta}_S(T) \subset 
 E^n\times \prod_{T\not\subset S} \Bl_{\Delta_T}^{\theta(g_T)} E^n.
$$

After identifying $\Delta_S$ with $E^{\mathbf{n}/S}$, $\Delta_T\cap \Delta_S$ becomes $\Delta_{T/S}^{\mathbf{n}/S}$
and $\mathring{\Delta}_S$ becomes $C_{\mathbf{n}/S}(E)$. Moreover, the map $g_{T/S}\in G_S$ coincides with
$g_T|_{\Delta_T\cap \Delta_S}$ when $S\cap T=\emptyset$ (or $g_{S\cup T}$ when $S\cap T\neq \emptyset$)
according to \eqref{diagram_compatible_g_S_family} and
the definition of $G_S$. Therefore $\tilde{\Delta}_S(T)$ can be identified with
$$
\Bl_{\Delta_{T/S}^{\mathbf{n}/S}}^{\theta(g_{T/S})} E^{\mathbf{n}/S}.
$$
After these identifications, $\iota_S$ becomes
$$
C_{\mathbf{n}/S}(E)\hookrightarrow E^{\mathbf{n}/S}\times 
\prod_{\substack{R\subset \mathbf{n}/S\\ T\subsetneq S, R=T/S }} \Bl_{\Delta_{R}^{\mathbf{n}/S}}^{\theta(g_{R})} E^{\mathbf{n}/S}.
$$
Notice that each $R\subset \mathbf{n}/S$ may appear more than once in the above Cartesian product. Fix $R$, the maps
$\iota_{S,T}$ are the same for all $T$ such that $T/S=R$. Thus by taking the diagonal of the product of repeated components
$\Bl_{\Delta_{R}^{\mathbf{n}/S}}^{\theta(g_{R})} E^{\mathbf{n}/S}$ in the above Cartesian product, we may view $\iota_S$
as the map
$$
C_{\mathbf{n}/S}(E)\hookrightarrow E^{\mathbf{n}/S}\times 
\prod_{{R\subset \mathbf{n}/S }} \Bl_{\Delta_{R}^{\mathbf{n}/S}}^{\theta(g_{R})} E^{\mathbf{n}/S}.
$$
Therefore  $\overline{\Imm\iota_S}\cong \bar{C}_{\mathbf{n}/S}^{G_S}(E^+)$ according to Definition 
\ref{Def_FM_cpt_family}.
\end{proof}
Now we see by taking the closure of $\mathring{\partial}_S$ and $\mathring{\Delta}_S$ in 
\eqref{diagram_partial_S_inner_Delta_S_inner}, we obtain a continuous map
$\partial_S\to \bar{C}_{\mathbf{n}/S}^{G_S}(E^+)$.

\begin{Proposition}\label{prop_partial_A_pullback}
The map $\partial_S\to \bar{C}_{\mathbf{n}/S}^{G_S}(E^+)$ is a fiber bundle with fiber type
$\bar{C}^\ast_S(\mathbb{R}^4)$. Indeed we have the following pullback diagram:
\begin{equation}\label{diagram_partial_A_pullback}
\begin{tikzcd}
   {\partial}_S \arrow[r]  \arrow[d]& \bar{C}^\ast_S(T^v E^+) \arrow[d]\\
 \bar{C}_{\mathbf{n}/S}^{G_S}(E^+) \arrow[r,"\bar{g}_S"]     &  E^+   
\end{tikzcd}
\end{equation}
where $\bar{g}_S$ is the following composition of maps:
$$
\begin{tikzcd}
\bar{C}_{\mathbf{n}/S}^{G_S}(E^+)\to E^{\mathbf{n}/S}\cong \Delta_S \arrow[r,"g_S"] &E^+.
\end{tikzcd}
$$
\end{Proposition}
\begin{proof}
According to \eqref{eq_partial_S_inner_family} and \eqref{eq_linear_on_N_Delta_S}, we have a pullback diagram:
\begin{equation*}
\begin{tikzcd}
   \mathring{\partial}_S \arrow[r]  \arrow[d]& C^\ast_S(T^v E^+) \arrow[d]\\
 \mathring{\Delta}_S \arrow[r,"g_S"]     &  E^+   \\
\end{tikzcd}
\end{equation*}
We define a space $Y$ using the pullback diagram
\begin{equation*}
\begin{tikzcd}
   Y\arrow[r]  \arrow[d]& \bar{C}^\ast_S(T^v E^+) \arrow[d]\\
 \bar{C}_{\mathbf{n}/S}^{G_S}(E^+) \arrow[r,"\bar{g}_S"]     &  E^+   
\end{tikzcd}
\end{equation*}
Then $Y$ is a compact manifold with boundary whose interior can be canonically identified with 
$\mathring{\partial}_S$.

According to the proofs of Propositions \ref{prop_Cn_bar_manifold} and 
Proposition \ref{prop_Cn_bar_manifold_family}, any point $\alpha\in \Delta_S$ has neighborhood $W\subset E^n$
such that $W$ is homeomorphic to an Euclidean space and all the diagonals in $W$ is mapped
to standard diagonals of the Euclidean space with standard normal bundle structures as smooth submanifolds. 
Since it is known
that the pullback diagram in the lemma holds in the smooth case \cite[Theorem 4.15]{Sinha}, we
have a homeomorphism $Y\cap p'^{-1}W\to \partial_S\cap p^{-1}W$ which locally extends the identification 
$\mathring{Y} \cong \mathring{\partial}_S$, where $p: \bar{C}_{\mathbf{n}}^G(E^+)\to E^n$ is the projection
map and $p'$ is the composition of maps
$$
Y\to \bar{C}_{\mathbf{n}/S}^{G_S}(E^+)\to E^{\mathbf{n}/S}\cong \Delta_S.
$$ 
These local extensions constructed from 
different $\alpha$ are compatible with each other since a continuous extension of a map on a dense subset to a Hausdorff space is unique. Therefore we can patch together these local extensions to obtain a homeomorphism from $Y$ to $\partial_A$.
\end{proof}

We have the following immediate corollary. 
\begin{Proposition}\label{prop_partial_S_0_product_decomposition}
If $0\in S$, then we have
$$
{\partial}_S\cong \bar{C}_{\mathbf{n}/S}^{G_S}(E^+)\times \bar{C}^\ast_S(T_\infty S^4).
$$
\end{Proposition}
\begin{proof}
If $0\in S$, then Condition (a) in Lemma \ref{lemma_g_S_existence_family} implies 
that $f_S(\vec{x})=g_S(\vec{x})=\infty$ for all $\vec{x}\in \Delta_S$. 
Therefore the image of $\bar{g}_S$ in \eqref{diagram_partial_A_pullback} is also the infinity section.
Thus 
$$
\bar{g}_S^\ast(T^vE^+)\cong  \bar{C}_{\mathbf{n}/S}^{G_S}(E^+) \times T_\infty S^4
$$
and Proposition \ref{prop_partial_A_pullback} implies 
$$
\partial_S\cong\bar{g}_S^\ast \bar{C}^\ast_S(T^vE^+)\cong \bar{C}_{\mathbf{n}/S}^{G_S}(E^+)
\times \bar{C}^\ast_S(T_\infty S^4).
$$
$$
\mathring{\partial}_S\cong \mathring{\Delta}_S\times C^\ast_S(T_\infty S^4). 
$$
\end{proof}

Suppose $0\notin S$. We define a continuous map $g_S':\Delta_S\to E^+$ by
\[
g_S'(\vec{x})=\begin{cases} g_S(\vec{x}) \quad &\text{if }|f_S(\vec{x})|\le M_0.\\
   \frac{|f_S(\vec{x})|-M_0}{M_1-M_0}   f_S(\vec{x}) +  \frac{M_1-|f_S(\vec{x})|}{M_1-M_0} g_S(\vec{x})                   \quad &\text{if } M_0<|f_S(\vec{x})| < M_1.\\
f_S(\vec{x})\quad &\text{if }|f_S(\vec{x})|\ge M_1.
\end{cases}
\]
where the linear combination of $f_S(\vec{x})$ and $g_S(\vec{x})$ is defined using the linear structure on
$\cN(\Delta)|_{f_S(\vec{x})}$. Here $M_0, M_1$ are the constants in (\ref{eq: V_0 outside E_D}) and (\ref{eq: N_delta outside ED}).

We claim that there is a canonical homeomorphism from
$\Bl_{\Delta_S}^{\theta(g_S)}E^n$ to $\Bl_{\Delta_S}^{\theta(g'_S)}E^n$ which extends the identity map
on the complements of the exceptional divisors. This allows us to use $g_S'$ instead of $g_S$ to define
$\bar{C}_{\mathbf{n}}^G(E^+)$. To see that claim, we note that $g_S'$ differs from $g_S$ only in the product region
$E^+-E_D$. These values do not really matter since blowups related to these values coincide with 
 blowups of smooth manifolds. Therefore we could replace the map $\bar{g}_S$ 
 in Proposition \ref{prop_partial_A_pullback} by $\bar{g}_S'$ defined by
 \begin{equation}\label{eq_def_g_S'}
\begin{tikzcd}
\bar{g}_S':\bar{C}_{\mathbf{n}/S}^{G_S}(E^+)\to E^{\mathbf{n}/S}\cong \Delta_S \arrow[r,"g_S'"] &E^+.
\end{tikzcd}
\end{equation}
A similar map $\bar{f}_S$ can be defined by replacing $g_S'$ on the last arrow by $f_S$.

We have a pullback diagram
\begin{equation}
\begin{tikzcd}
   \mathring{\partial}_S \arrow[r]  \arrow[d]& C^\ast_S(T^v E) \arrow[d]\\
 \mathring{\Delta}_S \arrow[r,"g_S'"]     &  E   \\
\end{tikzcd}
\end{equation}
We also have a trivialization 
\begin{equation}\label{eq: trivialization of T^vE}
T^vE\cong E\times \mathbb{R}^4,   \end{equation} 
that coincide with the standard trivialization of $T\mathbb{R}^{4}$ outside $E_{D}$.
From this, we obtain
\begin{equation}\label{eq: trivialization of g_S_TE}
g_S'^\ast T^vE\cong  \mathring{\Delta}_S\times \mathbb{R}^4,~
\mathring{\partial}_S\cong \mathring{\Delta}_S\times C^\ast_S(\mathbb{R}^4).
\end{equation}

\begin{Proposition}\label{prop_trivialization_pullback_TvE+}
There is a trivialization 
$$
\bar{g}_S'^\ast T^vE^+\cong \bar{C}_{\mathbf{n}/S}^{G_S}(E^+)\times \mathbb{R}^4
$$
of vector bundles over $\bar{C}_{\mathbf{n}/S}^{G_S}(E^+)$, which extends a re-scaling of (\ref{eq: trivialization of g_S_TE}).
\end{Proposition}
\begin{proof}
Recall that $\bar{C}_{\mathbf{n}/S}^{G_S}(E^+)$ is defined as the closure of the image of 
$$
\iota:C_{\mathbf{n}/S}(E)\hookrightarrow E^{\mathbf{n}/S}\times 
\prod_{R\subset \mathbf{n}/S} \Bl_{\Delta_{R}^{\mathbf{n}/S}}^{\theta(g_{R})} E^{\mathbf{n}/S}.
$$
We have
$$
\Bl_{\Delta_{\{0,S\}}^{\mathbf{n}/S}}^{\theta(g_{\{0,S\}})} E^{\mathbf{n}/S}
=\Bl_{\Delta_{\{0,S\}}^{\mathbf{n}/S}}^{\theta(f_{\{0,S\}})} E^{\mathbf{n}/S}.
$$ 
Since $f_{\{0,S\}}(\vec{x})=\infty$ for all $\vec{x}\in \Delta_{\{0,S\}}^{\mathbf{n}/S}$,
the linear structure on the normal bundle of $\Delta_{\{0,S\}}^{\mathbf{n}/S}$
is the pullback of the linear structure on the normal bundle of 
$\Delta_{\{0,S\}}^{\{0,S\}}\subset E^{\{0,S\}}$ under the obvious projection map.
In particular, we have
$$
\Bl_{\Delta_{\{0,S\}}^{\mathbf{n}/S}}^{\theta(g_{\{0,S\}})} E^{\mathbf{n}/S}
=\Bl_{\Delta_{\{0,S\}}^{\mathbf{n}/S}}^{\theta(f_{\{0,S\}})} E^{\mathbf{n}/S}\cong 
(E^+)^{\mathbf{n}/S-\{0,S\}}\times \Bl_{\Delta_{\{0,S\}}} E^{\{0,S\}}.
$$
The projection
$$
E^{\mathbf{n}/S}\times 
\prod_{R\subset \mathbf{n}/S} \Bl_{\Delta_{R}^{\mathbf{n}/S}}^{\theta(g_{R})} E^{\mathbf{n}/S}
\to \Bl_{\Delta_{\{0,S\}}^{\mathbf{n}/S}}^{\theta(g_{\{0,S\}})} E^{\mathbf{n}/S}
\to 
\Bl_{\Delta_{\{0,S\}}} E^{\{0,S\}}
$$
induces a continuous map
$$
\bar{C}_{\mathbf{n}/S}^{G_S}(E^+)\to \Bl_{\Delta_{\{0,S\}}} E^{\{0,S\}}.
$$

Notice that $E^{\{0,S\}}\cong E^+$ and  $\Delta_{\{0,S\}}$ corresponds to the 
infinity section of $E^+$. Hence $\Bl_{\Delta_{\{0,S\}}} E^{\{0,S\}}$ can be identified with
the blowup of $E^+$ along the infinity section.
The composition
\begin{equation*}
\begin{tikzcd}
\bar{C}_{\mathbf{n}/S}^{G_S}(E^+)\to \Bl_{\Delta_{\{0,S\}}} E^{\{0,S\}}\cong \Bl_\infty E^+ \arrow[r,"p_0"] &E^+.
\end{tikzcd}
\end{equation*}
is just the  map $\bar{f}_S$ defined below \eqref{eq_def_g_S'}.
According to Lemma \ref{lemma_framing_of_tangent_bundle_of_blowup_S^4} below, the pullback $p_0^\ast T^vE^+$ has a trivialization which extends (\ref{eq: trivialization of T^vE}) after re-scaling. 
Hence $\bar{f}_S^\ast T^v E^+$ also has a trivialization which extends the pullback of (\ref{eq: trivialization of T^vE}) after re-scaling. 

We have the decomposition 
$$
\bar{C}_{\mathbf{n}/S}^{G_S}(E^+)=\bar{f}_S^{-1}(E)\cup 
\bar{f}_S^{-1}\{x\in E^+| |x|\ge M_1\}.
$$
By definition, we have  $\bar{g}_S'({\bar{f}_S^{-1}(E)})\subset E$. Therefore, 
$\bar{g}_S'^\ast T^v E^+|_{\bar{f}_S^{-1}(E)}$ has a trivialization from the pullback of (\ref{eq: trivialization of T^vE}).
Again by definition, we have
$$
\bar{g}_S'|_{\{x\in E^+| |x|\ge M_1\}}=\bar{f}_S.
$$ 
Therefore, $\bar{g}_S'^\ast T^v E^+|_{\{x\in E^+| |x|\ge M_1\}}$ has a trivialization which extends 
 the pullback of (\ref{eq: trivialization of T^vE}) after re-scaling. By gluing these two trivilizations together after re-scaling, we obtain a trivialization of $\bar{g}_S'^\ast T^v E^+$ 
over $\bar{C}_{\mathbf{n}/S}^{G_S}(E^+)$.
\end{proof}

We view $S^4$ as the unit sphere in $\mathbb R^5$ and set $\infty=(1,0\cdots,0), \mathbf{0}=(-1,0\cdots,0)$. Using the stereographic projection
\begin{align*}
St:S^4-\{\infty\} &\to \mathbb{R}^4\\
(x_0,x_1,\cdots,x_4) &\mapsto (\frac{x_1}{1-x_0},\cdots, \frac{x_4}{1-x_0})
\end{align*}
we could identify $S^4$ with $\mathbb{R}^4\cup\{\infty\}$ and this $\mathbb{R}^4$ can also be identified with
$T_\mathbf{0}S^4$.
Similarly we have another 
stereographic projection
\begin{align*}
St':S^4-\{\mathbf{0}\} &\to \mathbb{R}^4\\
(x_0,x_1,\cdots,x_4) &\mapsto (\frac{x_1}{1+x_0},\cdots, \frac{x_4}{1+x_0})
\end{align*}
and this $\mathbb{R}^4$ can  be identified with $T_\infty S^4$.
\begin{Lemma}\label{lemma_framing_of_tangent_bundle_of_blowup_S^4}
Let $p:\Bl_\infty S^4\to S^4=\mathbb{R}^4\cup \{\infty\}$ be the blowup map. Then $p^\ast TS^4$ has
a trivialization whose restriction on $\Bl_\infty S^4-p^{-1}(\infty)=\mathbb{R}^4$ coincides with  the standard trivialization of $T\mathbb{R}^4$ after re-scaling.
\end{Lemma}
\begin{proof}
We could identify the exceptional divisor with the unit sphere $S(T_\infty S^4)$ and 
$\Bl_\infty S^4-\{\mathbf{0}\}$ with $S(T_\infty S^4)\times [0,\infty)$. The restriction of the blowup map to 
the complement of $\mathbf{0}$ can be described by
\begin{align*}
p':S(T_\infty S^4)\times [0,\infty) &\xrightarrow{f} T_\infty S^4 \xrightarrow{St'^{-1}} S^4 -\{\mathbf{0}\}\\
(v, t) &\mapsto tv   \mapsto (\frac{1-t^2}{1+t^2}, \frac{2tv}{1+t^2} )
\end{align*} 
Away from the exceptional divisor, we have
\begin{align*}
S(T_\infty S^4)\times (0,\infty) &\xrightarrow{\cong} T_\infty S^4-\{0\} \xrightarrow{St'^{-1}} S^4-
\{\mathbf{0},\infty \} \xrightarrow{St} \mathbb{R}^4-\{0\}\\
(v, t) &\mapsto tv   \mapsto (\frac{1-t^2}{1+t^2}, \frac{2tv}{1+t^2} ) \mapsto \frac{v}{t}
\end{align*} 
A direct calculation shows that 
$$
(T_{tv}(St\circ St'^{-1}) )^{-1} (\frac{\partial}{\partial x_i})=t^2(\frac{\partial}{\partial x_i}-
2\sum_{j=1}^4v_i v_j \frac{\partial}{\partial x_j})
$$
where $(\frac{\partial}{\partial x_i})$ is the standard framing of $T\mathbb{R}^4$. Now it is clear that
$$
(T_{tv}(St\circ St'^{-1}) )^{-1} (|x|^2\frac{\partial}{\partial x_i})=(\frac{\partial}{\partial x_i}-
2\sum_{j=1}^4v_i v_j \frac{\partial}{\partial x_j})
$$
extends to a framing on $p'^\ast TS^4$. Let $\rho:\mathbb{R}^4\to \mathbb{R}^+$ be a smooth function which
is equal to $|x|^2$ when $|x|$ is large, then the above discussion shows that $(\rho \frac{\partial}{\partial x_i})$
extends to a framing for $p^\ast TS^4$.
\end{proof}

An immediate corollary of Proposition \ref{prop_trivialization_pullback_TvE+} is the following.
\begin{Proposition}\label{prop_partial_S_no_0_product_decomposition}
If $0\not\in S$, then we have 
$$
\partial_S\cong \bar{C}_{\mathbf{n}/S}^{G_S}(E^+)\times \bar{C}_S^\ast(\mathbb{R}^4)
$$
\end{Proposition}

Given a subset $S\subset\mathbf{n}$, we define 
\[
V_S=\begin{cases} T_\infty S^4 \quad &\text{if }0\in S;\\
   \mathbb{R}^4                \quad &\text{if } 0\notin S.
\end{cases}
\]
Now Proposition \ref{prop_partial_S_0_product_decomposition} and 
Proposition \ref{prop_partial_S_no_0_product_decomposition}
can be summarized as
\begin{equation}\label{eq_partial_S_product_decomposition}
\varphi_{\mathbf{n},S}^E :\bar{C}_{\mathbf{n}/S}^{G_{S}}(E^+)\times \bar{C}^\ast_S(V_S) \cong \partial_S \bar{C}_{\mathbf{n}}^{G}(E^+)
.
\end{equation}
We call $\varphi_{\mathbf{n},S}^E$ the \emph{assembling homeomorphism} for the codimension-1 face $\bar{C}_{\mathbf{n}}^{G}(E^+)$.

Let $W=T^{v}E^{+}$. Given any any sets $A\subset D$ with $|A|\geq 2$, we have a similar assembling homeomorphism
\begin{equation}\label{eq: linear assembling  homeomorphism}
\varphi_{D,A}^W:\bar{C}^*_{D/A}(W)\times_{E^+} \bar{C}^\ast_A(W) \cong \partial_A \bar{C}_{D}(W)    
\end{equation}
The definition is straightforward since we are working with smooth configuration spaces.

We will also consider the assembling homeomorphism for codimension-2 faces. For any $A,D\subset \mathbf{n}$ with $|A|, |D|\geq 2$ and $A\cap D=\emptyset$, we let $\mathbf{n}/\{A,D\}$ be the quotient set 
of $\mathbf{n}$ that identifies all points in $A$ as a single point and identifies all points in $D$ as another point.
Similar to the discussion before Proposition \ref{prop_closure_Delta_S_Cn/S}, we can view diagonals
properly
included in $\Delta_{\{A,D\}}\cong E^{\mathbf{n}/\{A,D\}}$ as diagonals of  $E^{\mathbf{n}/\{A,D\}}$ and the subset
$G_{A,D}\subset G$ associated to these diagonals as a collection of maps for diagonals of $E^{\mathbf{n}/\{A,D\}}$. Hence
we can define $\bar{C}_{\mathbf{n}/\{A,D\}}^{G_{A,D}}(E^+)$.

\begin{Lemma}\label{Lem: codim-2 face (2)}
By iterating \eqref{eq_partial_S_product_decomposition} twice in two different orders, we have the following commutative diagram
\begin{equation*}
\begin{tikzcd}
\bar{C}_{\mathbf{n}/\{A,D\}}^{G_{A,D}}(E^{+})\times\bar{C}^{*}_{A}(V_{A})\times \bar{C}^{*}_{D}(V_{D})\ar[d,"\varphi^{E}_{\mathbf{n}/D,A}\times \id "]\ar[rr," \cong"]  && \bar{C}^{G_{A,D}}_{\mathbf{n}/\{A,D\}}(E^{+})\times\bar{C}^{*}_{D}(V_{D})\times \bar{C}^{*}_{A}(V_{A}) \ar[d,"\varphi^{E}_{\mathbf{n}/A,D}\times \id "]\\
\partial_{A} \bar{C}^{G_D}_{\mathbf{n}/D}(E^{+})\times \bar{C}^{*}_{D}(V_D) \ar[d,hook]  &&  \partial_{D} \bar{C}^{G_A}_{\mathbf{n}/A}(E^{+})\times \bar{C}^{*}_{A}(V_A) \ar[d,hook] \\
\bar{C}^{G_D}_{\mathbf{n}/D}(E^{+})\times\bar{C}^{*}_{D}(V_D) \ar[dr,"\varphi^{E}_{\mathbf{n},D}"]  &&  \bar{C}^{G_A}_{\mathbf{n}/A}(E^{+})\times \bar{C}^{*}_{A}(V_A) \ar[dl,"\varphi^{E}_{\mathbf{n},A}"'] \\
 & \bar{C}_{\mathbf{n}}^G(E^{+}) &\\
\end{tikzcd}
\end{equation*}
Furthermore, both compositions give the assembling homeomorphism 
\begin{equation}\label{eq: codim-2 assembling (2)}
\bar{C}_{\mathbf{n}/\{A,D\}}^{G_{A,D}}(E^{+})\times\bar{C}^{*}_{A}(V_A)\times \bar{C}^{*}_{D}(V_D)\xrightarrow{\cong }
\partial_{\{A,D\}} \bar{C}_{\mathbf{n}}^{G}(E^+)
\end{equation}
\end{Lemma}
\begin{proof}
The diagram can be obtained by iterating \eqref{eq_partial_S_product_decomposition} twice in two different 
orders. To show the commutativity we replace 
$$
\bar{C}_{\mathbf{n}/\{A,D\}}^{G_{A,D}}(E^{+})\times\bar{C}^{*}_{A}(V_{A})\times \bar{C}^{*}_{D}(V_{D})
$$
by its interior 
$$
{C}_{\mathbf{n}/\{A,D\}}(E)\times {C}^{*}_{A}(V_{A})\times {C}^{*}_{D}(V_{D}).
$$
The description \eqref{def_partial_F} 
implies that the composition maps  
from
$
{C}_{\mathbf{n}/\{A,D\}}(E)\times {C}^{*}_{A}(V_{A})\times {C}^{*}_{D}(V_{D})
$
to $\bar{C}_{\mathbf{n}}^G(E^{+})$ along the two paths of the diagram
are the same. Moreover, it is  
 a homeomorphism onto $\mathring{\partial}_{\{A,D\}} \bar{C}_{\mathbf{n}}^{G}(E^+)$. Since both
 $
\bar{C}_{\mathbf{n}/\{A,D\}}^{G_{A,D}}(E^{+})\times\bar{C}^{*}_{A}(V_{A})\times \bar{C}^{*}_{D}(V_{D})
 $
and $\partial_{\{A,D\}} \bar{C}_{\mathbf{n}}^{G}(E^+)$ are compact and Hausdorff, the 
original diagram commutes and
the composition map is a continuous map onto $\partial_{\{A,D\}} \bar{C}_{\mathbf{n}}^{G}(E^+)$. Since this map is also
injective, it is a homeomorphism onto $\partial_{\{A,D\}} \bar{C}_{\mathbf{n}}^{G}(E^+)$.

\end{proof}

We omit the proof of the following lemma since its proof is similar to that of Lemma \ref{Lem: codim-2 face (2)}.
\begin{Lemma}\label{Lem: codim-2 face (1)} Let $W=T^{v}E^{+}$. For any $A\subset D\subset \mathbf{n}$ with $|A|\geq 2$, by iterating Proposition \ref{prop_partial_A_pullback} we have the following commutative diagram 
\begin{equation*}
\begin{tikzcd}
\bar{C}_{\mathbf{n}/D}^{G_D}(E^{+})\times_{E^+}\bar{C}^{*}_{D/A}(W)\times_{E^{+}} \bar{C}^{*}_{A}(W)\ar[d,"\id\times \varphi^{W}_{D,A}"]\ar[rr," \varphi^{E}_{\mathbf{n}/A,D/A}\times \id"]  && \partial_{D/A}\bar{C}_{\mathbf{n}/A}^{G_A}(E^{+})\times_{E^{+}} \bar{C}^{*}_{A}(W)\ar[d,hook] \\
\bar{C}_{\mathbf{n}/D}^{G_D}(E^{+})\times_{E^+}\partial_{A} \bar{C}^{*}_{D}(W) \ar[d,hook] &&  \bar{C}_{\mathbf{n}/A}^{G_A}(E^{+})\times_{E^{+}} \bar{C}^{*}_{A}(W)\ar[d,"\varphi^{E}_{\mathbf{n},A}"]\\
\bar{C}_{\mathbf{n}/D}^{G_D}(E^{+})\times_{E^+} \bar{C}^{*}_{D}(W)  \ar[rr,"\varphi^{E}_{\mathbf{n},D}"] && 
\bar{C}_{\mathbf{n}}^G(E^{+}).
\end{tikzcd}
\end{equation*}
Furthremore, both compositions give the assembling homeomorphism 
\begin{equation}\label{eq: codim-2 assembling (1)}
\bar{C}_{\mathbf{n}/D}^{G_D}(E^{+})\times_{E^+}\bar{C}^{*}_{D/A}(W)\times_{E^{+}} \bar{C}^{*}_{A}(W)\xrightarrow{\cong}
\partial_{\{A,D\}}\bar{C}_{\mathbf{n}}^G(E^{+})= \partial_{A} \bar{C}_{\mathbf{n}}^G(E^{+})\cap  \partial_{D} \bar{C}_{\mathbf{n}}^G(E^{+})
\end{equation}
\end{Lemma}


\section{Construction of propagators}\label{section: propagators}
In this section, we will use $\cM$ to denote the set $\cM_\mathbf{n}$ defined below Equation \eqref{eq_normal_bundle_of_Delta_S}. And we let 
\[
\cM^{+}=\cM\cup \{\emptyset\}.
\]
As before, we let  $E_{D}\to B$ be a boundary-trivialized topological disk bundle with a framing. We fix a choice of observer system 
\[
G=\{g^{\mathcal{S}}_{S}: \Delta_{\mathcal{S}}\to E^{+}\mid \mathcal{S}\in\cM, S\in \mathcal{S}\}.
\] 
for $E^\mathbf{n}$.
Then for any $\mathcal{S}\in \cM^{+}$, we have a compactification $\bar{C}^{G_{\mathcal{S}}}_{\mathbf{n}/\mathcal{S}}(E^{+})$. 
Here we set $\mathbf{n}/\emptyset=\mathbf{n}$ and set $G_{\emptyset}=G$.

 The main purpose of this section is to inductively define a collection of maps
\[
\cP=\{P^{\mathcal{S}}_{[i],[j]}: \bar{C}^{G_\mathcal{S}}_{\mathbf{n}/\mathcal{S}}(E^{+})\to S^3\mid \mathcal{S}\in \cM^{+}, [i]\neq [j]\in \mathbf{n}/\mathcal{S}-\{[0]\}\}
\]called a ``propagator system''.  

\begin{Definition}\label{defi: sphere sequence} A \emph{sphere sequence}, denoted by $X(-)$, is a sequence 
\begin{equation}\label{eq: sequence of maps on S3}
X(1)\xrightarrow{h_1} X(2)\xrightarrow{h_2} X(3) \xrightarrow{h_3} \cdots
\end{equation}
of maps  between $X(i)=S^{3}$ such that the following conditions are satisfied
\begin{itemize}
\item For any $m$, the mapping degree $\deg(h_{m})$ is nonzero.
\item $h_{1}$ is even. I.e., 
\begin{equation}\label{eq: evenness of h1}
h_{1}(x)=h_{1}(-x),\ \forall x\in X(1)
\end{equation}
\item Consider the map 
\begin{equation}\label{eq: g_l}
g_{l}=h_{l}\circ \cdots \circ h_{1}: X(1)\to X(l).    
\end{equation}
Then for any integer $n$, there exists $l$ such that $n \mid \operatorname{deg}(g_{l})$. 
\end{itemize}
\end{Definition}

\begin{Definition} We say two maps $f: Y\to X(l)$ and $g: Y\to X(l')$ are \textbf{equivalent} if they are equal after being mapped to $X(m)$ for $m$ large enough. 
\end{Definition}

 We fix the identification 
\[
\bar{C}^{*}_{\mathbf{1}}(\bR^{4})\xrightarrow{\cong}X(1),\ [(x_0,x_1)]\mapsto \frac{x_{1}-x_{0}}{|x_{1}-x_{0}|}.
\]
This allows us to define the map 
\[
\psi_{i,j}:\bar{C}^{*}_{S}(\mathbb{R}^{4})\to X(1)
\]
for any two different elements $i,j$ in $S$. 

Next, we discuss the unusual coordinate on $\bar{C}^{*}_{S}(T_{\infty} S^{4})$. Here we again identify $S^{4}$ as the unit sphere in $\mathbb{R}^{5}$ and set $\infty=(1,0,\cdots,0)$. This gives the identification $\mathbb{R}^{4}\cong T_{\infty}S^{4}$ defined by $x\mapsto (0,x)$. 
We consider the map 
\[\tau_{\mathbf{n}}: C^{*}_{\mathbf{n}}(\bR^{4})\xrightarrow{\cong}  C^{*}_{\mathbf{n}}(T_{\infty}S^4)\]
defined by 
\begin{equation}\label{eq: inversion}
[(0,x_1,\cdots,x_n)]\mapsto [(0,\frac{x_1}{|x_1|^2},\cdots,\frac{x_n}{|x_n|^2})].
\end{equation}
The following lemma is proved in \cite[Section 2.3.3]{Wa21}.
\begin{Lemma}\label{lem: special coordinate} The map (\ref{eq: inversion}) extends to a homeomorphism $\tau_{\mathbf{n}}:  \bar{C}^{*}_{\mathbf{n}}(\bR^{4})\xrightarrow{\cong}  \bar{C}^{*}_{\mathbf{n}}(T_{\infty}S^4)$. 
\end{Lemma}

For a general finite ordered $S$, we define the homeomorphism \[\tau_{S}:\bar{C}^{*}_{S}(\bR^{4})\xrightarrow{\cong}  \bar{C}^{*}_{S}(T_{\infty}S^4)\]
via the unique order preserving bijection $S\leftrightarrow \mathbf{n}$ where $n=|S|-1$. 

\begin{Lemma}\label{lem: tau compatible with psi} Let $A$ be a subset of $S$ that contains $\min S$. Consider the maps 
\[\psi:\bar{C}^{*}_{S}(\bR^4)\to  \bar{C}^{*}_{A}(\bR^4) ,\quad \psi^{\infty}:\bar{C}^{*}_{S}(T_{\infty}S^4)\to  \bar{C}^{*}_{A}(T_{\infty}S^4)\]
induced by $A\hookrightarrow S$. Then the following diagram commutes
\begin{equation} \label{eq: tau compatible with psi}
\begin{tikzcd}
\bar{C}^{*}_{S}(T_{\infty}S^4) \arrow[d,"\tau_{S}^{-1} "] \arrow[r,"\psi^{\infty} "]  & \bar{C}^{*}_{A}(T_{\infty}S^4)\arrow[d,"\tau_{A}^{-1} "] \\
\bar{C}^{*}_{A}(\bR^4) \arrow[r,"\psi "]&  \bar{C}^{*}_{A}(\bR^4) .
\end{tikzcd}
\end{equation}
\end{Lemma}
 \begin{proof} For points in $C^{*}_{S}(T_{\infty}S^4)$, the diagram (\ref{eq: tau compatible with psi}) is straightforward to veriy. Since $C^{*}_{S}(T_{\infty}S^4)$ is dense in $\bar{C}^{*}_{S}(T_{\infty}S^4)$,  (\ref{eq: tau compatible with psi}) holds for all points in $\bar{C}^{*}_{S}(T_{\infty}S^4)$.
\end{proof}

For any $\cS\in \cM^{+}$ and any $A\subset \mathbf{n}/\cS$ with $|A|\geq 2$, we let \[\cS\cdot A=\{\pi^{-1}(a)\mid a\in (\mathbf{n}/\cS)/A, |\pi^{-1}(a)|\geq 2\}\in\cM^{+},\] where $\pi:\mathbf{n}\to (\mathbf{n}/\cS)/A$ is the quotient map. We use  
\[
p_1:  \partial_{A} \bar{C}^{G_{\cS}}_{\mathbf{n}/\cS}(E^{+})\to \bar{C}^{G_{\cS\cdot A}}_{\mathbf{n}/(\cS \cdot A)}(E^{+})
\]
to denote the projection via the assembling homeomorphism \eqref{eq_partial_S_product_decomposition}
(with $\mathbf{n}$ replaced by $\mathbf{n}/\mathcal{S}$ and $S=A$). We  also let 
\[
p'_{2}: \partial_{A} \bar{C}^{G_{\cS}}_{\mathbf{n}/\cS}(E^{+})\to \bar{C}^{*}_{A}(V_{A})
\]
be the projection via \eqref{eq_partial_S_product_decomposition}. When $[0] \notin A$, we set $p_2=p_2'$. When $[0] \in A$, 
we set \begin{equation}\label{eq: p_2}
 p_2=\tau_A^{-1}\circ p'_{2} .  
\end{equation} In any case, we have a map 
\[
p_2:  \partial_{A} \bar{C}^{G_{\cS}}_{\mathbf{n}/\cS}(E^{+})\to \bar{C}^{*}_{A}(\bR^{4}).
\]
\begin{Definition}\label{defi: propagator} Let $E$  be a boundary-trivialized, framed topological disk bundle and let $G$ be an observer system for $E^n$. For any $\cS\in \cM^{+}$, define the set 
\[
I(\cS):=\{([i],[j])\in (\mathbf{n}/\cS) \times (\mathbf{n}/\cS)\mid [i]\neq [0], [j]\neq [0], [i]\neq [j]\}.
\]
A propagator is a collection of maps 
\[P^{\cS}=\{P^{\cS}_{[i],[j]}: \bar{C}^{G_{\cS}}_{\mathbf{n}/\cS}(E^{+})\to X(l_\cS)\mid ([i],[j])\in I(\cS)\}\]
such that the following conditions holds for any $([i],[j])\in I(\cS)$ and any $A\subset \mathbf{n}/\cS$ with $|A|\geq 2$.
\begin{enumerate}
    \item $P^{\cS}_{[i],[j]}=P^{\cS}_{[j],[i]}$.
\item Suppose $[i]=[j]$ in $(\mathbf{n}/\cS)/A$. Then the restriction of $P^{\cS}_{[i],[j]}$ to $\partial_{A}\bar{C}^{G_{\cS}}_{\mathbf{n}/\cS}(E^{+})$ is equivalent to the composition
\begin{equation}\label{eq boundary condition: i,j the same}
\partial_{A}\bar{C}^{\cS}_{\mathbf{n}/\cS}(E^{+})\xrightarrow{p_2}  \bar{C}^{*}_{A}(\bR^{4})\xrightarrow{\psi_{[i],[j]}} X(1).    
\end{equation}
\item Suppose in $(\mathbf{n}/\cS)/A$, we have $[i]=[0]\neq [j]$. Then the restriction of $P^{\cS}_{[i],[j]}$ to $\partial_{A}\bar{C}^{G_{\cS}}_{\mathbf{n}/\cS}(E^{+})$ is equivalent to the composition
\begin{equation}\label{eq boundary condition: i,0 the same}
\partial_{A}\bar{C}^{\cS}_{\mathbf{n}/\cS}(E^{+})\xrightarrow{p_2} \bar{C}^{*}_{A}(\bR^{4})\xrightarrow{\psi_{[0],[i]}} X(1).
\end{equation}
\item  Suppose in $(\mathbf{n}/\cS)/A$, we have $[j]=[0]\neq [i]$. Then the restriction of $P^{\cS}_{[i],[j]}$ to $\partial_{A}\bar{C}^{G_{\cS}}_{\mathbf{n}/\cS}(E^{+})$ is equivalent to the composition
\begin{equation}\label{eq boundary condition: j,0 the same}
\partial_{A}\bar{C}^{G_{\cS}}_{\mathbf{n}/\cS}(E^{+})\xrightarrow{p_2} \bar{C}^{*}_{A}(\bR^{4})\xrightarrow{\psi_{[0],[j]}} X(1).
\end{equation}

\end{enumerate}
\end{Definition}
\begin{Definition}\label{defi: compatible propagator} Suppose  $\cM'$ is a subset of $\cM^{+}$ satisfying 
 $\cS\cdot A\in \cM'$ for any $\cS\in \cM'$, any $A\subset \mathbf{n}/\cS$ with $|A|\geq 2$.
We say a collection of propagators  $\{\cP^{\cS}\}_{\cS\in \cM'}$ is compatible if and any $i,j\in \mathbf{n}$ such that $([i],[j])\in I(\cS\cdot A)$, 
 the restriction of $P^{\cS}_{[i],[j]}$ on $\partial_{A}\bar{C}^{G_{\cS}}_{\mathbf{n}/\cS}(E^{+})$ is equivalent to the composition
\begin{equation}\label{eq boundary condition: i,j,0 different}
\partial_{A}\bar{C}^{G_{\cS}}_{\mathbf{n}/\cS}(E^{+})\xrightarrow{p_1} \bar{C}^{G_{\cS\cdot A}}_{\mathbf{n}/(\cS\cdot A)}(E^{+}) \xrightarrow{P^{\cS\cdot A}_{[i],[j]}} X(l_{\cS\cdot A}).    
\end{equation}
In the special case $\cM'=\cM^{+}$, we call a compatible collection of propagators \[\cP=\{\cP^{\cS}\}_{\cS\in \cM^{+}}\]
a \emph{propagator system}. Furthermore suppose $l_{\cS}=l$ for all $\cS\in \cM+$. Then we call $\cP$ a \emph{uniform propagator system} of degree $d$, where $d$ is the mapping degree of the map (\ref{eq: g_l}).
\end{Definition}

\begin{Definition} Let $G$ and $G'$ be two observer systems for $E$. Consider the propagator systems
\[\mathcal{P}=\{P^{\cS}_{[i],[j]}: \bar{C}^{G_{\cS}}_{\mathbf{n}/\cS}(E^{+})\to X(l_\cS)\}_{i,j,\cS}
\]
and 
\[\mathcal{Q}=\{Q^{\cS}_{[i],[j]}: \bar{C}^{G'_{\cS}}_{\mathbf{n}/\cS}(E^{+})\to X(l'_\cS)\}_{i,j,\cS}.
\]
A concordance from $(G,\mathcal{P})$ to $(G',\mathcal{Q})$ is a homotopy $\tilde{G}$ from $G$ to $G'$ together with a propagator system
\[
\mathcal{R}=\{R^{\cS}_{[i],[j]}: \bar{C}^{\tilde{G}_{\cS}}_{\mathbf{n}/\cS}(I\times E^{+})\to X(l''_\cS)\}_{i,j,\cS}.
\]
such that for any $i,j,\cS$, the restrictions of $R^{\cS}_{[i],[j]}$ to $\{0\}\times E^{+}$ and $\{1\}\times E^{+}$ are equivalent to $P^{\cS}_{[i],[j]}$ and $Q^{\cS}_{[i],[j]}$ respectively.
\end{Definition}

The main result of this section is the following theorem.
\begin{Theorem}\label{Thm: propagator}
Let $E_{D}$ be a boundary trivialized, framed topological disk bundle and let $G$ be an observer system in Lemma \ref{lemma_g_S_existence_family}. There exists a uniform propagator system $\mathcal{P}$ \[\mathcal{P}=\{P^{\cS}_{[i],[j]}: \bar{C}^{G_{\cS}}_{\mathbf{n}/\cS}(E^{+})\to X(l)\mid \cS\in \cM^{+}, ([i],[j])\in I(\cS)\}
.\]
 Furthermore, the pair $(G,\mathcal{P})$ is unique up to concordance.
\end{Theorem}
To prove Theorem \ref{Thm: propagator}, we start with some homotopy theoretic lemmas. Let $X(\infty)$ be the homotopy colimit (i.e. the mapping telescope) of this sequnce
\[X(\infty):=(\bigsqcup_{i\geq 1} ([i-1,i]\times X(i)))/\sim,\]
where $(i,x)\sim (i,h_{i}(x))$ for all $i\geq 1$ and $x\in X(i)$. Then 
\[
\pi_{m}(X(\infty))=\operatorname{colimit}(\pi_{m}(S^3)\xrightarrow{h_{1,*}}\pi_{m}(S^3)\xrightarrow{h_{2,*}}\cdots)
\]
Since $\pi_{m}(S^3)$ is torsion for any $m>3$,  we see that $X(\infty)$ is an Eilenberg-Maclane space $K(\bQ,3)$. We have a natural inclusion $X(i)\hookrightarrow X(\infty)$ defined by $x\mapsto (i-1,x)$. 

\begin{Lemma}\label{lem: extension} Let $(W,Y)$ be a cofibration which has the homotopy type of a finite CW pair. 
Let $1$ be the generator of $H^{3}(S^3;\bQ)$. Given a contiunuous map $f: Y\to X(n)$, suppose 
\begin{equation}
f^{*}(1)\in \operatorname{\ker}(H^{3}(Y;\bQ)\xrightarrow{\partial} H^{4}(W,Y;\bQ)).
\end{equation}
 Then for $m$ large enough, the extension problem 
\begin{equation*}
\begin{tikzcd}
Y \arrow[d,hook] \arrow[r,"f"]  & X(n) \arrow[d]\\
W\arrow[r, dotted] &  X(m)
\end{tikzcd}.
\end{equation*}
also has a solution.
\end{Lemma}
\begin{proof}
Since $(W,Y)$ is a cofibration, it suffices to solve the extension problem up to homotopy. Hence we may assume $(W,Y)$ is a finite CW pair. Since $X(\infty)$ is an Eilenberg-Maclane space, we have isomorphisms 
\[
[W,X(\infty)]\cong H^{3}(W;\bQ),\ [Y,X(\infty)]\cong H^{3}(Y;\bQ).
\] 
By our assumption, there exists a map $\tilde{f}$ that fits into the commutative diagram
\begin{equation*}
\begin{tikzcd}
Y \arrow[d,hook] \arrow[r,"f"]  & X(m) \arrow[d]\\
W\arrow[r, "\tilde{f}"] &  X(\infty)
\end{tikzcd}.
\end{equation*}
Since $W$ is compact, we have  
\[\tilde{f}(W)\subset (\bigsqcup_{1\leq i\leq m} ([i-1,i]\times X(i)))/\sim\]
Since $(\bigsqcup_{1\leq i\leq m} ([i-1,i]\times X(i)))/\sim$ deformation retracts to $X(m)$, the proof is finished.
\end{proof}

\begin{remark}
We will apply Lemma \ref{lem: extension} to the case where $W$ is a topological manifold and $Y$ is a union of codimension-1 faces. Since $Y$ has a collar, $(W,Y)$ is a cofibration. And since all topological manifolds are absolute neighborhood retracts, $(W,Y)$ has the homotopy type of a CW pair. (Here we use cellular approximation theorem and the mapping cylinder construction to turn $Y\hookrightarrow W$ into the inclusion map of a subcomplex).
\end{remark}

\begin{Lemma}\label{lem: relative homology}$H^{k}(\bar{C}^G_{\mathbf{n}}(E^{+}),\partial^{v}\bar{C}^G_{\mathbf{n}}(E^{+});\bZ)=0$ for $n\ge 2, k\leq 4$. Here $\partial^{v}\bar{C}^G_{\mathbf{n}}(E^{+})=\bar{C}^G_{\mathbf{n}}(E^{+})-C_{\mathbf{n}}(E^+)$ denotes the vertical boundary of $\bar{C}^G_{\mathbf{n}}(E^{+})$. 
\end{Lemma}
\begin{proof} When $n=2$, this was proved in \cite[Lemma 2.10]{Wa21} (in the smooth case) using the Serre spectral sequence. For our case we have
$$
H^{k}(\bar{C}^G_{\mathbf{n}}(E^{+}),\partial^{v}\bar{C}^G_{\mathbf{n}}(E^{+}))\cong
H^{k}_c ({C}_{\mathbf{n}}(E)).
$$
If $\partial B=\emptyset$, then by Poincar\'e duality we have
$$
H^{k}_c ({C}_{\mathbf{n}}(E))\cong H_{4n+\dim B-k}( C_\mathbf{n}(E)).
$$

Using the spectral sequence for the fibration 
$$
\mathbb{R}^4-\{n-1~\text{points}\}\to C_n(\mathbb{R}^4) \to C_{n-1}(\mathbb{R}^{4})
$$
 and induction we could deduce that $H_i(C_n(\mathbb{R}^4))=0$ unless $3|i$ and $i\le 3(n-1)$.
Since the fiber of $C_\mathbf{n}(E)\to B$ is $C_n(\mathbb{R}^4)$,  the spectral sequence for
this fiber bundle implies $H_j(C_\mathbf{n}(E))=0$ for all $j> 3(n-1)+\dim B$. Thus 
$H_{4n+\dim B-k}( C_\mathbf{n}(E))=0$ for $k\le 4$, which completes the proof when $\partial B=\emptyset$.

When $\partial B\neq \emptyset$ and $k\le 4$, we have  
$$
H^{k}(\bar{C}^G_{\mathbf{n}}(E^{+}),\partial \bar{C}^G_{\mathbf{n}}(E^{+}))\cong
H^{k}_c ({C}_{\mathbf{n}}(E)|_{B^\circ})\cong H_{4n+\dim B-k}( C_\mathbf{n}(E))=0
$$
and 
$$
H^{k}(\bar{C}^G_{\mathbf{n}}(E^{+})|_{\partial B},\partial^v \bar{C}^G_{\mathbf{n}}(E^{+})|_{\partial B})=0
$$
since $\partial B$ has no boundary. Now the lemma follows from the  cohomology long exact sequence for 
the three pairs
$$
(\bar{C}^G_{\mathbf{n}}(E^{+})|_{\partial B},\partial^v \bar{C}^G_{\mathbf{n}}(E^{+})|_{\partial B})
\to (\bar{C}^G_{\mathbf{n}}(E^{+}),\partial^{v}\bar{C}^G_{\mathbf{n}}(E^{+}) )
\to (\bar{C}^G_{\mathbf{n}}(E^{+}),\partial \bar{C}^G_{\mathbf{n}}(E^{+})).
$$

\end{proof}

For $3\leq k \leq n+1$, we let 
\[
\cM^{+}_{k}=\{\cS\in \cM^{+}\mid 3\leq |\mathbf{n}/\cS|\leq k\}.
\]

\begin{Lemma}\label{lem: propagator n=2 boundary}
For any $\cS\in \cM^{+}_3$ and $([i],[j])\in I(\cS)$,
there exists a unique continuous map \[P^{\cS,\partial}_{[i],[j]}:  \partial^{v}\bar{C}^{G_{\cS}}_{\mathbf{n}/\cS}(E^+)\to X(1)\] whose restriction to $\partial_{A}\bar{C}^{G_{\cS}}_{\mathbf{n}/\cS}(E^{+})$ equals the composition
\begin{equation}\label{eq: C^2 boundary}
\partial_{A}\bar{C}^{G_{\cS}}_{\mathbf{n}/\cS}(E^{+})\xrightarrow{p_{2}} \bar{C}_{\mathbf{n}/\cS}^{*}(\bR^4)\xrightarrow{\psi_{[k],[l]}} X(1).
\end{equation}
Here 
\[
([k],[l])=\begin{cases} ([i],[j])\text{ when }A=\{[0],[i],[j]\},\ \{[i],[j]\}\\
([0],[i])\text{ when }A=\{[0],[i]\}\\
([j],[0])\text{ when }A=\{[0],[j]\}\\
\end{cases}
\]
\end{Lemma}
\begin{proof} In the special case $B$ is a point and $E$ is smooth, this result was proved in \cite[Lemma 2.7]{Wa21}. For the general case, one needs to check that various definitions are compatible with each other on the intersections $\partial_{A}\bar{C}_{\mathbf{2}}(E^{+})$ for various $A$. Note that these intersections are all contained in $\bar{C}_{\mathbf{2}}(E^{+}-E_{D})$, which is just the product of $B$ with the configuration space $\bar{C}_{\mathbf{2}}(S^{4}-D^4)$. So the general case follows from the special case. 
\end{proof}

\begin{Corollary}\label{cor: propagator for 2 points} For any $\cS\in \cM^{+}_3$, there exists a propagator \[\cP^{\cS}=\{P^{\cS}_{[i],[j]}: \bar{C}^{G_{\cS}}_{\mathbf{n}/\cS}(E^+)\to X(l_\cS)\mid ([i],[j])\in I(\cS)\}.\] 
\end{Corollary}
\begin{proof} First we assume $[i]\leq [j]$. Then by Lemma \ref{lem: relative homology} and Lemma \ref{lem: extension}, the extension problem 
\begin{equation} \label{eq: extension C_2 to X_m}
\begin{tikzcd}
\partial^{v}\bar{C}^{G_{\cS}}_{\mathbf{n}/\cS}(E^+)\arrow[d] \arrow[r,"P^{\cS,\partial}_{[i],[j]}"]  & X(1) \arrow[d]\\
\bar{C}^{G_{\cS}}_{\mathbf{n}/\cS}(E^+) \arrow[r,dotted]&   X(m).
\end{tikzcd}
\end{equation}
has a solution when $m$ is larger then some constant $m_{i,j}$. We set $l_{\cS}=\max\limits_{i,j}m_{i,j}$ and let $P^{\cS}_{[i],[j]}$ be any solution.

When $[j]< [i]$, we set $P^{\cS}_{[i],[j]}=P^{\cS}_{[j],[i]}$. Note that the boundary condition is still satisfied because the map $h_{1}: X(1)\to X(2)$ is even. (See Definition \ref{defi: sphere sequence}).
\end{proof}

Now we generalize the above construction to $\cS\in \cM^{+}_{k}$ with $k>3$. Consider the constant map 
\[
q: \bar{C}^{*}_{\{0,1\}}(T_{\infty}S^{4})\to \infty\in S^4
\]
Since $q$ is a constant map, we have a canonical trivialization 
\[q^\ast TS^4\cong  \bar{C}^{*}_{\{0,1\}}(T_\infty S^4)
\times T_{\infty}S^{4}.\]
This gives the homeomorphism 
\[
\bar{C}^{*}_{\{1,2\}}(q^\ast TS^4)\cong\bar{C}^{*}_{\{0,1\}}(T_{\infty}S^{4})\times \bar{C}^{*}_{\{1,2\}}(T_{\infty}S^{4})
\]
On the other hand, under the obvious identification $\bar{C}^{*}_{\{0,1\}}(T_{\infty}S^{4})\cong \partial \operatorname{Bl}_{\infty}S^{4}$, the map $q$ is restricted from the blowup map $p:\Bl_\infty S^4\to S^4$. Therefore, by restricting the trivilization in Lemma \ref{lemma_framing_of_tangent_bundle_of_blowup_S^4}, we obtain a second trivialization 
\[q^\ast TS^4\cong \bar{C}^{*}_{\{0,1\}}(q^\ast TS^4)\times \bR^{4},\]
which induces another homeomorphism 
\[
\bar{C}^{*}_{\{1,2\}}(q^\ast TS^4)\cong\bar{C}^{*}_{\{0,1\}}(T_{\infty}S^{4})\times \bar{C}^{*}_{\{1,2\}}(\bR^4).
\]
Putting these two homeomorphisms together, we obtain the 
\begin{equation}\label{eq: homeomorphism induced by two trivilizations}
\bar{C}^{*}_{\{0,1\}}(T_{\infty}S^{4})\times \bar{C}^{*}_{\{1,2\}}(T_{\infty}S^{4})\cong\bar{C}^{*}_{\{0,1\}}(T_{\infty}S^{4})\times \bar{C}^{*}_{\{1,2\}}(\bR^{4}).
\end{equation}
We let \[
\xi: \bar{C}^{*}_{\{0,1\}}(T_{\infty}S^{4})\times \bar{C}^{*}_{\{1,2\}}(T_{\infty}S^{4})\to \bar{C}^{*}_{\{1,2\}}(\bR^{4})
\]
be the composition of the homeomorphism (\ref{eq: homeomorphism induced by two trivilizations}) with the projection. 
\begin{Lemma} \label{lem: continuity lemma}
The following diagram commutes 
\begin{equation} \label{diagram: differences between two infinity points}
\begin{tikzcd}
\bar{C}^{*}_{\{0,1\}}(T_{\infty}S^4)\times \bar{C}^{*}_{\{1,2\}}(T_{\infty}S^4)\ar[rr, "\varphi^{T_{\infty}S^{4}}_{\mathbf{2},\{1,2\}}"] \ar[d,"\xi"]& &\bar{C}^{*}_{\mathbf{2}}(T_{\infty}S^4)\ar[d,"\tau_{\mathbf{2}}^{-1}"]\\
\ \bar{C}^{*}_{\{1,2\}}(\bR^{4})\ar[r,"\cong","\psi_{1,2}" '] &  X(1) &  \bar{C}^{*}_{\mathbf{2}}(\bR^4)\ar[l,"\psi_{1,2}"]\\ \end{tikzcd}
\end{equation}
Here we identify $\bar{C}^{*}_{\{0,1\}}(T_{\infty}S^4)$ with $\bar{C}^{*}_{\mathbf{2}/\{1,2\}}(T_{\infty}S^4)$ via the order-preserving bijection $\{0,1\}\cong \mathbf{2}/\{1,2\}$. 
\end{Lemma}

\begin{proof} 
This is a special case of Lemma \ref{lem: propagator n=2 boundary}. It states that when $n=2$, $\cS=\emptyset$, the map $P^{\cS,\partial}_{1,2}$ is well-defined at the intersection between the faces $\partial_{\mathbf{2}}\bar{C}_{\mathbf{2}}(S^{4})$ and $\partial_{\{1,2\}}\bar{C}_{\mathbf{2}}(S^{4})$. 
\end{proof}

We define \[
\xi': \bar{C}^{*}_{\{0,1\}}(T_{\infty}S^{4})\times \bar{C}^{*}_{\{1,2\}}(T_{\infty}S^{4})\to \bar{C}^{*}_{\{0,1\}}(\bR^{4})
\]
as the composition of the projection
$$
\bar{C}^{*}_{\{0,1\}}(T_{\infty}S^{4})\times \bar{C}^{*}_{\{1,2\}}(T_{\infty}S^{4})\to \bar{C}^{*}_{\{0,1\}}(T_{\infty}S^{4})
$$
and the map $\tau_{\mathbf{1}}^{-1}:\bar{C}^{*}_{\{0,1\}}(T_{\infty}S^{4})\to \bar{C}^{*}_{\{0,1\}}(\mathbb{R}^{4})$.
\begin{Lemma} \label{lem: continuity lemma 2}
The following diagram commutes 
\begin{equation} \label{diagram: differences between two infinity points 2}
\begin{tikzcd}
\bar{C}^{*}_{\{0,1\}}(T_{\infty}S^4)\times \bar{C}^{*}_{\{1,2\}}(T_{\infty}S^4)\ar[rr, "\varphi^{T_{\infty}S^{4}}_{\mathbf{2},\{0,1\}}"] \ar[d,"\xi'"]& &\bar{C}^{*}_{\mathbf{2}}(T_{\infty}S^4)\ar[d,"\tau_{\mathbf{2}}^{-1}"]\\
\ \bar{C}^{*}_{\{0,1\}}(\bR^{4})\ar[r,"\cong","\psi_{0,1}" '] &  X(1) &  \bar{C}^{*}_{\mathbf{2}}(\bR^4)\ar[l,"\psi_{1,2}"]\\ \end{tikzcd}
\end{equation}
Here we identify $\bar{C}^{*}_{\{1,2\}}(T_{\infty}S^4)$ with $\bar{C}^{*}_{\mathbf{2}/\{0,1\}}(T_{\infty}S^4)$ via the order-preserving bijection $\{1,2\}\cong \mathbf{2}/\{0,1\}$. 
\end{Lemma}
\begin{proof} This is another special case of Lemma \ref{lem: propagator n=2 boundary}, which states that when $n=2$, $\cS=\emptyset$, the map $P^{\cS,\partial}_{1,2}$ is well-defined at the intersection between the faces $\partial_{\mathbf{2}}\bar{C}_{\mathbf{2}}(S^{4})$ and $\partial_{\{0,1\}}\bar{C}_{\mathbf{2}}(S^{4})$. 
\end{proof}
To  construct $\cP^{\cS}$ inductively over $|\mathbf{n}/\cS|$, 
we need the following proposition. 

\begin{Proposition}\label{prop: boundary propagator}For $3\leq k\leq n$, let $\{\cP^{\cS}\mid \cS\in \cM^{+}_k \}$
be a compatible collection of propagators. Then any $\cS\in \cM^{+}_{k+1}-\cM^{+}_{k}$ and any $([i],[j])\in I(\cS)$, then there exists $l\gg 0$ and a map 
 \[
P^{\cS,\partial}_{[i],[j]}: \partial^{v}\bar{C}^{G_{\cS}}_{\mathbf{n}/\cS}(E^{+})\to X(l)
\] 
such that for any $A\subset \mathbf{n}/\cS$ with $|A|\geq 2$, the restriction of $P^{\cS,\partial}_{[i],[j]}$
 to $\partial_{A}\bar{C}^{G_{\cS}}_{\mathbf{n}/\cS}(E^{+})$ equals $P^{\cS,A}_{[i],[j]}$, where $P^{\cS,A}_{[i],[j]}$ is defined as follows.
\begin{itemize}
\item Suppose $([i],[j])\in I(\cS\cdot A)$, then we define $P^{\cS,A}_{[i],[j]}$ as the composition 
\[
\partial_{A}\bar{C}^{G_{\cS}}_{\mathbf{n}/\cS}(E^+)\xrightarrow{(\ref{eq boundary condition: i,j,0 different})} X(l_{\cS\cdot A})\to X(l)
\]
\item Suppose $[i]=[j]$ in $\mathbf{n}/(\cS\cdot A)$, then we define $P^{\cS,A}_{[i],[j]}$ as the composition 
\[
\partial_{A}\bar{C}^{G_{\cS}}_{\mathbf{n}/\cS}(E^+)\xrightarrow{(\ref{eq boundary condition: i,j the same})} X(1)\to X(l).
\]

\item Suppose $[i]=[0]\neq [j]$ in $\mathbf{n}/(\cS\cdot A)$, then we define $P^{\cS,A}_{[i],[j]}$ as the composition 
\[
\partial_{A}\bar{C}^{G_{\cS}}_{\mathbf{n}/\cS}(E^+)\xrightarrow{(\ref{eq boundary condition: i,0 the same})} X(1)\to X(l).
\]
\item Suppose $[j]=[0]\neq [i]$ in $\mathbf{n}/(\cS\cdot A)$, then we define $P^{\cS,A}_{[i],[j]}$ as the composition 
\[
\partial_{A}\bar{C}^{G_{\cS}}_{\mathbf{n}/\cS}(E^+)\xrightarrow{(\ref{eq boundary condition: j,0 the same})} X(1)\to X(l).
\]
\end{itemize}
\end{Proposition}
\begin{proof} It suffices to show that any $A,D\subset \mathbf{n}/\cS$ with $|A|,|D|\geq 2$, and any $([i],[j])\in I(\cS)$, we have  
\[P^{\cS,A}_{[i],[j]}(x)=P^{\cS,D}_{[i],[j]}(x), \quad \forall x\in \partial_{\{A,D\}}\bar{C}^{G_{\cS}}_{\mathbf{n}/\cS}(E^+).\]
Proof of this claim is a case-by-case check using Lemma \ref{Lem: codim-2 face (1)} or Lemma \ref{Lem: codim-2 face (2)}. Since most of the cases are similar, we only focus on two typical cases.

First, we suppose \[\{[0],[i]\}\subset A\subset D,\quad [j]\in D\setminus A.\]
By Definition \ref{defi: propagator}, the map $P^{\cS,D}_{[i],[j]}$ equals the composition 
\[
\partial_{D}\bar{C}^{G_{\cS}}_{\mathbf{n}/\cS}(E^+)\xrightarrow{p_2} \bar{C}^{*}_{D}(\bR^{4})\xrightarrow{\psi_{[i],[j]}} X(1)\to X(l).
\]
The assembling homeomorphism (\ref{eq: codim-2 assembling (1)}) provides a decomposition
\begin{equation}\label{eq: codim2 decomposition}
\partial_{\{A,D\}}\bar{C}^{G_{\cS}}_{\mathbf{n}/\cS}(E^+)\xrightarrow{\cong}\bar{C}^{G_{\cS\cdot D}}_{\mathbf{n}/(\cS\cdot D)}(E^{+})\times\bar{C}^{*}_{D/A}(T_{\infty}S^4)\times\bar{C}^{*}_{A}(T_{\infty}S^4).
\end{equation}
We use $\textrm{pj}_{k}$ to denote the projection to the $k$-th component. Note that the composition 
\[
\partial_{\{A,D\}}\bar{C}^{G_{\cS}}_{\mathbf{n}/\cS}(E^+)\hookrightarrow \partial_{D}\bar{C}^{G_{\cS}}_{\mathbf{n}/\cS}(E^+)\xrightarrow{p_2}\bar{C}^*_{D}(\bR^4)
\]
equals the composition 
\[
\partial_{\{A,D\}}\bar{C}^{G_{\cS}}_{\mathbf{n}/\cS}(E^+)\xrightarrow{(\textrm{pj}_{2},\textrm{pj}_{3})} \bar{C}^{*}_{D/A}(T_{\infty}S^4)\times\bar{C}^{*}_{A}(T_{\infty}S^4)\xrightarrow{\varphi^{T_{\infty}S^4}_{D,A}} \bar{C}^{*}_{D}(T_{\infty}S^4)\xrightarrow{\tau^{-1}_{D}}\bar{C}^{*}_{D}(\mathbb{R}^4).
\]
Hence the restriction of $P^{\cS,D}_{[i],[j]}$ to $\partial_{\{A,D\}}\bar{C}^{G_{\cS}}_{\mathbf{n}/\cS}(E^+)$ can be rewritten as
\begin{equation}\label{eq: partial_A,D 1}
\begin{split}
\partial_{\{A,D\}}\bar{C}^{G_{\cS}}_{\mathbf{n}/\cS}(E^+)\xrightarrow{(\textrm{pj}_{2},\textrm{pj}_{3})} \bar{C}^{*}_{D/A}(T_{\infty}S^4)\times\bar{C}^{*}_{A}(T_{\infty}S^4)\\ \xrightarrow{\varphi^{T_{\infty}S^4}_{D,A}} \bar{C}^{*}_{D}(T_{\infty}S^4)
\xrightarrow{\tau^{-1}_{D}}\bar{C}^{*}_{D}(\mathbb{R}^4)\xrightarrow{\psi_{[i],[j]}} X(1)\to X(l).
\end{split}
\end{equation}
Consider the map $\bar{C}^{*}_{D}(T_{\infty}S^4)\to \bar{C}^{*}_{\mathbf{2}}(T_{\infty}S^4)$
induced by the injection $\mathbf{2}\to D$ sending $0, 1, 2$ to $[0], [i], [j]$ respectively. Then Lemma \ref{lem: continuity lemma 2} implies the following commutative diagram 
\begin{equation*}
\begin{tikzcd}
\bar{C}^{*}_{D/A}(T_{\infty}S^4)\times \bar{C}^{*}_{A}(T_{\infty}S^4)\ar[rrr, "\varphi^{T_{\infty}S^{4}}_{D,A}"] \ar[d,"\text{projection}"]&  & &\bar{C}^{*}_{D}(T_{\infty}S^4)\ar[d,"\tau_{D}^{-1}"]\\
\bar{C}^{*}_{A}(T_{\infty}S^4)\ar[r,"\tau^{-1}_{A}"]& \bar{C}^{*}_{A}(\bR^{4})\ar[r,"\psi_{[0],[i]}" '] &   X(1) &  \bar{C}^{*}_{D}(\bR^4)\ar[l,"\psi_{[i],[j]}"]\\ \end{tikzcd}    
\end{equation*}
So we can rewrite (\ref{eq: partial_A,D 1}) as the composition 
\begin{equation}\label{eq: partial_A,D 2}
\partial_{\{A,D\}}\bar{C}^{G_{\cS}}_{\mathbf{n}/\cS}(E^+)\xrightarrow{\textrm{pj}_{3}} \bar{C}^{*}_{A}(T_{\infty}S^4)
\xrightarrow{\tau^{-1}_{A}}\bar{C}^{*}_{A}(\mathbb{R}^4)\xrightarrow{\psi_{[0],[i]}} X(1)\to X(l).
\end{equation}
By a similar argument, we see that (\ref{eq: partial_A,D 2}) also equals the restriction of $P^{\cS,A}_{[i],[j]}$ to  $\partial_{\{A,D\}}\bar{C}^{G_{\cS}}_{\mathbf{n}/\cS}(E^+)$. This finishes the proof.

All other cases are similar to the previous one, except when 
\[
[i],[j]\in A\subset D,\quad [0]\in D\setminus A.
\]
In this case, the assembling homeomorphism still provides the decomposition (\ref{eq: codim2 decomposition}). But the maps $P^{\cS,A}_{[i],[j]}$ and $P^{\cS,D}_{[i],[j]}$ have different expressions. To describe them, we consider the map
\[
\psi': \bar{C}^{*}_{D/A}(T_{\infty}S^4)\to \bar{C}^{*}_{\mathbf{1}}(T_{\infty}S^4) 
\]
induced by the map $\mathbf{1}\to D/A$ given by $0\mapsto [0], 1\mapsto [A]$. And we consider the map 
\[
\psi'': \bar{C}^{*}_{A}(T_{\infty}S^4)\to \bar{C}^{*}_{\{1,2\}}(T_{\infty}S^4)
\]
induced by the map $\{1,2\}\to A$ given by $1\mapsto [i], 2\mapsto [j]$. Given a point \[x\in  \partial_{\{A,D\}}\bar{C}^{G_{\cS}}_{\mathbf{n}/\cS}(E^+),\] we use $(\alpha,\beta,\gamma)$ to denote its coordinates under the decomposition \eqref{eq: codim2 decomposition}. Unwinding definitions, we see that $P^{\cS,A}_{[i],[j]}(x)$ and  $P^{\cS,D}_{[i],[j]}(x)$ are equal to the image of $(\psi'(\beta),\psi''(\gamma))$ under the composition 
\[
\bar{C}^{*}_{\{0,1\}}(T_{\infty}S^4)\times \bar{C}^{*}_{\{1,2\}}(T_{\infty}S^4)\xrightarrow{f_{A}, f_{D}} X(1)\to X(l) 
\]
Here $f_{A}, f_{D}$ are the  two routes in the diagram (\ref{diagram: differences between two infinity points}). More precisely, $f_{A}$ 
equals the composition 
\[
\bar{C}^{*}_{\{0,1\}}(T_{\infty}S^4)\times \bar{C}^{*}_{\{1,2\}}(T_{\infty}S^4)\xrightarrow{\xi}\bar{C}^{*}_{\{1,2\}}(\bR^{4})\xrightarrow{\psi_{1,2}}X(1)\to X(l)
\]
and $f_{D}$ equals the composition
\[
\begin{split}
\bar{C}^{*}_{\{0,1\}}(T_{\infty}S^4)\times \bar{C}^{*}_{\{1,2\}}(T_{\infty}S^4)&\xrightarrow{\varphi^{T_{\infty}S^{4}}_{[2],\{1,2\}}}\bar{C}^{*}_{\{0,1,2\}}(T_{\infty}S^4)\\
&\xrightarrow{\tau_{[2]}}\bar{C}^{*}_{\{0,1,2\}}(\bR^4)\xrightarrow{\psi_{1,2}}X(1)\to X(l)
\end{split}.
\]
Since (\ref{diagram: differences between two infinity points}) commutes, we have $P^{\cS,A}_{[i],[j]}(x)=P_{[i],[j]}^{\cS,[D]}(x)$. 
\end{proof}

 The following result allows us to construct $\cP^{\cS}$ inductively over $|\mathbf{n}/\cS|$.

\begin{Corollary}\label{cor: propagator induction} Given $3\leq k\leq n$, let $\{\cP^{\cS}\mid \cS\in \cM^{+}_k \}$
be a compatible collection of propagators. Then there exists a collection of propagators 
\[\{\cP^{\cS}\mid \cS\in \cM^{+}_{k+1}-\cM^{+}_{k} \}
\]
such that the union 
\[\{\cP^{\cS}\mid \cS\in \cM^{+}_k \}\cup \{\cP^{\cS}\mid \cS\in \cM^{+}_{k+1}-\cM^{+}_{k}
\}\]
is again a compatible collection of propagators.
\end{Corollary} \begin{proof} For $\cS\in \cM^{+}_{k+1}- \cM^{+}_{k}$ and $([i],[j])\in I(\cS)$, we let
\[
P^{\cS,\partial}_{[i],[j]}: \partial^{v}\bar{C}^{G_{\cS}}_{\mathbf{n}/\cS}(E^{+})\to X(l)
\]
be the map provided by Proposition \ref{prop: boundary propagator}. Then we apply Lemma \ref{lem: relative homology} and Lemma \ref{lem: extension} to solve the extension problem 
\begin{equation} \label{eq: extension C_n to X_m}
\begin{tikzcd}
\partial^{v}\bar{C}^{G_{\cS}}_{\mathbf{n}/\cS}(E^{+})\arrow[d,hook] \arrow[r,"P^{\cS,\partial}_{[i],[j]}"]  & X(l_{n-1}) \arrow[d]\\
\bar{C}^{G_{\cS}}_{\mathbf{n}/\cS}(E^{+}) \arrow[r,dotted, "P^{\cS}_{[i],[j]} "]&   X(l_\cS)
\end{tikzcd}
\end{equation}
for $l_{\cS}$ large enough. When $i>j$, we simply define $P^{\cS}_{[i],[j]}=P^{\cS}_{[j],[i]}$.
\end{proof}
By Corollary \ref{cor: propagator induction}, given any observer system $G$, we can inductively  construct a propagator system $\cP$. By composing with maps in the sphere sequence, we can further turn $\cP$ into a uniform propagator. Thus the proof of Theorem \ref{Thm: propagator} is completed by the following proposition.

\begin{Proposition}\label{prop: propagator unique up to concordance}The pair $(G,\cP)$ is unique up to concordance.
\end{Proposition}
\begin{proof} Consider observer systems $G,G'$ and propagator systems \[\mathcal{P}=\{P^{\cS}_{[i],[j]}: \bar{C}^{G_{\cS}}_{\mathbf{n}/\cS}(E^{+})\to X(l)\}_{i,j,\cS}
\]
and 
\[\mathcal{Q}=\{Q^{\cS}_{[i],[j]}: \bar{C}^{G'_{\cS}}_{\mathbf{n}/\cS}(E^{+})\to X(l')\}_{i,j,\cS}.
\]
We may assume $l=\l'$ by composing with maps in there sphere sequence.  By Lemma \ref{lemma_g_S_existence_family}, there exits a homotopy $\tilde{G}$ from $G$ to $G'$. 

For $\cS\in \cM^{+}$, we consider following subspace of  $\bar{C}^{\tilde{G}_{\cS}}_{\mathbf{n}/\cS}(I\times E^+)$
\[
W_{\cS}:= \bar{C}^{\tilde{G}_{\cS}}_{\mathbf{n}/\cS}(\{0\}\times E^+)\cup  
\partial^{v}\bar{C}^{\tilde{G}_{\cS}}_{\mathbf{n}/\cS}(I\times E^+)
\cup\bar{C}^{\tilde{G}_{\cS}}_{\mathbf{n}/\cS}(\{1\}\times E^+).
\]
Then by Lemma \ref{lem: relative homology}, we have
\begin{equation}\label{eq: homology vanishing}
H^{4}(\bar{C}^{\tilde{G}_{\cS}}_{\mathbf{n}/\cS}(I\times E^+), W_{\cS};\bZ)
=0.    
\end{equation}

For $\cS\in \cM^{+}_{3}$, we apply Lemma \ref{lem: propagator n=2 boundary} to $I\times E^{+}$ and get the map 
\[
\partial^{v}\bar{C}^{\tilde{G}_{\cS}}_{\mathbf{n}/\cS}(I\times E^+)\to X(1).
\]
By composing with the map $X(1)\to X(l)$, we get the map
\begin{equation}\label{eq: RS}
R^{\cS,\partial}_{[i],[j]}:  \partial^{v}\bar{C}^{\tilde{G}_{\cS}}_{\mathbf{n}/\cS}(I\times E^+)\to X(l).    
\end{equation}
Consider the map
\[
P^{\cS,\partial}_{[i],[j]} \cup R^{\cS,\partial}_{[i],[j]}\cup Q^{\cS,\partial}_{[i],[j]}: W^{\cS}\to  X(l).
\]
By Lemma \ref{lem: extension} and (\ref{eq: homology vanishing}), the extension problem 
\begin{equation}\label{diagram: homotopy extension} 
\begin{tikzcd}
W_{\cS}\arrow[d,hook] \arrow[rr,"P^{\cS}_{[i],[j]} \cup R^{\cS,\partial}_{[i],[j]}\cup Q^{\cS}_{[i],[j]}"]&  & X(l) \arrow[d]\\
\bar{C}^{\tilde{G}_{\cS}}_{\mathbf{n}/\cS}(I\times E^+) \arrow[rr,dotted, "R^{\cS}_{[i],[j]}"]& &  X(l'').
\end{tikzcd}
\end{equation}
has a solution when $l''$ is large enough. 

Suppose we have defined $R^{\cS}_{[i],[j]}$ for all $\cS\in \cM^{+}_{k}$. Then for any $\cS\in \cM^{+}_{k+1}- \cM^{+}_{k}$ and any $([i],[j])\in I(\cS)$, we apply Proposition \ref{prop: boundary propagator} to the bundle $I\times E^+$ and define the map (\ref{eq: RS}). Then we define $R^{\cS}_{[i],[j]}$ by solving the extension problem (\ref{diagram: homotopy extension}). This finishes the inductive construction of the propagator system  
\[
\mathcal{R}=\{R^{\cS}_{[i],[j]}\mid \cS\in \cM^{+}, ([i],[j])\in I(\cS)\}.
\]
The pair $(\tilde{G},\mathcal{R} )$ gives the desired concordance from $(G,\cP)$ to $(G',\mathcal{Q})$.
\end{proof}


\section{Kontsevich integral for topological families}\label{section: Kontsevich integral}
In this section, we define configuration space integrals for boundary trivialized, framed topological disk bundles. 

\subsection{Definition of Kontsevich's characteristic classes} Consider a boundary trivialized, framed disk bundle 
\[D^{4}\to E_{D}\xrightarrow{\pi } B\]
over a closed, oriented topological manifold $B$. We fix an positive even integer $n$. We also fix a choice of an observer system 
\[
G:=\{g^{\cS}_{S}: \Delta_{\cS}\to E^{+}\mid S\in \cS, \cS\in \cM\}
\]
and a uniform propagator system 
\[
\cP:=\{P^{\cS}_{[i],[j]}: \bar{C}^{G_{\cS}}_{\bfn/\cS}(E^{+})\to X(l)\}
\] 
of degree $d$. (See Definition \ref{defi: propagator}.)

Consider a decorated, oriented admissible graph $(\Gamma,\rho,o)$ with $p$ vertices and $q$ edges such that $n=2q-2p$. 
 Via the decoration $\rho$, we have the bijection 
\begin{equation}\label{eq: v-gamma as quotient}
\tilde{\rho}: \{1,\cdots, n\}/\cS  \xrightarrow{\cong}  V(\Gamma)
\end{equation}
for some $\cS\in \cM $ such that $[0]\notin S$ for any $S\in \cS$. For each edge $e\in E(\Gamma)$, we pick any direction on $e$ and let
\[i(e), j(e)\in \{1,\cdots, n\}/\cS\subset \bfn/\cS
\]
be elements that correspond to the initial vertex and the end vertex of $e$ respectively. 
 Let 
$X(\Gamma)=\prod_{e\in E(\Gamma)}X(l)$. Then we can define the map 
\begin{equation}\label{eq: graph propagator}
P(\Gamma,\rho):= \prod_{e\in E(\Gamma)} P^{\cS}_{i(e),j(e)}: \bar{C}^{G_{\cS}}_{\bfn/\cS}(E^{+})\to X(\Gamma).
\end{equation}
This map does not depend on direction of the edges because $P^{\cS}_{[i],[j]}$ is symmetric in $i,j$. (See Definition \ref{defi: propagator}.)

Recall that the Čech cochain complex of a space $Y$ can be defined as the direct limit
\[
C_{\cech}^{*}(Y):=\lim_{\substack{\longrightarrow\\ \text{open cover }\mathcal{U} \text{ of }Y}} C_{\cech}^{*}(\mathcal{U}).
\] 

\begin{Definition}\label{Defi: standard cocyle} We say a cocycle $\xi\in C_{\cech}^{3q}(\prod_{q}S^3)$ is standard if the following conditions are satisfied:
\begin{itemize}
    \item $\xi$ represents the standard generator in $H^{3q}(\prod_{q}S^3;\mathbb{Q})$.
    \item Let $\Sigma_{q}$ be the group of permutations on $\{1,\cdots, q\}$. Then $\xi$ is symmetric under the action of $\Sigma_{q}$. That means 
    \begin{equation}\label{eq: permutation symmetry}
    g^{*}(\xi)=\operatorname{sign}(g)\cdot \xi,\quad \forall g\in \Sigma_{q}.
    \end{equation}
Here $\operatorname{sign}(g)=1$ if $g$ is even and $\operatorname{sign}(g)=-1$ if $g$ is odd. This can be achieved by replacing $\xi$ with
\[\frac{1}{q!}\sum_{g\in \Sigma_{q}}(\operatorname{sign}(g)\cdot g^{*}(\xi)).\]
\end{itemize}
\end{Definition}
We pick a standard Čech cocyle $\xi$. Consider a bijection 
\[\iota:E(\Gamma)\xrightarrow{\cong} \{1,\cdots, q\}\] which respects the orientation $o$. This gives a homeomorphism 
\[
X(\Gamma)\xrightarrow{\cong }\prod_{q}S^3
\]
so we can pull back $\xi$ to $\xi_{o}\in C_{\cech}^{3q}(X(\Gamma))$. By (\ref{eq: permutation symmetry}), $\xi_{o}$ is independent with the bijection $\iota$ and we have $\xi_{o}=-\xi_{\bar{o}}$.

\begin{Definition}Let $(\Gamma,\rho,o)$ be an admissible, decorated, oriented graph and let $D^4\to E_{D}\xrightarrow{\pi} B$ be a boundary trivialized disk bundle over a closed manifold $B$ of dimension $m(\Gamma)+\frac{n(\Gamma)}{2}$, equipped with a framing on it vertical normal tangent bundle. We define the configuration space integral 
\[
I^{\Top}_{(\Gamma,\rho,o)}(\pi, \cP, G, \xi)=\frac{1}{d^{\frac{3n(\Gamma)-2m(\Gamma)}{2}}}\langle P(\Gamma, \rho)^{*}(\xi_{o}), [\bar{C}^{G_{\cS}}_{\bfn/\cS}(E^+)]\rangle\in \bQ
\]
Here $d$ denotes the degree of $\mathcal{P}$.
\end{Definition}

In general, $I^{\Top}_{(\Gamma,\rho,o)}(\pi, \cP, G, \xi)$ depends on choices of various auxiliary data. To get an actual invariant of the bundle, we need to take certain linear combination of graphs.

\begin{Definition}
Suppose $B$ is a closed manifold of dimension $m+\frac{n}{2}$. For any cochain \[\eta=\sum_{1\leq i\leq l} a_{i}\cdot (\Gamma_{i},\rho_{i},o_{i})\in \widetilde{\cG}^{m,n}
 \]
We define 
 \[I^{\Top}_{\eta}(\pi):=\sum_{1\leq i\leq l}a_{i}\cdot  I^{\Top}_{(\Gamma_{i},\rho_{i},o_{i})}(\pi, \cP, G, \xi)\in \bQ.\]
For a cohomology class $[\eta]\in H^{m,n}(\widetilde{\cG})$, we pick any representative $\eta$ and define  
\[
I^{\Top}_{[\eta]}(\pi)=I^{\Top}_{\eta}(\pi).
\]
 \end{Definition}

The main result of this section is the following theorem.
\begin{Theorem}\label{thm: I-top welldefined} $I^{\Top}_{[\eta]}(\pi)\in \bQ$ doesn't depend on auxiliary choices so is an invariant of the framed disk bundle $\pi$. Furthermore, if two such bundles $\pi_{1},\pi_2$ are framed cobordant to each other, then one has   $I^{\Top}_{[\eta]}(\pi_1)=I^{\Top}_{[\eta]}(\pi_2)$.
\end{Theorem}

\subsection{Boundary behavior of $P(\Gamma,\rho)$}
To prove Theorem \ref{thm:  I-top welldefined}, we study the boundary behavior of the map 
$P(\Gamma,\rho)$ for any $(\Gamma,\rho,o)\in \cG^{m,n}$. We identify   $V(\Gamma)$ with $\{1,\cdots,n\}/\cS$ using (\ref{eq: v-gamma as quotient}). Note that many arguments here are adapted from Watanabe's argument \cite{Wa21} to the Čech setting. For a subset $A\subset \{1,\cdots,n\}/\cS$, we let $\Gamma_{A}$ be the subgraph of $\Gamma$ consisting of all vertices in $A$ and edges between them, and we let $\Gamma/A$ be the quotient graph obtained by collapsing all vertices in $A$ and all edges between them. Using the bijection
\[\rho/A: V(\Gamma/A)\xrightarrow{\cong} \{1,\cdots, n\}/(\cS\cdot A)\subset \bfn/(\cS\cdot A),\] we can define the map 
\[
P(\Gamma/{A},\rho/A):= \prod_{e\in E(\Gamma/A)} P^{A}_{i(e),j(e)}: \bar{C}^{G_{\cS\cdot A}}_{\bfn/(\cS\cdot A)}(E^{+})\to X(\Gamma/A).
\]
Using the bijection
\begin{equation}\label{eq: rhoA}
\rho_{A}: V(\Gamma_{A})\xrightarrow{\cong} A.
\end{equation}
we can define the map 
\[
\psi(\Gamma_{A},\rho_{A}):\bar{C}^{*}_{A}(\mathbb{R}^{4})\xrightarrow{\prod_{e\in E(\Gamma_A)} \psi _{i(e),j(e)}}  \prod_{e\in E(\Gamma_A)} X(1)\xrightarrow{\prod g_{l}} X(\Gamma_{A}).
\]

\begin{Lemma}\label{Lem: boundary behavior 1} For any $A\subset n/\cS$ with $[0]\notin A$ and $|A|\geq 2$, the following diagram commutes 
\begin{equation*} 
\begin{tikzcd}[column sep = large]
\partial_{A}\bar{C}^{G_{\cS}}_{\bfn/\cS}(E^{+})\ar[rrr,"P(\Gamma{,}\rho)"]\arrow[d,"(p_{1}{,}p_{2})"] & & & X(\Gamma)\arrow[d,"\cong"]\\
\bar{C}^{G_{(\cS\cdot A
)}}_{\bfn/(\cS\cdot A)}(E^{+})\times \bar{C}^{*}_{A}(\bR^{4}) \ar[rrr,"P(\Gamma/A{,}\rho/A)\times \psi(\Gamma_{A}{,}\rho_{A})"]& &  & X(\Gamma/A)\times X(\Gamma_A)\\
\end{tikzcd}
\end{equation*}
Here the homeomorphism
$X(\Gamma)\cong X(\Gamma/A)\times X(\Gamma_{A})$
is induced by the natural decomposition 
$E(\Gamma)=E(\Gamma/A)\sqcup E(\Gamma_{A})$.
\end{Lemma}
\begin{proof}
For any edge $e\in E(\Gamma/A)\subset E(\Gamma)$, by (\ref{eq boundary condition: i,j,0 different}), we have the following commutative diagram  

\begin{equation*} 
\begin{tikzcd}[column sep = large]
\partial_{A}\bar{C}^{G/\cS}_{\bfn/\cS}(E^{+})\ar[rrr,"P^{\cS}_{i(e),j(e)}"]\arrow[d,"p_{1}"] & & & X(l)\arrow[d,equal]\\
\bar{C}^{G_{(\cS\cdot A
)}}_{\bfn/(\cS\cdot A
)}(E^{+})\ar[rrr,"P^{\cS\cdot A}_{i(e),j(e)} "]& &  & X(l)\\
\end{tikzcd}
\end{equation*}
For any edge $e\in E(\Gamma_A)$, by (\ref{eq boundary condition: i,j the same}), we have the following commutative diagram 
\begin{equation*} 
\begin{tikzcd}[column sep = large]
\partial_{A}\bar{C}^{G_{\cS}}_{\bfn/\cS}(E^{+})\ar[rrr,"P^{\cS}_{i(e),j(e)}"]\arrow[d,"p_{2}"] & & & X(l)\arrow[d,equal]\\
\bar{C}^{*}_{A}(\bR^{4}) \ar[rr,"\psi_{i(e),j(e)} "]&  & X(1)\ar[r,"g_{l}"] & X(l)\\
\end{tikzcd}
\end{equation*}
The proof is finished by putting these diagrams together.    
\end{proof}

Now suppose $A$ is a subset of $\bfn/\cS$ that contains $[0]$. Let $A^{c}=\bfn/\cS -A$. Then (\ref{eq: v-gamma as quotient}) induces an bijection
\[
\rho_{A^{c}}: V(\Gamma_{A^{c}})\xrightarrow{\cong} \bfn/(\cS\cdot A) -\{[0]\}\subset \bfn/(\cS\cdot A).
\]
So we can define the map
\[
P(\Gamma_{A^{c}},\rho_{A^{c}}):= \prod_{e\in E(\Gamma_{A^{c}})} P^{A}_{i(e),j(e)}: \bar{C}^{G_{A}}_{\bfn/(\cS\cdot A)}(E^{+})\to X(\Gamma_{A^{c}}).
\]
We also consider the bijection 
\begin{equation}\label{eq: rho/Ac}
\rho/A^{c}: V(\Gamma/A^{c})\xrightarrow{\cong} A    
\end{equation}
which sends the unique collapsed vertex to $[0]$ and sends all other vertices of $V(\Gamma/A^{c})$ to themselves.  This bijection allows us to define the map 
\[
\psi(\Gamma/A^{c},\rho/A^{c}):\bar{C}^{*}_{A}(\mathbb{R}^{4})\xrightarrow{\prod_{e\in E(\Gamma/A^{c})} \psi _{i(e),j(e)}}  \prod_{e\in E(\Gamma/A^{c})} X(1)\xrightarrow{\prod g_{l}} X(\Gamma/A^{c})
\]

\begin{Lemma}\label{Lem: boundary behavior 2} For any $A\subset \bfn/\cS$ with $[0]\in A$ and $|A|\geq 2$, the following diagram commutes 
\begin{equation*} 
\begin{tikzcd}[column sep = large]
\partial_{A}\bar{C}^{G_{\cS}}_{\bfn/ \mathcal{S}}(E^{+})\ar[rrr,"P(\Gamma{,}\rho)"]\arrow[d,"(p_{1}{,}p_{2})"] & & & X(\Gamma)\arrow[d,"\cong"]\\
\bar{C}^{G_{\cS\cdot A}}_{\bfn/(\cS\cdot A)}(E^{+})\times \bar{C}^{*}_{A}(\bR^{4}) \ar[rrr,"P(\Gamma_{A^{c}}{,}\rho_{A^{c}})\times \psi(\Gamma/A^{c}{,}\rho/A^{c})"]& &  & X(\Gamma_{A^{c}})\times X(\Gamma/A^{c})\\
\end{tikzcd}
\end{equation*}
Here the homeomorphism
$X(\Gamma)\cong X(\Gamma_{A^{c}})\times X(\Gamma/A^{c})$
is induced by the natural decomposition 
$E(\Gamma)=E(\Gamma_{A^{c}})\sqcup E(\Gamma/A^{c})$.
\end{Lemma}

\begin{proof} The proof is the same as Lemma \ref{Lem: boundary behavior 1}.
\end{proof}

\begin{Proposition}\label{prop: boundary vanishing |A|>2} Suppose the base $B$ has dimension $m(\Gamma)+\frac{n(\Gamma)}{2}+1$. Then for any $A\subset \bfn/\cS$ with $|A|\geq 3$, we have 
\[
\langle P(\Gamma, \rho)^{*}(\xi), [\partial_{A}\bar{C}^{G}_{\bfn}(E^+)]\rangle =0.
\]
\end{Proposition}
\begin{proof} The proof is adapted from \cite[Lemma E.3]{Wa21} using Čech cohomology. 
 We set 
\[
\Gamma_{1}=
\begin{cases}\Gamma/{A}&\text{ if }[0]\notin A\\
\Gamma_{A^{c}}&\text{ if }[0]\in A
\end{cases}
\]
and we set 
\[\Gamma_{2}=
\begin{cases}\Gamma_{A}&\text{ if }[0]\notin A\\
\Gamma/A^{c}&\text{ if }[0]\in A
\end{cases}
\]
In both cases, we have a decomposition 
\[
E(\Gamma)\cong E(\Gamma_1)\sqcup E(\Gamma_2).
\]
Recall that $\Gamma$ has $p$ vertices and $q$ edges. We also assume $\Gamma_{2}$ has $p'$ vertices and $q'$ edges. Then $\Gamma_1$ has $p+1-p'$ vertices and $q-q'$ edges. Moreover, we have the bijection (\ref{eq: rhoA}) or (\ref{eq: rho/Ac}) between $V(\Gamma_2)$ and $A$. This implies $p'=|A|\geq 3$. 

There are several cases to consider. 
Suppose all vertices of $\Gamma_{2}$ are at least 3-valent. Then we have 
\begin{equation}\label{eq: m'n'}
2q'\geq 3p'.
\end{equation}
By Lemma \ref{Lem: boundary behavior 1} and Lemma \ref{Lem: boundary behavior 2}, the map  $P(\Gamma,\rho)|_{\partial_{A}\bar{C}_{\bfn/\mathcal{S}}(E^{+})}$ has two factorizations 
\begin{equation*}
\begin{tikzcd} &X(\Gamma_1)\times \bar{C}^{*}_{A}(\bR^{4})\ar[dr] & \\ 
\partial_{A}\bar{C}^{G}_{\bfn/\mathcal{S}}(E^{+})\ar[ur]\ar[dr]&& X(\Gamma)\\
&\bar{C}^{G_{\cS\cdot A}}_{\bfn/\cS\cdot A}(E^{+})\times  X(\Gamma_2)\ar[ur]& 
\end{tikzcd}
\end{equation*}
We can compute that 
\[
\Dim (X(\Gamma_1)\times \bar{C}^{*}_{A}(\bR^{4})) =3q+(4p'-5-3q') 
\]
and  
\[
\Dim (\bar{C}^{G_{A}}_{\bfn/A}(E^{+})\times  X(\Gamma_2))=3q-(4p'-5-3q').
\]
So $P(\Gamma, \rho)^{*}(\xi)\neq 0$ only if $4p'-5-3q'=0$. But then (\ref{eq: m'n'})  implies $p'+10\leq 0$. This is impossible.

Suppose some vertex $u\in A$ of $\Gamma_{2}$ is $2$-valent. Let $u$ be connected to the vertices $v,w\in A$ by edges $e_1$ and $e_2$ respectively.  We consider the homeomorphism 
\[
\bar{\nu}:  \bar{C}^{*}_{A}(\bR^{4}) \to \bar{C}^{*}_{A}(\bR^{4}) \]
that extends the map  $\nu: C_{A}(\bR^{4})\to C_{A}(\bR^{4})$ defined by 
\[
 (\cdots, x_{v},\cdots, x_{u}, \cdots, x_{w},\cdots )\mapsto (\cdots, x_{v},\cdots,x_{v}+x_{w}-x_{u},\cdots, x_{w},\cdots ).
\]
Consider the map 
\[
P_{A}:=(p_{1},p_{2})^{-1}\circ P(\Gamma,\rho): \bar{C}^{G_{\cS\cdot A}}_{\bfn/(\cS\cdot A)}(E^{+})\times \bar{C}^{*}_{A}(\bR^{4})\to X(\Gamma) 
\]
Then the following diagram commutes  
\begin{equation*}
\begin{tikzcd}
\bar{C}^{G_{\cS\cdot A}}_{\bfn/(\cS\cdot A)}(E^{+})\times \bar{C}^{*}_{A}(\bR^{4})\ar[d,"\id\times \bar{\nu}"]\ar[r,"P_{A}"]& X(\Gamma)\ar[d,"g"]\\
\bar{C}^{G_{\cS\cdot A}}_{\bfn/(\cS\cdot A)}(E^{+})\times \bar{C}^{*}_{A}(\bR^{4})\ar[r,"P_{A}"]& X(\Gamma)\\
\end{tikzcd}
\end{equation*}
Here $g$ is the map that switches the $e_{1}$ and $e_{2}$ factors. Since $\id\times \bar{\nu}$ is orientation preserving, we have  
\[
\langle P_{A}^{*}(\xi), [\bar{C}^{G_{\cS\cdot A}}_{\bfn/(\cS\cdot A)}(E^{+})\times \bar{C}^{*}_{A}(\bR^{4})]\rangle =\langle P_{A}^{*}(g^{*}(\xi)), [\bar{C}^{G_{\cS\cdot A}}_{\bfn/(\cS\cdot A)}(E^{+})\times \bar{C}^{*}_{A}(\bR^{4})]\rangle
\] 
But since $g^{*}(\xi)=-\xi$ by (\ref{eq: permutation symmetry}), we see that $\langle P_{A}^{*}(\xi), [\bar{C}^{G_{\cS\cdot A}}_{\bfn/(\cS\cdot A)}(E^{+})\times \bar{C}^{*}_{A}(\bR^{4})]\rangle$ must be zero.

Suppose some vertex $v\in A$ of $\Gamma_{2}$ is $1$-valent, connected to another vertex $u$. We will show that the pullback $\xi_{A}:=P^{*}_{A}(\xi)$ equals zero. Take any point 
\[p=(\cdots, x_{u},\cdots, x_{v},\cdots)\in C_{\bfn/\cS\cdot A}(E^{+})\times C^{*}_{A}(\mathbb{R}^{4}).\]
Let $L\subset \bR^{4}$ be the plane that passes through  $x_{v}$ and is orthogonal to $x_{u}-x_{v}$. Because of the symmetry 
\[
P_{A}(\cdots, x_{u},\cdots, x_{v},\cdots)=P_{A}(\cdots, x_{u},\cdots, tx_{v}+(1-t)x_{u},\cdots), \forall t\in (0,1)
\] 
we see that in a neighborhood $U$ of $p$, the map $P_{A}$ factors through the manifold 
\[C_{\bfn/(\cS\cdot A)}(E^{+})\times C^{*}_{A-\{v\}}(\mathbb{R}^{4})\times L
\]
which has dimension $3q-1$. So $\xi_{A}|_{U}=0$. This shows that $\xi_{A}$ equals zero on $C_{\bfn/(\cS\cdot A)}(E^{+})\times C^{*}_{A}(\mathbb{R}^{4})$. Taking closure, we get $\xi_{A}=0$. The same argument show that $\xi_{A}=0$ if $\Gamma_{2}$ has a $0$-valent vertex.
\end{proof}
\begin{Proposition}\label{prop: boundary vanishing |A|=2} Suppose the base $B$ has dimension $m(\Gamma)+\frac{n(\Gamma)}{2}+1$. Then for any $A\subset \bfn$ with $|A|= 2$, we have 
\[
\langle P(\Gamma, \rho)^{*}(\xi), [\partial_{A}\bar{C}^{G_{\cS}}_{\bfn/\cS}(E^+)]\rangle =0.
\]
if at least one of the following conditions is satisfied:
\begin{enumerate}
\item $[0]\in A$;
\item  $\Gamma/A$ has multiple edges between two vertices;
\item $\Gamma$ has no edge that connects the the  two elements of $A$.
\end{enumerate}
\end{Proposition}
\begin{proof} We use the same notations as in the proof of Proposition \ref{prop: boundary vanishing |A|>2}.

(1) Let $A=\{[0],v\}$. Since $v$ is at least 3-valent in $\Gamma$, there are multiple edges of the graph $\Gamma_2=\Gamma/A^{c}$ that connect the vertex $v$ and the vertex $[A^{c}]$. We pick two of these edges $e,e'$. By Lemma \ref{Lem: boundary behavior 2}, the following diagram commutes 
\[\begin{tikzcd}
\bar{C}^{G_{\cS\cdot A}}_{\bfn/(\cS\cdot A)}(E^{+})\times \bar{C}^{*}_{A}(\bR^{4})\ar[d,"\id"]\ar[r,"P_{A}"]&  X(\Gamma)\ar[d,"g"]\\
\bar{C}^{G_{\cS\cdot A}}_{\bfn/(\cS\cdot A)}(E^{+})\times \bar{C}^{*}_{A}(\bR^{4})\ar[r,"P_{A}"]& X(\Gamma)\\
\end{tikzcd}
\]
Here $g$ is the permutation that switches the factors $e$ and $e'$. Since 
$g^{*}(\xi)=-\xi$, we see that 
\[
P(\Gamma, \rho)^{*}(\xi)|_{\partial_{A}\bar{C}^{G_{\cS}}_{\bfn/\cS}(E^+)}=P^{*}_{A}(\xi)=0.
\]

The proof of (2) is similar. By (1), we may assume $[0]\notin A$. Then $\Gamma_{1}=\Gamma/A$ has multiple edges between two vertices. Then we use the permutation symmetry of $\xi$ to show that $P^{*}_{A}(\xi)=0$.

In case (3), we may again assume $[0]\notin A$. The graph $\Gamma_{2}=\Gamma_{A}$ has no edges, so $P_{A}$ factors through the manifold $\bar{C}^{G_{\cS\cdot A}}_{\bfn/(\cS\cdot A)}(E^{+})$, whose dimension is strictly smaller than the degree of $\xi$. So $P_{A}^{*}(\xi)=0$ again.
\end{proof}
In later discussions, we will use an edge $e$ of $\Gamma$ to denote the subset $A\subset \bfn/\cS$ consisting of its end points. 

\subsection{Behavior under cobordism} Now we study the change of $I^{\Top}_{(\Gamma,\rho,o)}(\pi, \cP, G, \xi)$
under framed cobordisms. 

Let $B_{0}, B_{1}$ be closed topological manifolds of dimension $m+\frac{n}{2}$ and let $\tilde{B}$ be an oriented cobordism from $B_0$ to $B_1$. Let $\tilde{\pi}:\tilde{E}_{D}\to \tilde{B}$ be a boundary trivialized disk bundle with a framing. Let $\pi_{i}: E^{i}_{D}\to B_{i}$ be the restriction of $\tilde{\pi}$. 

Consider an observer system 
\[
\widetilde{G}:=\{g^{\cS}_{S}: \Delta_{\cS}\to \tilde{E}^{+}\mid S\in \cS, \cS\in \cM\}
\]
and uniform propagator system 
\[
\widetilde{\cP}:=\{P^{\cS}_{[i],[j]}: \bar{C}^{G_{\cS}}_{\bfn/\cS}(\tilde{E}^{+})\to X(l)\}
\] 
of degree $d$. By restricting $\widetilde{G}$ and $\widetilde{\cP}$, we obtain a propagator sysmtem $G^{i}$ and an observer system $\cP^{i}$ for the bundle $\pi_{i}$.

\begin{Proposition}\label{prop: cobordism difference} Let $(\Gamma,\rho,o)\in \cG^{m,n}$ be an oriented decorated graph. As before, we use $\rho$ to identify $V(\Gamma)$ with $\bfn/\cS$ for some $\cS$. Then one has
\[
I^{\Top}_{(\Gamma,\rho,o)}(\pi_{1}, G^1,\cP^1, \xi)-I^{\Top}_{(\Gamma,\rho,o)}(\pi_{0}, G^0,\cP^0, \xi)=\frac{1}{d^{\frac{3n-2m}{2}}}\sum_{e\in E'(\Gamma)} 
\langle P(\Gamma, \rho)^{*}(\xi), [\partial_{e}\bar{C}^{G_{\cS}}_{\bfn/\cS}(\tilde{E}^+)]\rangle
\] 
\end{Proposition}
\begin{proof} Wa have a map 
\[
P(\Gamma, \rho):\bar{C}^{G_{\cS}}_{\bfn/\cS}(\tilde{E}^+)\to X(\Gamma).
\]
Since $\bar{C}^{G_{\cS}}_{\bfn/\cS}(\tilde{E}^+)$ is a topological manifold of dimension $3|E(\Gamma)|+1$ and it boundary has a decomposition \[\bar{C}^{G^{0}_{\cS}}_{\bfn/\cS}(E^{+}_0)\bigcup \bar{C}^{G^{1}_{\cS}}_{\bfn/\cS}(E^{+}_1)\bigcup (\bigcup\limits_{A\subset \bfn/\cS, |A|\geq 2} \partial_{A}\bar{C}^{\tilde{G}_{\cS}}_{\bfn/\cS}(\tilde{E}^+)), \]
 we see that 
\[
I^{\Top}_{(\Gamma,o)}(\pi_{1}, G^1,\cP^1, \xi)-I^{\Top}_{(\Gamma,o)}(\pi_{0}, G^0,\cP^0, \xi)=\sum_{A\subset \bfn, |A|\geq 2}\langle P(\Gamma, \rho)^{*}(\xi), [\partial_{A}\bar{C}^{G_{\cS}}_{\bfn/\cS}(\tilde{E}^+)]\rangle 
\]
By Proposition \ref{prop: boundary vanishing |A|>2} and Proposition \ref{prop: boundary vanishing |A|=2}, only those 2-element subset $A$ that correspond to $e\in E'(\Gamma)$ can contribute the right hand side. This finishes the proof.
\end{proof}
\begin{Proposition}\label{prop: cobordism invariance, propagator}  Given any cocycle 
\[
\eta\in \ker(\widetilde{\delta}^{m,n}: \widetilde{\cG}^{m,n} \to \widetilde{\cG}^{m+1,n}),\]
we have 
\[
I^{\Top}_{\eta}(\pi_{0}, G^0,\cP^0, \xi)=I^{\Top}_{\eta}(\pi_{1}, G^1,\cP^1, \xi)
\] 
\end{Proposition}
\begin{proof} 
We let 
\[
\eta =\sum_{1\leq i\leq l} a_{i}\cdot (\Gamma_{i},\rho_{i},o_{i}).
\]
Then by Proposition \ref{prop: cobordism difference}, we have 
\begin{equation}\label{eq: boundary difference}
I^{\Top}_{\eta}(\pi_{1}, G^1,\cP^1, \xi)-I^{\Top}_{\eta}(\pi_{0}, G^0,\cP^0, \xi)=\frac{1}{d^{\frac{3n-2m}{2}}}\sum_{\substack{1\leq i\leq l\\
e\in E'(\Gamma_{i})}} 
a_{i}\cdot \langle P(\Gamma_{i}, \rho_{i})^{*}(\xi), [\partial_{e}\bar{C}^{G_{\cS_i}}_{\bfn/\cS_{i}}(\tilde{E}^+)]\rangle
\end{equation}
Note the  following observation: given $e_{i}\in E'(\Gamma_{i})$ and $e_{j}\in E'(\Gamma_{j})$, suppose we have an isomorphism between the decorated graphes
\[
\tau: (\Gamma_{i}/e_{i},\rho_{i}/e_{i})\xrightarrow{\cong}  (\Gamma_{j}/e_{j},\rho_{j}/e_{j}).
\]
Then we have 
\[
\langle P(\Gamma_{i}, \rho_{i})^{*}(\xi), [\partial_{e}\bar{C}^{G_{\cS}}_{\bfn/\cS_{i}}(\tilde{E}^+)]\rangle=\pm \langle P(\Gamma_{j}, \rho_{j})^{*}(\xi), [\partial_{e}\bar{C}^{G_{\cS}}_{\bfn/\cS_{j}}(\tilde{E}^+)]\rangle.\]
Here we take $+$ if $\tau$ is orientation preserving and we take $-$ if $\tau$ is orientation reversing. Since $\eta$ is a cocycle, contributions of terms on the right side of (\ref{eq: boundary difference}) cancel each other. This finishes the proof.
\end{proof}

\subsection{Invariance} 
\begin{Proposition}\label{prop: independence with P and G} Let $\pi:E_D\to B$ be a boundary trivialized, framed disk bundle over a manifold $B$ of dimension $m+\frac{n}{2}$. Let $\eta\in \widetilde{\cG}^{m,n}$ be a cocycle. Then for any standard Čech cocycle $\xi$ and any two choices of observer-propagator systems $(G^0,\cP^0)$, and $(G^1,\cP^1)$, we have 
\[
I^{\Top}_{\eta}(\pi, G^0,\cP^0, \xi)=I^{\Top}_{\eta}(\pi, G^1,\cP^1, \xi).
\] 
\end{Proposition}
\begin{proof}By Proposition \ref{prop: propagator unique up to concordance}, one can find a concordance 
from $(G^0,\widetilde{\cP}^{0})$ to $(G^1,\widetilde{\cP}^{1})$, given by an observer-propagator system $(\widetilde{G},\widetilde{\cP})$ for the bundle
\[
\tilde{\pi}: \widetilde{E}_D= [0,1]\times E_D\to [0,1]\times B.
\]
By Proposition \ref{prop: cobordism invariance, propagator}, we see that 
\[
I^{\Top}_{\eta}(\pi, G^0,\cP^0, \xi)=I^{\Top}_{\eta}(\pi, G^1,\cP^1, \xi).
\] 
\end{proof}

To prove the independence with $\xi$, we need the following lemma.

\begin{Lemma}\label{lem: Čech cocycle concordance} Let $\xi^0,\xi^{1}\in C^{3q}_{\cech}(\prod_{q}S^3)$ be two standard cocycles. Then there exists a cocycle 
\[\widetilde{\xi}\in C^{3q}_{\cech}([0,1]\times \prod_{q}S^3)\]
that pulls back to $\xi^{i}$ along $\{i\}\times \prod_{q}S^3$ and satisfies 
\[g^{*}(\widetilde{\xi})=\operatorname{sign}(g)\cdot \widetilde{\xi},\quad\forall g\in \Sigma_{q}.
\]
\end{Lemma}
\begin{proof}
 we let $\xi^1=\xi^{0}+d\alpha$ for some $\alpha\in C^{3q-1}_{\cech}(\prod_{q}S^3)$. Then we can find $\tilde{\alpha}\in C^{3q-1}_{\cech}([0,1]\times \prod_{q}S^3)$ that pulls back to $0$ along $\{0\}\times  \prod_{q}S^3$ and pulls back to $\alpha$ along $\{1\}\times  \prod_{q}S^3$. Let \[p:[0,1]\times \prod_{q}S^3\to \prod_{q}S^3\] be the projection map. Let $\widetilde{\xi}'=p^{*}(\xi^{0})+\widetilde{\alpha}$ and let 
\[\widetilde{\xi}=\frac{1}{q!}\sum_{g\in \Sigma_{q}}(\operatorname{sign}(g)\cdot g^{*}(\widetilde
{\xi}')).\]
Then $\widetilde{\xi}$ has the desired properties.
\end{proof}

\begin{Proposition}\label{prop: independence with xi}  Let $\pi:E_D\to B$ be a boundary trivialized, framed disk bundle over a manifold $B$ of dimension $m+\frac{n}{2}$ and let $\eta\in \widetilde{\cG}^{m,n}$ be a cocycle. Then for any observer-propagator system $(G,\cP)$ and any two choices of standard cocycles $\xi^0$ and $\xi^1$, we have  
\[
I^{\Top}_{\eta}(\pi, G,\cP, \xi^0)=I^{\Top}_{\eta}(\pi, G,\cP, \xi^1).
\] 
\end{Proposition}
\begin{proof} Let $q=\frac{3n-2m}{2}$. By Lemma \ref{lem: Čech cocycle concordance}, we can find a cocycle 
\[
\widetilde{\xi}\in C^{3q}([0,1]\times \prod_{q}S^3)
\]
that pulls back to $\xi^0$ and $\xi^1$ on its two boundary components. 
We let $\eta= \sum_{1\leq i\leq k}a_{i}(\Gamma_{i},\rho_{i},o_{i})$. We use $\rho_{i}$ to identify $V(\Gamma_{i})$ with $\bfn/S_{i}$ for some $\cS_{i}$. For each $i$, we pick a bijection $E(\Gamma_{i})\to \{1,\cdots, q\}$ that represents $o_{i}$. This induces a map 
\[
[0,1]\times X(\Gamma_{i})\to [0,1]\times \prod_{q}S^3
\]
We let $\widetilde{\xi}_{o_i}$ be the pullback of $\widetilde{\xi}$ under this map.

For each $i$, we consider the map \[
\widehat{P}(\Gamma_{i},\rho_{i}): \bar{C}^{G_{\cS_{i}}}_{\bfn/\cS_{i}}([0,1]\times E^+)=[0,1]\times \bar{C}^{G_{\cS_{i}}}_{\bfn/\cS_{i}}(E^+)\xrightarrow{\id\times P(\Gamma_{i},\rho_{i})} [0,1]\times X(\Gamma_{i})
\]
Then similar to the proof of Proposition \ref{prop: cobordism difference}, we have 
\begin{equation}\label{eq: Čech difference}
\begin{split}
I^{\Top}_{(\Gamma_{i},\rho_{i},o_{i})}(\pi, G,\cP, \xi^{1})-&
I^{\Top}_{(\Gamma_{i},\rho_{i},o_{i})}(\pi, G,\cP, \xi^{0})=\\ &\frac{1}{d^q}\sum_{e\in E'(\Gamma_{i})} 
\langle \widehat{P}(\Gamma_{i}, \rho_{i})^{*}(\widetilde{\xi}_{o_{i}}), [\partial_{e}\bar{C}^{G_{\cS_{i}}}_{\bfn/\cS_{i}}(\tilde{E}^+)]\rangle
    \end{split}
\end{equation}

Since $\eta$ is a cocycle, terms on the right hand side of (\ref{eq: Čech difference}) cancel when we do linear combination with coefficients $a_{i}$. This finishes the proof.
\end{proof}

By Proposition \ref{prop: independence with P and G} and Proposition \ref{prop: independence with xi}, for any cocycle $\eta$, the quantity $I^{\Top}_{\eta}(\pi, G,\cP, \xi)$
is independent with the choice of $(G,\cP,\xi)$. So we will simply write $I^{\Top}_{\eta}(\pi)$.
\begin{Proposition}\label{prop: independence with X} For any cocyle $\eta\in \widetilde{\cG}^{m,n}$, the quantity $I^{\Top}_{\eta}(\pi)$ is independent with the choice of the sphere sequence (\ref{eq: sequence of maps on S3}).
\end{Proposition}
\begin{proof} We consider two sphere sequences $X(-)=\{h_{i}:X(i)\to X(i+1)\}_{i\geq 1}$ and $X'(-)=\{h'_{i}:X(i)\to X(i+1)\}_{i\geq 1}$. Given  a bundle $\pi: E\to B$, we pick an observer system $G$, a propagator system
\[
\cP=\{P^{\cS}_{[i],[j]}: \bar{C}^{G_{S}}_{\bfn/\cS}(E^+)\to X(l)\}
\]
for $X(-)$, and a propagator system 
\[
\cP'=\{P^{\cS,'}_{[i],[j]}: \bar{C}^{G_{S}}_{\bfn/\cS}(E^+)\to X(l')\}
\]
for $X'(-)$. We can find a nonzero degree maps $\rho: X(l)\to S^3$ and $\rho':X(l')\to S^3$ such that the maps 
\[
X(1)\xrightarrow{h_{l}\circ\cdots\circ h_1} X(l)\xrightarrow{\rho} S^3 \quad \text{and}\quad
X(1)\xrightarrow{h'_{l'}\circ\cdots\circ h'_1} X(l')\xrightarrow{\rho'} S^3
\]
are homotopic. We let $\{\bar{h}_{1}^{t}:X(1)\to S^3\}$ be such a homotopy and we consider a family $\{X^{t}(-)\}_{0\leq t\leq 1}$ of sphere sequences,  given by  
\[
X(1)\xrightarrow{\bar{h}^t} X(2)\xrightarrow{\bar{h}_2} X(3)\xrightarrow{\bar{h}_3} \cdots.
\]
Here $\bar{h}_2, \bar{h}_3,\cdots $ can be any maps of large degrees. By repeating the proof of Proposition \ref{prop: propagator unique up to concordance}, we may construct a family of maps 
\[\cP^{t}=\{P^{\cS,t}_{[i],[j]}: \bar{C}^{G_{\cS}}_{\bfn/\cS}(E^+)\to X(l_{\cS})\}\]
such that $P^{\cS,t}_{[i],[j]}$ depends continuously on $t$ and such that $\cP^{t}$ is a propagator system for the sphere sequence $X^{t}(-)$. The same cobordism argument as Proposition \ref{prop: cobordism invariance, propagator} shows that for any  standard Čech cocycle $\xi$, the quantity $I^{\Top}_{\eta}(\pi,G,\cP^{t},\xi)$ is independent with $t$.  By the definition of $\bar{h}^{t}$, we see that 
$I^{\Top}_{\eta}(\pi,G,\cP^{0},\xi)=I^{\Top}_{\eta}(\pi,G,\cP,\xi)$
and 
$I^{\Top}_{\eta}(\pi,G,\cP^{1},\xi)=I^{\Top}_{\eta}(\pi,G,\cP',\xi)$. This finishes the proof.
\end{proof}

To this end, we have proved the invariance of $I^{\Top}_{\eta}(\pi,G,\cP,\xi)$ under different choices of $\pi,G,\cP$ and sphere sequence. This justify the notation  $I^{\Top}_{\eta}(\pi)$. 

It remains to prove the invariance of  $I^{\Top}_{\eta}(\pi)$ under difference representatives of $\eta$ in the same cohomology class. We start with some preparations. We fix a cocycle $\alpha\in C^{3}_{\cech}(S^3)$  
 that represents the standard generator of $H^{3}(S^3)$. And we fix a cochain   $\beta\in C^{5}_{\cech}(S^3\times S^3)$ such that \[d\beta=\operatorname{pj}_{1}(\alpha)\wedge \operatorname{pj}_{2}(\alpha)+ \operatorname{pj}_{2}(\alpha)\wedge \operatorname{pj}_{1}(\alpha).\]  
Here $\operatorname{pj}_{i}:S^3\times S^3\to S^3$ denotes projection to the $i$-th component. Fix an integer $q\geq 1$. For any $1\leq i\neq j\leq q$, we let $\alpha_{i}\in C^{3}_{\cech}(\prod_{q} S^3)$ be the  pull back of $\alpha$ along the projection to the $i$-th component and let $\beta_{i,j}\in C^{5}_{\cech}(\prod_{q} S^3)$ be the pull back of $\beta$ along the $(i,j)$-th component. For any element $g$ of the permutation group $\Sigma_{q}$, we let 
\[\xi(g)=\alpha_{g(1)}\wedge  \cdots \wedge \alpha_{g(q)}\in C^{3q}_{\cech}(\prod_{q} S^3). \]
and we let 
\[
\xi_{q}=\frac{1}{q!}\sum_{g\in \Sigma_{q}}\operatorname{sign}(g)\cdot \xi(g)
\]
It's straightforward to see that $\xi_{q}$ is a standard cocycle. 

For each $1\leq i\leq q$, we consider the subset $\Sigma_{q,i}=\{g\in\Sigma_{q}\mid g(i)=1\}$. Then we define 
\[
\xi_{q,i}=\frac{1}{(q-1)!}\sum_{g\in \Sigma_{q,i}}\operatorname{sign}(g)\cdot \xi(g).
\]

Although $\xi_{q}$ and $\xi_{q,i}$ are not equal, their difference can be expressed as 
\[\xi_{q}-\xi_{q,i}=d\gamma_{q,i} \]

Here $\gamma_{q,i}\in C^{3q-1}_{\cech}(\prod_{q} S^3)$ is a polynomial in $\alpha_{*}$ and $\beta_{*,*}$. For our discussion, the specific form of $\gamma_{q,i}$ is not important. We just need the symmetric property
\begin{equation}\label{eq: gamma symmetric}
g^{*}(\gamma_{q,i})=\operatorname{sign}(g)\cdot \gamma_{q,i},\quad \forall g\in \Sigma_{q,i}.   \end{equation}

\begin{Proposition}\label{prop: invariance under eta} Suppose $[\eta]=[\eta']\in H^{m,n}(\widetilde{\cG})$. Then 
$I^{\Top}_{\eta}(\pi)=I^{\Top}_{\eta'}(\pi)$
for any framed, boundary trivialized disk bundle $\pi:E_{D}\to B$ over a closed manifold $B$ of dimension $m+\frac{n}{2}$ .
\end{Proposition}

\begin{proof}
By linearity, it suffices to show that $I^{\Top}_{\eta}(\pi)=0$ if $\eta=\widetilde{\delta}^{m-1,n}(\Gamma,\rho,o)$ for some decorated, oriented graph $(\Gamma,\rho,o)$. We let \[p=|V(\Gamma)|=n-m+1,\quad q=|E(\Gamma)|=\frac{3n}{2}-m+1.\] We use $\rho$ to fix an identification between $V(\gamma)$ and $\bfn/\cS$ for some $\cS$. And we fix a bijection $i:E(\gamma)\to \{1,\cdots, q\}$ that is compatible with $o$. This gives a diffeomorphism $X(\Gamma)\cong \prod_{q} S^3.$
Using this diffeomorphism, we pull back $\xi_{q}$ and obtain $\xi_{o}\in C^{3q}_{\cech}(X(\Gamma))$. For each $e\in E(\Gamma)$, we also pull back $\xi_{q,i(e)}$ and $\gamma_{q,i(e)}$ to $\xi_{o,e}$ and $\gamma_{o,e}$ respectively. Then we have \begin{equation}\label{eq: dgamma}
\xi_{o}-\xi_{o,e}=d\gamma_{o,e}.    
\end{equation}

We fix an observer-propagator system $(G,\cP)$ for the bundle $\pi$. And we consider the map 
\[
P(\Gamma, \rho):\bar{C}^{G_{\cS}}_{\bfn/\cS}(E^+)\to X(\Gamma).
\]
Note that $\dim(\bar{C}^{G_{\cS}}_{\bfn/\cS}(E^+))=3q+1=\dim(X(\Gamma))+1$.
By the decomposition \[\partial \bar{C}^{G_{\cS}}_{\bfn/\cS}(E^+)=\bigcup\limits_{A\subset \bfn/\cS, |A|\geq 2} \partial_{A}\bar{C}^{\tilde{G}_{\cS}}_{\bfn/\cS}(E^+),\]
we see that 
\begin{equation}\label{eq: Čech boundary vanishing 1}
\sum_{A\subset \bfn/\cS, |A|\geq 2}\langle P(\Gamma, \rho)^{*}(\xi_{o}), [\partial_{A}\bar{C}^{G_{\cS}}_{\bfn/\cS}(E^+)]\rangle=0.    
\end{equation}
By Proposition \ref{prop: boundary vanishing |A|>2} and Proposition \ref{prop: boundary vanishing |A|=2}, only those faces that correspond to some $e\in E'(\Gamma)$ can have nonzero contributions to (\ref{eq: Čech boundary vanishing 1}). Therefore, we get \begin{equation}\label{eq: Čech boundary vanishing 2}
\sum_{e\in E'(\Gamma)}\langle P(\Gamma, \rho)^{*}(\xi_{o}), [\partial_{e}\bar{C}^{G_{\cS}}_{\bfn/\cS}(E^+)]\rangle=0.
\end{equation} 
Next, we note that (\ref{eq: dgamma}) implies 
\begin{equation}\label{eq: Čech boundary vanishing 5}
\langle P(\Gamma, \rho)^{*}(\xi_{o}-\xi_{o,e}),[\partial_{e}\bar{C}^{G_{\cS}}_{\bfn/\cS}(E^+)]\rangle= \langle P(\Gamma, \rho)^{*}(\gamma_{o,e}),[\partial\partial_{e}\bar{C}^{G_{\cS}}_{\bfn/\cS}(E^+)]\rangle    
\end{equation}

By the decomposition 
\[\partial \partial_{e} \bar{C}^{G_{\cS}}_{\bfn/\cS})(E^+)=\bigcup\limits_{A\subset \bfn/(\cS\cdot e), |A|\geq 2} \partial_{A}\partial_{e} \bar{C}^{G_{\cS}}_{\bfn/\cS}(E^+),\]
we obtain 
\begin{equation}\label{eq: Čech boundary vanishing 4}
\langle P(\Gamma, \rho)^{*}(\gamma_{o,e}),[\partial\partial_{e}\bar{C}^{G_{\cS}}_{\bfn/\cS}(E^+)]\rangle=\sum_{A\subset \bfn/(\cS\cdot e), |A|\geq 2}\langle P(\Gamma, \rho)^{*}(\gamma_{o,e}), [\partial_{A}\partial_{e} \bar{C}^{G_{\cS}}_{\bfn/\cS}(E^+)]\rangle
\end{equation}
By (\ref{eq: gamma symmetric}),  $\gamma_{o,e}$ is symmetric under any permutations of $E(\Gamma)$ that fixes $e$. Therefore, we can repeat the same argument in the proof of Proposition \ref{prop: boundary vanishing |A|>2} and Proposition \ref{prop: boundary vanishing |A|=2} to conclude that 
\[
\langle P(\Gamma, \rho)^{*}(\gamma_{o,e}), \partial_{A}\partial_{e} \bar{C}^{G_{\cS}}_{\bfn/\cS}\rangle=0
\]
unless $A$ corresponds to some $e'\in E'(\Gamma/e)$. So we can rewrite (\ref{eq: Čech boundary vanishing 4}) as  
\begin{equation}\label{eq: Čech boundary vanishing 6}
\langle P(\Gamma, \rho)^{*}(\gamma_{o,e}),[\partial\partial_{e}\bar{C}^{G_{\cS}}_{\bfn/\cS}(\tilde{E}^+)]\rangle=\sum_{e'\in E'(\Gamma/e)}\langle P(\Gamma, \rho)^{*}(\gamma_{o,e}), [\partial_{e'}\partial_e \bar{C}^{G_{\cS}}_{\bfn/\cS}]\rangle
\end{equation}

Combining (\ref{eq: Čech boundary vanishing 5}) and (\ref{eq: Čech boundary vanishing 6}), we get 
\begin{equation}\label{eq: Čech boundary vanishing 8}
\sum_{e\in E'(\Gamma)}\langle P(\Gamma,\rho)^{*}(\xi_{o}-\xi_{o,e}),[\partial_{e}\bar{C}^{G_{\cS}}_{\bfn/\cS}(E^+)]\rangle=\sum_{(e,e')\in B}\langle P(\Gamma,\rho)^{*}(\gamma_{o,e}),[\partial_{e'}\partial_e\bar{C}^{G_{\cS}}_{\bfn/\cS}(E^+)]\rangle
\end{equation}
Here 
 \[B=\{(e,e')\in E(\Gamma)\times E(\Gamma)\mid e\neq e'\text{ and } \Gamma/(e\cup e')\text{ is admissible}\}.\]
Notice that both
$\partial_{e'}\partial_e \bar{C}^{G_{\cS}}_{\bfn/\cS}(E^{+})$ and
$\partial_{e}\partial_{e'} \bar{C}^{G_{\cS}}_{\bfn/\cS}(E^{+})$
 can be identified with 
$$
\partial_{\{e.e'\}}\bar{C}^{G_{\cS}}_{\bfn/\cS}(E^{+})=
\bar{C}^{G_{\cS}}_{\bfn/(\cS\cdot e\cup e')}(E^{+})
\times \bar{C}_e^\ast(\mathbb{R}^4)
\times  \bar{C}_{e'}^\ast(\mathbb{R}^4)
$$
according to \eqref{eq_partial_S_product_decomposition},
but their orientations are different. Let 
$\tilde{g}_{e,e'}$ be the diffeomorphism
on $\partial_{\{e.e'\}}\bar{C}^{G_{\cS}}_{\bfn/\cS}(E^{+})$ that switches the components 
$\bar{C}_e^\ast(\mathbb{R}^4)$ and
 $\bar{C}_{e'}^\ast(\mathbb{R}^4)$. Then 
 $\tilde{g}_{e,e'}$ can be viewed as an orientation preserving diffeomorphism from $\partial_{e'}\partial_e \bar{C}^{G_{\cS}}_{\bfn/\cS}(E^{+})$ to
$\partial_{e}\partial_{e'} \bar{C}^{G_{\cS}}_{\bfn/\cS}(E^{+})$ and
 we have
the following commutative diagram
\begin{equation*}
\begin{tikzcd}
\partial_{e'}\partial_e \bar{C}^{G_{\cS}}_{\bfn/\cS}(E^{+})\ar[rr,"P(\Gamma{,}\rho)"]\arrow[d,"\tilde{g}_{e,e'}"] & & X(\Gamma)\arrow[d,"g_{e,e'}"]\\
\partial_{e}\partial_{e'}\bar{C}^{G_{\cS}}_{\bfn/\cS}(E^{+})\ar[rr,"P(\Gamma{,}\rho)"]&  & X(\Gamma)\\
\end{tikzcd}
\end{equation*}
where $g_{e,e'}$ is the permutation that switches the components $e$ and $e'$.  Since $g_{e,e'}$ pulls back $\gamma_{o,e'}$ to $-\gamma_{o,e}$, we see that 
\[
\langle P(\Gamma,\rho)^{*}(\gamma_{o,e}),[\partial_{e'}\partial_{e}\bar{C}^{G_{\cS}}_{\bfn/\cS}(E^+)]\rangle=-
\langle P(\Gamma,\rho)^{*}(\gamma_{o,e'}),[\partial_{e}\partial_{e'}\bar{C}^{G_{\cS}}_{\bfn/\cS}(E^+)]\rangle
\]
Therefore, the right hand side of (\ref{eq: Čech boundary vanishing 8}) cancels each other and we get 
\[\sum_{e\in E'(\Gamma)}\langle P(\Gamma,\rho)^{*}(\xi_{o}),[\partial_{e}\bar{C}^{G_{\cS}}_{\bfn/\cS}(E^+)]\rangle=\sum_{e\in E'(\Gamma)}\langle P(\Gamma,\rho)^{*}(\xi_{o,e}),[\partial_{e}\bar{C}^{G_{\cS}}_{\bfn/\cS}(E^+)]\rangle
\]
So (\ref{eq: Čech boundary vanishing 2}) gives 
\[\sum_{e\in E'(\Gamma)}\langle P(\Gamma,\rho)^{*}(\xi_{o,e}),[\partial_{e}\bar{C}^{G_{\cS}}_{\bfn/\cS}(E^+)]\rangle=0.
\]
The compatibility condition (\ref{eq boundary condition: i,j,0 different}) implies 
\[
I^{\Top}_{(\Gamma/e,\rho/e,o/e)}(\pi)=\frac{1}{d^q}\langle P(\Gamma,\rho)^{*}(\xi_{o,e}), [\partial_{e}\bar{C}^{G_{\cS}}_{\bfn/\cS}(E^+)]\rangle.
\]
So we get  
\[I^{\Top}_{\eta}(\pi)=\sum_{e\in E'(\Gamma)} I^{\Top}_{(\Gamma/e,\rho/e,o/e)}(\pi)=0.
\]
This finishes the proof.
\end{proof}

\begin{proof}[Proof of Theorem \ref{thm: I-top welldefined}] Independence with the choice of a observer-propagator system is proved in Proposition \ref{prop: independence with P and G}. Independence with the Čech cocycle is proved in Proposition \ref{prop: independence with xi}. Independence with  the sphere sequence is proved in Proposition \ref{prop: independence with X}. By Proposition \ref{prop: invariance under eta}, different representatives $\eta$ of the same homology class give the same invariant $I^{\Top}_{[\eta]}(\pi)$. And by Proposition \ref{prop: cobordism invariance, propagator}, the invariant $I^{\Top}_{[\eta]}(\pi)$ is unchanged under any framed cobordism.
\end{proof}

\subsection{Extension to formally smooth bundle over spheres}
We end this section  by extending the definition of $I^{\Top}_{[\eta]}$ to formally smooth disk bundles over spheres.

\begin{Lemma}
\label{cor: cobordism invariance} Let $[\eta]\in H^{m,n}(\widetilde{\cG})$ be a cohomology class and let $\pi: E_{D}\to S^{m+n}$ be a framed disk bundle. Let  $f^{*}(\pi): E'_{D}\to S^{m+n}$ be the corresponding bundle obtained by pulling back $\pi$ via a map $f:S^{m+n}\to S^{m+n}$. Then we have 
\begin{equation}\label{eq: I-top pullback}
I^{\Top}_{[\eta]}(f^{*}(\pi))=\deg(f)\cdot I_{[\eta]}^{\Top}(\pi).
\end{equation}
\end{Lemma}
\begin{proof}
Without loss of generality, we may assume $\deg(f)\geq 0$. By attaching 1-handles, one may find a cobordism $W$ from $\bigsqcup_{\deg(f)}S^{m+n}$ to $S^{m+n}$ and a map $\widetilde{f}: W\to S^{m+n}$ that restricts to the identity on each incoming boundary component and restricts to $f$ on the outgoing boundary. We pull back $\pi$ via $\widetilde{f}$ and apply Theorem \ref{thm: I-top welldefined} to finish the proof.
\end{proof}

\begin{Lemma}\label{lem: linear structure to framing}
Let $\pi: E_{D}\to S^{k}$ ($k>0$) be a disk bundle equipped with a linear structure $\theta$ on its vertical tangent microbundle. Then there exists a nonzero degree map $f:S^{k}\to S^{k}$ such that the pull-back linear structure $f^{*}(\theta)$ can be promoted to a framing $\widetilde{\theta}$  on $f^{*}(\pi)$. Furthermore, given any two such maps $f_{1}, f_{2}:S^{k}\to S^{k}$, let $\widetilde{\theta}_{i}$ be the corresponding vertical framing on $f^{*}_{i}(\pi)$. Then there exists a nonzero-degree map $g_{1},g_{2}: S^{k}\to S^{k}$ such that 
\begin{equation}\label{eq: deg}
\deg(f_{1})\cdot\deg(g_{1})=\deg(f_{2})\cdot\deg(g_{2})
\end{equation}
and the framed bundles \[((g_{1}\circ f_{1})^{*}(\pi), g_{1}^{*}(\widetilde{\theta}_{1})),\quad ((g_{2}\circ f_{2})^{*}(\pi), g_{2}^{*}(\widetilde{\theta}_{2}))
\] are framed cobordant to each other. 
\end{Lemma}
\begin{proof}
The linear structure $\theta$ is an isomorphism between $\cT^{v}E$ and a vector bundle $V\to E$. Since $\theta$ is standard in $E\setminus E_{D}$, we have a trivilization $\rho: V|_{E\setminus E_{D}}\to \bR^{4}\times (E\setminus E_{D})$ that coincide (under $\theta$) with the standard linear structure on $E\setminus E_{D}$. And a  framing is an extension of $\rho$ to the whole bundle $V$.

We fix a base point $b\in S^k$ and use $F_{b}\cong \bR^{4}$ to denote the fiber of $E$ over $b$. Then $\theta|_{F_{b}}$ corresponds to a reduction of the structure group of $\cT(\bR^{4})$ from $\Top(4)$ to $O(4)$. And this reduction is standard on $\bR^{4}\setminus D^{4}$. Such reduction is given by a map $(\bR^{4},\bR^{4}\backslash D^4)\to (\Top(4)/O(4),*)$. Since $\pi_{4}(\Top(4)/O(4))=0$, we see that $\theta|_{F_{b}}$ is homotopic to the standard linear structure on $F_{b}$. This implies the trivilization $\rho$ can be extended to $F_{b}$. 

The next step is to extend $\rho$ from $(E\setminus E_{D})\cup F_{b}$ to $E$. Since the latter is obtained from the former by attaching a single $(k+4)$-cell, the obstruction for such extension is an element in the group \[H^{k+4}(E,(E\setminus E_{D})\cup F_{b}; \pi_{k+3}(O(4));\mathbb{Z}).\]
Since $\pi_{k+3}(O(4))$ is torsion, such obstruction class can be eliminated by pulling back with a map $g:S^{k}\to S^{k}$ with a suitable degree. This proves the first assertion.

Next, we assume we have two pairs $(f_{i},\widetilde{\theta}_{i})$ for $i=1,2$. We take suitable $g_{i}:S^{k}\to S^{k}$ such that (\ref{eq: deg}) is satisfied. Then we have a homotopy $h:[0,1]\times S^{k}  \to S^{k}$ from $g_{1}\circ f_{1}$ to $g_{2}\circ f_{2}$. By pulling back $(\pi,\theta)$ via $h$, we get a disk bundle over $I\times S^k$ equipped with vertical linear structure $\widetilde{\theta}$. This provides a cobordism between $(g_{1}\circ f_{1})^{*}(\theta)$ and $(g_{2}\circ f_{2})^{*}(\theta)$. The restriction of $\widetilde{\theta}$ to $(\partial I\times E)\cup (I\times F_{b})$ is already a vertical framing. The obstruction to extend this framing to $I\times E$ is given by an element in 
\[
H^{k+5}(I\times E, (\partial I\times E)\cup (I\times F_{b});\pi_{k+4}(O(4))).
\]
Again, since $\pi_{k+4}(O(4))$ is torsion, we can eliminate this obstruction class by replacing $g_{i}$ with $g\circ g_{i}$ for a suitable $g:S^k\to S^k$ of nonzero degree.
\end{proof}

\begin{Definition} Let $[\eta]\in H^{m,n}(\widetilde{\cG})$ be a cohomology class in the decorated graph complex. Let $k=\frac{n}{2}+m$ and let $\pi:E\to S^{k}$ be a topological bundle equipped formal smooth structure $\theta$. We pick a map $f:S^{k}\to S^{k}$ of nonzero degree such that $f^{*}(\theta)$ can be lifted to a framing on $f^{*}(\pi)$. Then we define
\[
I_{[\eta]}^{\Top}(\pi):=\frac{1}{\deg(f)}I^{\Top}_{[\eta]}(f^{*}(\pi)).
\]
This is well-defined by Corollary \ref{cor: cobordism invariance} and Lemma \ref{lem: linear structure to framing}.
\end{Definition}

\begin{Lemma}\label{lemma: concordance invaraince}
Let $[\eta]\in H^{m,n}(\widetilde{\mathcal{G}})$ be a cohomology class. Let $\widetilde{B}$ be the $(p+q)$-th punctured $S^{k+1}$, treated as an oriented cobordism from $\sqcup_{p}S^{k}$ and $\sqcup_{q}S^{k}$. Let $\widetilde{\pi}:\widetilde{E}_{D}\to \widetilde{B}$ be a formally smooth disk bundle, equipped with a vertical linear structure $\xi$. We use $\pi_{i}:E^{i}_{D}\to S^2$ (resp. $\pi'_{j}:E^{j}_{D}\to S^2$) to denote the restriction of $\widetilde{\pi}$ to the $i$-th component (resp. the $j$-th outgoing component). Then one has 
\begin{equation}\label{eq: concordant invariance}
\sum_{1\leq i\leq p}I_{[\eta]}^{\Top}(\pi_{i})=\sum_{1\leq j\leq q}I_{[\eta]}^{\Top}(\pi'_{j}).    
\end{equation}
In particular, $I^{\Top}$ is invariant under concordance of formal smooth disk bundles over spheres.
 \end{Lemma}
\begin{proof}
The only obstruction of lifting the vertical linear structure on $\widetilde{\pi}$ to a vertical framing is an element in the group
\[
H^{k+5}(\widetilde{E}_{D},\partial^{v}\widetilde{E}_{D};\pi_{k+4}(O(4)))\cong H^{k+1}(\widetilde{B};\pi_{k+4}(O(4)))=0.
\]
So we can lift $\xi$ to a vertical framing. So $\sqcup_{1\leq i\leq p} \pi_{i}$ to $\sqcup_{1\leq j\leq q} \pi'_{j}$ are framing cobordant to each other. So (\ref{eq: concordant invariance}) follows from the cobordism invariance of $I^{\Top}$ (see Theorem \ref{thm: I-top welldefined}).
\end{proof}

\section{Equivalence of two definitions for smooth families}\label{section: equivalence of definitions}

In \cite{Wa18,Wa21}, Watanabe defined a version of Kontsevich's characteristic classes $I(-)$ for smooth disk bundles equipped with a framing. The purpose of this section is to prove that the invariants $I(-)$ and $I^{\Top}(-)$ coincide with each other.

Throughout this section, we assume the following setting:  
let $D^{4}\to E_{D}\xrightarrow{\pi} B$ be a smooth disk bundle over a closed, smooth manifold $B$ of dimension $m+n$, equipped with  a framing on its vertical tangent bundle. We fix a Riemannian metric on $E_{D}$ and use the exponential map to define a framing on the vertical tangent microbundle. Since the bundle $\pi$ is smooth, for any based, finite set $A$, the configuration space $C_{A}(E^{+})$ has a canonical compactification $\bar{C}_{A}(E^{+})$ and it coincide with our definition $\bar{C}^{G}_{A}(E^{+})$, as long as the observer system $G$ is smooth. Actually, in this case, the linear structure $\theta(g_{S})$ that appears in (\ref{eq_def_iota_family}) is compatible with the smooth structure. So the space $\Bl_{\Delta_S}^{\theta(g_S)} E^n$ is just the blow up of $E^n$ along the smooth submanifold $\Delta_S$. Furthermore, given any two non-base elements $i\neq j\in A$, the map 
\[
C_{A}(E^{+})\to C_{\mathbf{2}}(E^{+}),\quad (x_{a})_{a\in A}\mapsto (\infty,x_{i},x_{j})
\]
extends to a smooth map 
\begin{equation}\label{eq: configuration space forgetful map}
\phi_{i,j}: \bar{C}_{A}(E^{+})\to \bar{C}_{\mathbf{2}}(E^{+}).
\end{equation}

Given any decorated, oriented admissible $(\Gamma,\rho,o)$ with $p=n-m$ vertices and $q=\frac{3n}{2}-m$ edges. The the decoration $\rho$ induces a total order on $V(\Gamma)$. We use this order to fix the identification $V(\Gamma)\cong \{1,\cdots, p\}$. We also fix an order on $E(\Gamma)\leftrightarrow \{1,\cdots, q\}$ that is compatible with $o$. For each $1\leq l\leq q$, we let $i(l), j(l)$ be the start point and the end point of the $l$-th edge $e_{l}$. Then we define the map 
\begin{equation}\label{eq: phi-gamma}
\phi_{(\Gamma,\rho,o)}:=\prod_{1\leq k\leq q}\phi_{i(e_{k}),j(e_{k})}: \bar{C}_{\bfp}(E^+)\to \prod_{q}\bar{C}_{\mathbf{2}}(E^+).
\end{equation}
By Lemma \ref{lem: propagator n=2 boundary}, 
there is a map 
\[P^{\partial}: \partial\bar{C}_{\mathbf{2}}(E^{+})\to S^3=X(1)\]
whose restriction to $\partial_{A} \bar{C}_{\mathbf{2}}(E^{+})$ equals the composition 
\[
\partial_{A} \bar{C}_{\mathbf{2}}(E^{+})\xrightarrow{p_2} \bar{C}^{*}_{A}(\bR^4)\xrightarrow{\psi_{k,l}} S^3
\]
where
\[
(k,l)=\begin{cases} (1,2)\text{ when }A=\{0,1,2\},\ \{1,2\}\\
(0,1)\text{ when }A=\{0,1\}\\
(2,0)\text{ when }A=\{2,0\}\\
\end{cases}.
\]
Let $\alpha=\Vol\in \Omega^{3}(S^3)$ be the standard volumn form. It is proved in \cite[Lemma 2.10]{Wa21} that $P^{\partial,*}(\alpha)$ can be extended to a closed form $\omega\in \Omega^{3}(\bar{C}_{\mathbf{2}}(E^{+}))$. We let 
\[\omega(q):=\operatorname{pj}_{1}^\ast(\omega)\wedge\cdots \wedge \operatorname{pj}_{q}^\ast(\omega)\in \Omega^{3q}(\prod_{q}\bar{C}_{\mathbf{2}}(E^{+})).\]
Now we consider a homology class  $[\tilde{\eta}]\in H^{m,n}(\widetilde{\cG})$, represented by a cocycle 
\[\widetilde{\eta}=\sum_{l=1}^{k}a_{i}(\Gamma_{l},\rho_{l},o_{l}),\] where $(\Gamma_{l},\rho_{l},o_{l})$ is a decorated, oriented admissible graph with $p$ vertices and $q$ edges. We let $[\eta]=f_{*}([\widetilde{\eta}])$ be it image under the forgetful map $f_*$ (see (\ref{eq: forgetful})). In \cite{Wa21}, Watanabe defined the Kontsevich integral as follows 
\[
I_{[\eta]}(\pi):=\sum^{k}_{i=1}a_{i}\cdot \int_{\bar{C}_{n}(E^+)}\phi^{*}_{(\Gamma_{i},\rho_{i},o_{i})}(\omega(q)). 
\]
\begin{Lemma}\label{lem: invariance of differential form} Let $\alpha_{1},\cdots, \alpha_{q}\in \Omega^{3}(S^3)$ be differential forms that are invariant under the antipodal map on $S^3$ and satisfy $\int_{S^3}\alpha=1$. Let $\omega_{i}\in \Omega^{3}(\bar{C}_{\mathbf{2}}(E^+))$ be any closed 3-form that extends  $P^*(\alpha_{i})$. Consider the differential form
\begin{equation}\label{eq: omega'q}
\omega'(q):=\frac{1}{q!}\sum_{\sigma\in\Sigma_{q} }\operatorname{pj}_{1}^\ast(\omega_{\sigma(1)})\wedge\cdots \wedge \operatorname{pj}_{q}^\ast(\omega_{\sigma(q)})).
\end{equation}
Then one has 
\begin{equation*}
I_{[\eta]}(\pi)=  \sum^{k}_{i=1}a_{i}\cdot \int_{\bar{C}_{n}(E^+)}\phi^{*}_{(\Gamma_{i},\rho_{i},o_{i})}(\omega'(q)). 
\end{equation*}
\end{Lemma}
\begin{proof} This is a mild generalization of \cite[Theorem 2.15 (2)]{Wa21}, which proves the case $\alpha_i=\alpha$ for all $i$. (See \cite[Page 105, second paragraph]{Wa21}). The only special feature used in the proof is the invaraiance of $\alpha_{i}$ under the antipodal map so it  can be generalized verbatimly in our situation.

For completeness, we sketch the argument here: By our assumption, there exists closed 3-form \[\widetilde{\alpha}_{i}\in \Omega^3([0,1]\times S^3)\] that restricts to $\alpha$ on $\{0\}\times S^3$ and restricts to $\alpha_{i}$ on on $\{1\}\times S^3$. By taking average, we may further assume that $\widetilde{\alpha}_{i}$ is invariant under the antipodal map $[0,1]\times S^3$ defined by $(t,x)\mapsto (t,-x)$. Because the inclusion \[\partial ([0,1]\times \bar{C}_{\mathbf{2}}(E^+))\hookrightarrow [0,1]\times\bar{C}_{\mathbf{2}}( E^+)\] induces a surjection on $H^3(-)$, there exists a closed 3-form  $\widetilde{\omega}_{i}\in \Omega^{3}([0,1]\times \bar{C}_{\mathbf{2}}(E^{+}))$ that satisfies the boundary condition 
\[
\begin{split}
\widetilde{\omega}_{i}|_{\{0\}\times \bar{C}_{\mathbf{2}}(E^{+})}&=\omega,\\
\widetilde{\omega}_{i}|_{\{1\}\times \bar{C}_{\mathbf{2}}(E^{+})}&=\omega_{i}\\
\widetilde{\omega}_{i}|_{[0,1]\times \partial\bar{C}_{\mathbf{2}}(E^{+})}&=(\id\times P^{\partial})^{*}(\widetilde{\alpha}_{i}).
\end{split}
\]
We let 
\[
\widetilde{\omega}(q):=\frac{1}{q!}\sum_{\sigma\in\Sigma_{q} }\operatorname{pj}_{1}^\ast(\widetilde{\omega}_{\sigma(1)})\wedge\cdots \wedge \operatorname{pj}_{q}^\ast(\widetilde{\omega}_{\sigma(q)}))\in \Omega^{3q}([0,1]\times \prod_{q}\bar{C}_{\mathbf{2}}(E^{+})).
\]
By Stokes' formula, one has 
\[
\int_{\partial ([0,1]\times \bar{C}_{\bfp}(E^{+}))}(\id\times \phi_{(\Gamma_{i},\rho_{i},o_{i})})^*(\widetilde{\omega}(q))=0.
\]
This can be rewritten as 
\begin{equation}\label{eq: boundary vanishing}
\int_{\bar{C}_{\bfp}(E^{+})}\phi_{(\Gamma_{i},\rho_{i},o_{i})}^{*}(\omega(q)-\omega'(q))=\sum_{A\subset \bfp, |A|\geq 2}\int_{[0,1]\times \partial_{A}\bar{C}_{\bfp}(E^{+})} (\id\times \phi_{(\Gamma_{i},\rho_{i},o_{i})})^*(\widetilde{\omega}(q)).    
\end{equation}
As proved in \cite[Lemma E.3, E.4, E.5]{Wa21}, the right hand side of (\ref{eq: boundary vanishing}) can be nonvanishing only when $A$ corresponds to an edge $e\in E'(\Gamma_{i})$. Since $\eta$ is a cocyle,   these terms cancel each other when we do linear combination with coefficient $a_{i}$. So we get
\[
\sum^{k}_{i=1}a_{i}\int_{\bar{C}_{\bfp}(E^{+})}\phi_{(\Gamma_{i},\rho_{i},o_{i})}^{*}(\omega(q))=
\sum^{k}_{i=1}a_{i}\int_{\bar{C}_{\bfp}(E^{+})}\phi_{(\Gamma_{i},\rho_{i},o_{i})}^{*}(\omega'(q)).
\]
\end{proof}
Next, we fix a specific propagator system $\cP$ as follows: First, we apply Lemma \ref{lem: extension} to find a map $P: \bar{C}_{\mathbf{2}}(E^{+})\to X(l)$ for some $l\gg 0$ that fits into the following commutative diagram 
\begin{equation*}
\begin{tikzcd}
\partial \bar{C}_{\mathbf{2}}(E^{+}) \arrow[d,hook] \arrow[r,"P^{\partial}"]  & X(1) \arrow[d,"g_{l}"]\\
\bar{C}_{\mathbf{2}}(E^{+}) \arrow[r,"P"] &  X(l)
\end{tikzcd}.
\end{equation*}
We let $d$ be the mapping degree of $g_{l}$. For any $\cS\in \cM$ and any $[i]\neq [j]\in \bfn/\cS-\{[0]\}$, we define the propagator $P^{\cS}_{[i],[j]}$ as the composition 
\[
\bar{C}_{\bfn/\cS}(E^{+})\xrightarrow{\phi_{[i],[j]}}\bar{C}_{\mathbf{2}}(E^{+})\xrightarrow{P}X(l)
\]
The map $\phi_{[i],[j]}$ is defined in (\ref{eq: configuration space forgetful map}). With such propagator system, the map 
$P(\Gamma_{l},\rho_{l})$ defined in 
(\ref{eq: graph propagator}) is exactly the composition 
\[
\bar{C}_{\bf n/\mathcal S}(E^{+})\xrightarrow{\phi_{(\Gamma_{l},\rho_{l},o_{l})}} \prod_{q} \bar{C}_{\mathbf{2}}(E^{+})\xrightarrow{\prod_{q}P} \prod_{q} S^3.
\]
By Sard's theorem, we can pick a generic point $(x_{1},\cdots, x_{q})\in \prod_{q} S^3$ such that for any $\sigma\in \Sigma_{q}$ and for any $l\in \{1,\cdots,k\}$, $(x_{\sigma(1)},\cdots, x_{\sigma(q)})$ is a regular value of the smooth map $P(\Gamma_{l},\rho_{l})$ and it doesn't belong to the image of $\partial \bar{C}_{\bfp}(E^{+})$.

We let $\alpha'_{i}\in \Omega^{3}(S^3)$ be a form supported near $x_{i}$ with $\int_{S^3}\alpha'_{i}=1$. And we let $\alpha_{i}=\frac{1}{d}g^{*}_{l}(\alpha'_{i})$. Then $\omega_{i}=\frac{1}{d}P^{*}(\alpha'_{i})$ is an extension of $P^{\partial,*}(\alpha_{i})$. By Lemma \ref{lem: invariance of differential form}, we have 
\[
I_{[\eta]}(\pi)=\sum_{1\leq l\leq k}a_{l}\cdot \int_{\bar{C}_{\bfp}(E^{+})} \omega'(q),
\]
where $\omega'(q)$ is defined in (\ref{eq: omega'q}). 

On the other hand, we let $\xi'\in C^{3q}_{\cech}(\prod_{q}S^3)$ be a Čech cocycle supported near $(x_{1},\cdots, x_{q})$ with $\int_{\prod_{q}S^3}\xi'=1$. And we let
\[\xi=\frac{1}{q!}\sum_{\sigma\in \Sigma_{q}}(\operatorname{sign}(\sigma)\cdot \sigma^{*}(\xi')).\]

\begin{Lemma}\label{lem: integrating form= integrating Čech cocycle}For any $1\leq l\leq k$, we have
\begin{equation}\label{eq: integrating form= integrating Čech cocycle}
\int_{\bar{C}_{\bfp}(E^{+})} \phi^*_{(\gamma_{l},\rho_{l},o_{l})}\omega'(q)=\frac{1}{d^q}\langle P(\Gamma_{l},\rho_{l})^{*}(\xi),[\bar{C}_{\bfp}(E^{+})]\rangle    
\end{equation}
\end{Lemma}
\begin{proof} For any $\sigma\in \Sigma_{q}$, we use $\vec{x}_{\sigma}$ to denote the point $(x_{\sigma(1)},\cdots, x_{\sigma(q)}))$.  We define the local degree of $P(\Gamma_{l},\rho_{l})$ at $\vec{x}_{\sigma}$ to be the signed count of points in $P(\Gamma_{l},\rho_{l})^{-1}(\vec{x}_{\sigma})$. We denote this by $\operatorname{deg}(P(\Gamma_{l},\rho_{l}),\vec{x}_{\sigma})$.
Since this local degree can be computed by integrating Čech cocycle or a differential form, one can check that both sides of (\ref{eq: integrating form= integrating Čech cocycle}) equal
\[
\frac{1}{d^q\cdot q!}\sum_{\sigma\in \sigma_{q}} \operatorname{deg}(P(\Gamma_{l},\rho_{l}),\vec{x}_{\sigma}).
\]
This finishes the proof. 
\end{proof}
By our definition of $I^{\Top}$ and Lemma \ref{lem: invariance of differential form}, we have proved the following main result of this section.

\begin{Proposition}\label{proposition: I^W=I^top} Let $[\widetilde{\eta}]\in H^{m,n}(\widetilde{\cG})$ be a homology class and let $\eta\in H^{m,n}(\cG)$ be its image under the forgetful map (\ref{eq: forgetful}). Let $D^4\to E_{D}\xrightarrow{\pi} B$ be a smooth disk bundle over a closed manifold $B$ of dimension $m+\frac{n}{2}$. Suppose $B$ is a sphere, or suppose $E_D$ is equipped with a  framing. Then one has $I_{[\eta]}(\pi)=I^{\Top}_{[\widetilde{\eta}]}(\pi)$. 
\end{Proposition}


\section{Applications to topology of diffeomorphism groups}\label{section: applications}
In this section, we prove Theorem \ref{thm: main} and its applications. We start with some preparations. 
\subsection{Obstruction theory} We briefly recall some basic facts about obstruction theory and Moore-Postnikov tower. See \cite{Hatcher} for more details. This will be used in the proof of Theorem \ref{thm: Sm noncontractible}.

Let $Y\to Z$ be a principal fibration whose fiber is the the Eilenberg-Maclane space $K(G;n)$ for some abelian $G$ and some $n\geq 1$. Let $A\hookrightarrow B$ be a cofibration. Given maps $h:B\to Z$ and $g: A\to Y$, such that the following diagram commutes 
\[
\xymatrix{A\ar@{^{(}->}[d]\ar[r]^{g}& Y\ar[d]\\
B\ar[r]^{h} & Z
},
\]
one can define a canonical obstruction class \[
o(g,h)\in H^{n+1}(B,A;G)
\]
such that $o(g,h)=0$ if and only if the relative lifting problem
\[
\xymatrix{A\ar@{^{(}->}[d]\ar[r]^{g}& Y\ar[d]\\
B\ar[r]^{h}\ar@{-->}[ur]^{\widetilde{h}} & Z
}
\]
has a solution. In this case, any two choices of $\widetilde{h}$ are related by to each other by twisting with element in $H^{n}(B,A;G)$. In particular,  then  $\widetilde{h}$ is unique up to homotopy if $H^{n}(B,A;G)=0$.

The obstruction class is functorial. More precisely, given the following commutative diagram 
\[
\xymatrix{A'\ar@{^{(}->}[d]\ar[r]^{g'}&A\ar@{^{(}->}[d]\ar[r]^{g}& Y\ar[d]\ar[r]^{g''}&Y'\ar[d]\\B'\ar[r]^{h'}&
B\ar[r]^{h} & Z\ar[r]^{h''}& Z'
},
\]
such that $A'\hookrightarrow B'$ is another cofibration, and that the rightmost part is a pullback square. We have 
\begin{equation}\label{eq: obstruction functorial}
o(g''\circ g\circ g, h''\circ h\circ h')=(h{'})^*(o(g,h))\in H^{n}(B',A';G).    
\end{equation}

Now we consider a simple fibration $q:C\to D$ whose fiber $F$ is a simple space. This fibration has a 
principal Moore-Postnikov tower
\[
C\simeq Z_{\infty}\to \cdots \to Z_{n+1}\xrightarrow{ p_{n+1}} Z_{n}\to \cdots \to Z_{0}=D.
\]
Here $p_{n+1}: Z_{n+1}\to Z_{n}$ is the homotopy fiber of a map $Z_{n}\to K(\pi_{n+1}(F),n+2)$, represented by an element \[\bfk_{n}(C,D)\in H^{n+2}(Z_{n}; \pi_{n+1}(F))\] called the $n$-th $\bfk$-invariant. The Moore-Postnikov tower and the $\bfk$-invariants are both functorial under pullbacks.

Consider the relative lifting problem 
\[
\xymatrix{A\ar@{^{(}->}[d]\ar[r]^{f}& C\ar[d]\\
B\ar[r]^{l}\ar@{-->}[ur] & D.
}
\]
We let $l_0=l$ and inductively lift $l_n: B\to Z_n$ to $l_{n+1}: B\to Z_{n+1}$. The corresponding obstruction is denoted as 
\[
o_{n}(f,l_n)\in H^{n+2}(B,A;\pi_{n+1}(F)).
\]
In general, different choices of $l_n$ gives different class $o_{n}(f,l_n)$. We define 
\[
o_{n}(f,l):=\{o_{n}(f,l_n)\mid l_n \text{ lifts }l\}\subset H^{n+2}(B,A;\pi_{n+1}(F)).
\]
(We set $o_{n}(f,l)=\emptyset$ when a lift $l_n$ doesn't exist.) 
A lift $l_{n+1}$ exists if and only if $0\in o_{n}(f,l)$.

We say $o_{n}(f,l)$ has \emph{no indeterminacy} if it contains at most one element. The following lemma summarize some special cases where there are no indeterminacy.

\begin{Lemma}\label{lem: obstruction independence with lift} $o_{n}(f,l)$ has no indeterminacy one of the following condition holds
\begin{enumerate}
\item $\bfk_{n}(C,D)=0$, or
\item $\bfk_n(C,D)$ is pulled back from a homology class 
$\bfk'\in H^{n+2}(Z_{n'};\pi_{n+1}(F))$ for some $0\leq n'\leq n$ that satisfies 
\begin{equation}\label{eq: relative homology vanishing}
H^{k}(B,A,\pi_{k}(F))=0,\quad \forall 1\leq k\leq n'.    
\end{equation}
\end{enumerate}
\end{Lemma}
\begin{proof} (1) Since $\bfk_{n}(C,D)=0$, the fibration $Z_{n+1}\to Z_{n}$ is trivial. So we have a projection to fiber map $\pj: Z_{n+1}\to K(\pi_{n+1}F,n+1)$ that fits into the  pullback square 
\[
\xymatrix{Z_{n+1}\ar[d]^{p_{n+1}}\ar[r]^-{\pj}& K(\pi_{n+1}F,n+1)\ar[d]\\
Z_n\ar[r]& \ast
}.
\]
By functoriality of the obstruction class (\ref{eq: obstruction functorial}), the class $o_{n}(g,l_n)$ equals the obstruction class of the lifting problem 
\[
\xymatrix{A\ar[r]^-{f}\ar@{^{(}->}[d]& C\ar[r] & Z_{n+1} \ar[r]^-{\pj}& K(\pi_{n+1}F,n+1)\ar[d]\\
B\ar@{-->}[urrr]\ar[rrr]& & &\ast
},
\]
so it only depends on $f$.

(2) $\bfk'$ gives a map $Z_{n'}\to K(\pi_{n+1}F,n+2)$. We denote the homotopy fiber by $W$. Then we have the pullback square of principal fibrations 
\[
\xymatrix{Z_{n+1}\ar[rr]\ar[d]& & W\ar[d]\\
Z_{n}\ar[rr]^-{p_{n'+1}\circ\cdots\circ p_{n}}& &Z_{n'}
}.
\]
By functorality of obstruction classes, $o_{n}(g,l_{n})$ equals the obstruction class of the lifting problem
\[
\xymatrix{A\ar[r]^-{f}\ar@{^{(}->}[d]& C\ar[r] & Z_{n+1} \ar[r]& W\ar[d]\\
B\ar@{-->}[urrr]\ar[rrr]^{p_{n'+1}\circ\cdots\circ p_{n}\circ l_{n}}& & &Z_{n'}
},
\]
So $o_{n}(g,l_{n})$ depends on $l_n$ only through the composition \[p_{n'+1}\circ\cdots\circ p_{n}\circ l_{n}: B\to Z_{n'}.\] This is a lift of $l$. However, by our assumption (\ref{eq: relative homology vanishing}), such a lift is unique up to homotopy.
\end{proof}
\subsection{Construction of bundles and linear structures}  
As we mentioned in the introduction, given any admissible, oriented trivalent graph $(\Gamma,o)$ with $2k$ vertices, Watanabe constructed a smooth bundle 
\[
D^4\to E^{\Gamma}_{D}\to S^{k}
\]
which is framed and boundary trivialized. In this subsection, we construct various bundles from $E^{\Gamma}_{D}$. 

First, any positive integer $n$, we let \[D^4\to E^{n\Gamma}_D\to S^k\] be the pull back of $D^4\to E^{\Gamma}_D$ under a degree-$n$ map $S^k\to S^k$. Consider an integer coefficient linear combination 
\[
\eta=\sum^{k}_{i=1} n_{i}(\Gamma_{i},o_{i})
\]
of admissible, oriented trivalent graphs with $2k$ vertices. Then $\eta$ represents an integer class in $\cA^{\even}_{k}$. By gluing together the bundles $E^{n_{i}\Gamma_{i}}_{D}\to S^{k}$ along the fiber over the base point, we obtain a bundle over $\vee_{n}S^{k}$. Then we 
pull back this bundle via the quotient map $S^{k}\to \vee_{n}S^{k}$ and obtain the bundle 
\[
D^4\to E^{\eta}_{D}\to S^{k}.
\]
We use $\alpha_{\eta}\in \pi_{k}(\BDiff_{\partial}(D^4))$ to denote the correspond element. 
Given a smooth 4-manifold $X$ (possibly with boundary) and a smooth embedding $D^4\hookrightarrow \Inte(X)$, we let 
\[
X\to E^{\eta}_{X}\to S^k
\]
be the bundle obtained from $E^{\eta}_{D}\to S^{k}$ by attaching the product bundle $S^k\times (X-\Inte(D^4))$. Note that the bundle $D^4\to E^{\eta}_{D}\to S^k$ has a canonical trivialization as a topological bundle because $\Homeo_{\partial}(D^4)$ is contractible. So it induces a smooth structure on the restriction of the trivial bundle 
\[
D^4\to D^{k+1}\times D^4\to D^{k+1}
\]
to $\partial D^{k+1}$. This smooth structure, together with the standard smooth structure on $\partial D^4\times D^{k+1}\to D^{k+1}$, induces a linear structure  on the microbundle bundle 
\[
\cT^{v}(D^{k+1}\times D^4)|_{\partial(D^{k+1}\times D^4)},
\]
We denote this linear structure by $l(\eta,\partial D^{k+5})$. Given any smooth $X$ with an embedding $D^4\hookrightarrow \Inte X$, consider the product bundle $X\to D^{k+1}\times X\to D^{k+1}$. Then we can use the linear structure on $\cT X$ to extend $l(\eta,\partial D^{k+5})$ to a linear structure on $\cT^{v}(D^{k+1}\times X)|_{D^{k+1}\times X-\Inte D^{k+5}}$, which we denote as \[l(\eta, D^{k+1}\times X-\Inte D^{k+5}).\] We use $l(\eta, \partial(D^{k+1}\times X))$ to denote the restriction \[l(\eta, D^{k+1}\times X-\Inte D^{k+5})|_{\partial(D^{k+1}\times X)}.\]

\subsection{Proof of the main theorem and its applications} Now we are ready to prove the main result and its applications.
\begin{proof}[Proof of Theorem \ref{thm: main}] Consider the trivial bundle 
\[
D^4\to E_{D}=(S^m\times D^4)\to S^m.
\]
The vertical tangent microbundle $\cT^{v}E_{D}$ has a canonical trivilization. So concordance classes of framings on $E_{D}$ that are standard on $(S^m\times \partial D^4)\cup (\{\infty\}\times D^4)$ one-one correspond to homotopy classes of maps 
\[
(D^{m+4},S^{m+3})\to (\Top(4),*).
\]
By Lemma \ref{lemma: concordance invaraince}, we have a well-defined map 
\begin{equation}\label{eq：I-top-eta}
I^{\Top}_{i_*[\eta]}:\pi_{m+\frac{n}{2}+4}(\Top(4)/O(4))\otimes \bQ\to \bQ \end{equation}
for any $[\eta]\in H^{m,n}(\mathcal{G})$. Here $i_*$ is given in (\ref{eq: average}). By Lemma \ref{lemma: concordance invaraince}, the map (\ref{eq：I-top-eta}) is a group homomorphism and is linear with respect to $[\eta]$. So we can rewrite it as a map 
\[
I^{\Top}_{m,n}: \pi_{m+\frac{n}{2}+4}(\Top(4))\otimes \bQ\to H_{m,n}(\mathcal{G}).
\]
Specializing to the case $m=0$ and $k=\frac{n}{2}$, we get the map 
\[
I^{\Top}: \pi_{k+4}(\Top(4))\otimes \bQ\to \mathcal{A}^{\even}_{k}.
\]
Commutativity of the diagram (\ref{eq: diagram I-top=I-W}) follows from the equivalence of $I$ and $I^{\Top}$ for smooth bundles (see Proposition \ref{proposition: I^W=I^top}).
\end{proof}

\begin{proof}[Proof of Corollary \ref{cor: Top and Homeo}] Watanabe proved that the map $I: \pi_{k}(\BDiff^{\Fra}_{\partial}(D^4))\otimes \bQ\to \cA^{\even}_{k}$ is surjective by showing that 
$I(\alpha_{\Gamma})=[(\Gamma,o)]$.  Here $\alpha_{\Gamma}$ is the element that classifies the bundle $E^{\Gamma}_D\to S^k$. By commutativity of the diagram (\ref{eq: diagram I-top=I-W}), the map $I^{\Top}$ is surjective as well. This directly implies Part (1).

The action of $\Homeo(S^4)$ on $S^4$ gives us a fibration 
$$
\Top(4) \to \Homeo(S^4) \to  S^4
$$
where the second map is defined by the action of $\Homeo(S^4)$ on  $\infty\in S^4$.
Since $\pi_i(S^4)\otimes \bQ=0$ when $k\ge 8$, the long exact sequence of homotopy groups implies 
$$
\dim \pi_k(\Homeo(S^4))\otimes \bQ =\dim \pi_k(\Top(4))\otimes \bQ \ge \dim(\cA^{\even}_{k-4})
$$
when $k\ge 8$. The standard action of $O(5)$ on $S^4$ gives us a fibration $O(4)\to O(5)\to S^4$ which
fits into the following commutative diagram
\begin{equation}
\begin{tikzcd}
  O(4)\arrow[r]\arrow[d] & O(5) \arrow[r] \arrow[d] & S^4\arrow[d,"="] \\
 \Top(4)\arrow[r] & \Homeo(S^4) \arrow[r]   & S^4  
\end{tikzcd}
\end{equation} 
By taking the homotopy groups, we obtain the commutative diagram
\begin{equation*}
\xymatrix{
 \pi_7^\bQ(S^4) \ar[r] \ar[d]& \pi_6^\bQ(O(4)) \ar@{->}[r]\ar[d] & \pi_6^\bQ(O(5)) \ar[d] \ar[r]& 0 \ar[d]\\
 \pi_7^\bQ(S^4) \ar[r] & \pi_6^\bQ(\Top(4)) \ar@{->}[r] & \pi_6^\bQ(\Homeo(S^4)) \ar[r] & 0
}
\end{equation*}
where the two rows are exact and $\pi_i^\bQ$ denotes the rational homotopy group.
 Since $\pi_6^\bQ(O(4))=0$,
the diagram implies the map $\pi_7^\bQ(S^4)\to \pi_6^\bQ(\Top(4))$ is zero. Hence 
$$
\dim \pi_6^\bQ(\Homeo(S^4))=\dim \pi_6^\bQ(\Top(4)) \ge \dim(\cA^{\even}_{2}\otimes \bQ)=1.
$$
This completes the proof of Part (2) 
since
$$
\dim\cA^{\even}_{0}=\dim\cA^{\even}_{1}=\dim\cA^{\even}_{3}=0.
$$

Now we focus on the identity component $\Top^+(4)$ of $\Top(4)$. 
Since $\Top^+(4)$ is homotopy equivalent to the loop space $\Omega \BTop^+(4)$ and $\BTop^+(4)$ is simply connected, 
 there is a rational homotopy equivalence
\begin{equation}\label{eq_Top4_K_decomposition}
\Top^+(4) \simeq_\bQ \prod_{k\ge 2 } K(\pi_k^\bQ(\Top^+(4)); k) 
\end{equation}
according to \cite{FHT-rational_homotopy}*{Corollary of Proposition 16.7}.
The linear surjection 
\[
I^{\Top}: \pi_k^\bQ(\Top^+(4))\cong \pi_{k}(\Top(4))\otimes \bQ\to\cA^{\even}_{k-4}\] 
induces a map $F_k:K(\pi_k^\bQ(\Top^+(4)); k) \to K(\cA^{\even}_{k-4}; k)$. One can also find 
a right inverse of $I^{\Top}$ which induces a right inverse of $F_k$ (in the homotopy sense). Therefore $F_k$ induces 
a surjection on homology and an injection on cohomology. By combining these $F_k$ we obtain 
a map
$$
F: \prod_{k\ge 2 } K(\pi_k^\bQ(\Top^+(4)); k) \to \prod_{k\ge 2 } K(\cA^{\even}_{k-4}; k) 
$$
which induces an injection on cohomology. Since the cohomology group of the right-hand-side is just the symmetric algebra generated by $V$, this completes the proof of the $\Top(4)$ case in Part (3). The proof for the $\Homeo(S^4)$ case is similar.
Part (4) follows directly from Part (3).
\end{proof}

Let $\Gamma_4$ be the complete graph with $4$ vertices and let $o$ be any orientation on $\Gamma$. Recall that $\cA^{\even}_{2}$ is generated by $(\Gamma_4,o)$. We let $\theta_0\in V_6=(\cA^{\even}_{2})^\ast$ be the dual basis. According to 
Corollary \ref{cor: Top and Homeo}, $\theta_0$ also determines an element in $H^6(\Top^+(4);\bQ)$, which we still denote by $\theta_0$. 
We can also view $\theta_0$ as an element in $E^2_{0,6}\cong H^6(\Top^+(4);\bQ)$ of 
the (rational) cohomology spectral sequence 
for the fibration $\Top^+(4)\to \ETop^+(4)\to \BTop^+(4)$.
\begin{Proposition}
The element $\theta_0\in E^2_{0,6}$ is transgressive and $\theta_1:=\tau(\theta_0)$ is a non-zero element in $E_{7,0}^7\cong H^7(\BTop^+(4);\bQ)$ where $\tau$ denotes the transgression.
\end{Proposition}
\begin{proof}
\begin{figure}
\begin{tikzpicture}
  \matrix (m) [matrix of math nodes,
    nodes in empty cells,nodes={minimum width=5ex,
    minimum height=6ex,outer sep=-5pt},
    column sep=1ex,row sep=1ex]{
          6     &   \theta_0   &     &     & \\
          3     &  \bQ^2 &  \bQ^4 &  & \\
          0     &  \bQ  & \bQ^2 &    & \\
    \quad\strut &   0  &  4  &  7 & \strut \\};
  \draw[-stealth] (m-1-2.south east)-- (m-2-3.north west) ;
\draw[thick] (m-1-1.east) -- (m-4-1.east);
\draw[thick] (m-4-1.north) -- (m-4-5.north) ;
\end{tikzpicture}
\caption{}\label{figure_spectral_sequence_BTop4}
\end{figure}

According to \eqref{eq_Top4_K_decomposition} and the fact that $\Top^+(4)/SO(4)$ is rationally 4-connected, we have
$$
H^i(\Top^+(4);\bQ)=0~\text{for}~i\le 2~\text{or}~i=4,~ H^3(\Top^+(4);\bQ)= H^3(SO(4);\bQ)=\bQ^2
$$
and
$$
H^6(\Top^+(4);\bQ)= H^6(SO(4);\bQ)\oplus  
(\pi_6^\mathbb{Q}(\Top^+(4))^\ast.
$$
Notice that $\theta_0\in(\pi_6^\mathbb{Q}(\Top^+(4))^\ast $.
By Hurewicz theorem, we have 
$$
H^i(\BTop^+(4);\bQ)=0~\text{for}~i\le 3,~ H^4(\BTop^+(4);\bQ)=\bQ^2.
$$
It is straightforward to check that there is no non-zero entry mapped to $E^i_{7,0}$ in the spectral sequence
 (see Figure \ref{figure_spectral_sequence_BTop4})
until $i=7$. Therefore we have
$E^7_{7,0}=E^2_{7,0}\cong H^7(\BTop^+(4);\bQ)$.
The first  differential on $E^i_{0,6}$ with non-zero target is $$d_4: E^4_{0,6}=H^6(\Top^+(4);\bQ)\to E^4_{4,3}.$$
By  construction we see that the pullback of $\theta_0$ to $H^6(SO(4);\bQ)$ is zero. By comparing with the spectral sequence
for $SO(4)\to ESO(4)\to BSO(4)$ we deduce that $d_3(\theta_0)=0$.
The next differential on $E^i_{0,6}$ with non-zero target
is the transgression $\tau$. Therefore we conclude that $\theta_0$ is transgressive. 
Since $\ETop^+(4)$ is contractible, $\tau(\theta_0)$ must be non-zero. 
\end{proof}

Now we assume that $\BO\to \BTop(4)$ is a cofibration. (This can be obtained by replacing $\BTop(4)$ with the mapping cylinder.)  Since 
\[H^{6}(\BO;\bQ)=H^{7}(\BO;\bQ)=0,\] the inclusion $(\BTop(4),\emptyset)\to (\BTop(4),\BO)$ induces an isomorphism \[H^{7}(\BTop(4);\bQ)\cong H^{7}(\BTop(4),\BO;\bQ).\] We let $\theta_2$ be the image of $\theta_1$ under this isomorphism. 

Let $A\hookrightarrow B$ be a cofibration and let\[
\mathbb{R}^4\to E\xrightarrow{\pi} B.
\]
be a topological bundle. We assume $E|_{A}$ is has a linear structure $l$. Then the pair $(\pi, l)$ is classified by (the homotopy class) of a map 
\[
f_l:(B,A)\to (\BTop(4),\BO).
\]
We use $\theta(l)\in H^{7}(B,A;\bQ)$ to denote $f_{l}^{*}(\theta_2)$. Note that $l$ can be extended to a linear structure on $E$ only if $\theta(l)=0$. 

Now we consider the  special case $(B,A)=(D^{7}, \partial D^7)$.  Then $f_l$ represents a class in the group $\pi_{7} (\BTop(4),\BO)\otimes \mathbb{Q}$. We let $\alpha_l\in \pi_{7}(\BTop(4))\otimes \bQ$ be the image of $[f_l]$ under the isomorphism 
\[
\pi_{7}(\BTop(4),\BO)\otimes \bQ\cong \pi_{7}(\BTop(4))\otimes \bQ.
\]
\begin{Lemma}\label{lem: computation of theta_l on D7}    
$\langle \theta(l(n\Gamma_4,\partial D^7)),[D^7]\rangle =n$.
\end{Lemma}
\begin{proof}
Given any $\alpha\in \pi_{7}(\BTop(4))\otimes \bQ$, we let $\alpha'\in H_{7}(\BTop(4);\bQ)$ be its Hurewicz image. Since $\theta_0$ is the dual basis of $[(\Gamma_4,o)]$, we see that  
\begin{equation}\label{eq: theta_1 paired with alpha}
I(\alpha)=
\langle \theta_1,\alpha'\rangle \cdot [(\Gamma_4,o)]\in \cA^{\even}_{2}.    
\end{equation}
Setting $\alpha=\alpha_{l}$ in (\ref{eq: theta_1 paired with alpha}), we get
\[I^{\Top}(\alpha_{l})=\langle \theta(l),[D^7]\rangle \cdot [(\Gamma_4,o)].\]
Now we let $l=l(n\Gamma_4,\partial D^7)$. By Theorem \ref{thm: main} and Watanabe's calculation that $I(\alpha_{nl})=n[(\Gamma_4,o)]$, we see that $\langle \theta(l),[D^7]\rangle=n$.
\end{proof}
\begin{Lemma}\label{lem: theta(l) calculation} Let $A$ be compact topological manifold of dimension 7 and let $D^{7}\hookrightarrow A$ be an orientation preserving embedding. Let $B=A\setminus \Inte(D^7)$. Consider a topological bundle $\mathbb{R}^4\to E\to A$ and a linear structure  $l(B)$ on $E|_{B}$. Let $l(\partial A))$ and $l(\partial D^7)$ be the restriction of $l(B)$ to $\partial A$ and $\partial D^7$ respectively. Suppose $l(\partial D^7)=l(n\Gamma_4,\partial D^7)$ under some trivialization of $E|_{D^7}$. Then we have \[\langle\theta(l(B)),[A]\rangle=\langle\theta(l(\partial A)),[A]\rangle=n.\] In particular $l(\partial A)$ is not extendable over $A$ if $n\neq 0$.
\end{Lemma}
\begin{proof} This follows from Lemma \ref{lem: computation of theta_l on D7} and the naturality of $\theta(-)$ under the pull back along the embedding $(D^7,\partial D^7)\hookrightarrow (A,B)$. 
\end{proof}
Applying Lemma \ref{lem: theta(l) calculation} to the linear structure $l(n\Gamma_4,\partial (D^3\times M))$, we immediately get the following corollary.
\begin{Corollary}\label{cor: l-gamma4 nonextendable}
For any $n\neq 0$, the linear structure $l(n\Gamma_4,\partial (D^3\times M))$ on $\cT^{v}(D^3\times M)|_{\partial(D^3\times M) }$ does not extend to $\cT^{v}(D^3\times M)$.    
\end{Corollary}

\begin{Lemma}\label{lem: beta-vanishing}
Let $\alpha\in \pi_{k}(\cM^{s}(D^4\# n(S^2\times S^2)))$ be the element given by the bundle $E^{\eta}_{D^4\# n(S^2\times S^2))}\to S^{k}$. Let $\beta\in H_{k}(\cM^{s}(D^4\# n(S^2\times S^2));\bQ)$ be the Hurewicz image of $\alpha$. Then $\beta=0$ when $n$ is large enough.
\end{Lemma}
\begin{proof}
In \cite[Theorem 1.5]{Galatius2017}, Galatius and Randal-Williams proved a generalised Madsen–Weiss theorem. As a corollary, the rational cohomology ring of the stable moduli space $\cM^{s}(D^{4}\# \infty(S^2\times S^2))$
is generated by the Morita-Miller-Mumford $\kappa$-classes. Since $E^{\eta}_{D^4\# n(S^2\times S^2))}$ is obtained by gluing a trivial bundle to the framed bundle $E^{\eta}_{D}\to S^k$, these $\kappa$-classes evaluates zero on $\beta$. Hence we see that $\beta=0$ in 
$$\lim_{\substack{\rightarrow\\n}} H_{k}(\cM^{s}(D^{4}\# n(S^2\times S^2);\bQ)).$$
\end{proof}

Now we are ready to prove Theorem \ref{thm: MMM class}.
\begin{proof}[Proof of Theorem \ref{thm: MMM class}] We let $\theta=\theta_1$. Since $H^{7}(\BO;\mathbb{Q})=0$, the pull back of $\theta$ to $H^{7}(\BO;\bQ)$ is zero. This implies that $\kappa_{\theta}(\pi)=0$ if the structure group of $\pi$ can be reduced to $O(4)$. It remains to show that $\kappa_{\theta}(\pi)$ is nontrivial for some bundle $\pi$. By Lemma \ref{lem: beta-vanishing},
 for suitable $m,n>0$ , the map $S^{2}\to \cM^{s}(D^{4}\# n(S^2\times S^2))$ which  classifies the bundle 
\[
D^{4}\# n(S^2\times S^2)\to E^{m\Gamma_{4}}_{D^{4}\#n(S^2\times S^2)} \to S^2
\]
is null homologous. This means the bundle $E^{m\Gamma_{4}}_{D^{4}\#n(S^2\times S^2)} \to S^2$ bounds a smooth, boundary trivialized bundle 
\[
D^{4}\# n(S^2\times S^2)\to E'\to \mathring{Y}
\]
over a punctured 3-manifold $\mathring{Y}=Y-\Inte(D^4)$. Let $X$ be a closed, oriented, smooth 4-manifold and let $\mathring{X}=X-\Inte(D^4)$. By gluing $E'\to \mathring{Y}$ with the product bundle $\mathring{X}\times \mathring{Y}\to \mathring{Y}$, we obtain a smooth bundle 
\[
X\# n(S^2\times S^2)\to E'_{X}\to \mathring{Y}
\]
that is bounded by the bundle 
\[
X\# n(S^2\times S^2)\to E^{m\Gamma_4}_{X\# n(S^2\times S^2)}\to S^2.
\]
The trivialization of the topological bundle $E^{m\Gamma_4}_{D}\to S^2$ extends to a trivialization of $E^{m\Gamma_4}_{X\# n(S^2\times S^2)}\to S^2$. This allows us to glue $E'_{X}$ with $D^{3}\times (X\# n(S^2\times S^2))$ and form the topological bundle  
\[
X\# n(S^2\times S^2)\to E\xrightarrow{\pi} Y.
\]
Note that we have an embedding 
\[
D^{7}=D^{3}\times D^4 \hookrightarrow D^{3}\times (X\# n(S^2\times S^2)\hookrightarrow E
\]
and the bundle $E\to Y$ has a  smooth structure on the complement of $\Inte(D^7)$. Let $l$ be the induced linear structure on $\cT^v E|_{E-\Inte(D^7)}$. Then $l|_{\partial D^7}$ equals $l(m\Gamma_4,\partial D^7)$. By Lemma \ref{lem: theta(l) calculation}, we have $\langle \theta(l), [E]\rangle =m$.
This implies \[\langle \kappa_{\theta}(\pi),[Y]\rangle=\langle \theta(\pi), [E]\rangle =\langle \theta(l), [E]\rangle=m\neq 0. \] This finishes the proof.
\end{proof}

Now we proceed to the proof of Theorem \ref{thm: Sm noncontractible}.  We have the following more precise result.
\begin{Theorem}\label{thm: bundle nonextension}  Let $X$ be a compact, orientable smooth 4-manifold and let $\eta$ be a nonzero integer class in $\cA^{\even}_{k}$. Suppose one of the following conditions hold
\begin{enumerate}
    \item $k=2$;
    \item $\partial X\neq \emptyset$;
    \item $\sigma(X)=0$.
\end{enumerate}
Then the linear structure $l(\eta,\partial(D^{k+1}\times X))$ on $\cT^{v}(D^{k+1}\times X)|_{\partial(D^{k+1}\times X)}$ can not be extended to $\cT^{v}(D^{k+1}\times X)$. 
\end{Theorem}

\begin{proof}[Proof of Theorem \ref{thm: Sm noncontractible}] The linear structure $l(\eta,\partial(D^{k+1}\times X))$  is induced by a smooth structure on the bundle 
$\partial D^{k+1}\times X \to \partial D^{k+1}$. Such a smooth structure is classified by a map  $s_{\eta}: \partial D^{k+1} \to \BDiff_{\partial}(X)$.   
Since the linear structure can not be extended, the smooth structure can not be extended as well. So the following extension problem 
\begin{equation*}
 \xymatrix{\partial D^{k+1}\ar[r]^{s_{\eta}}\ar@{^{(}->}[d]& \BDiff_{\partial}(X)\ar[d]^{i^{X}}\\
D^{k+1}\ar@{-->}^{\not\exists}[ur]\ar[r] & \BHomeo_{\partial}(X)}      
\end{equation*}
has no solution. The above commutative diagram gives a canonical nullhomotopy $h$ of $i^{X}\circ s_{\eta}$. By definition of $\Sm(X)$, the pair $(s_{\eta},h)$ gives a nontrivial element in $\pi_{k}(\Sm(X))\otimes \bQ$.
\end{proof}

Now we prove Theorem \ref{thm: bundle nonextension}. The case $k=2$ is simply a restatement of Corollary \ref{cor: l-gamma4 nonextendable} so has been proved. Proof of the other two cases will occupy the rest of this section.  

We consider the universal bundle 
\[SO(4)\to \ESO\to \BSO.\] Let $F=\Top(4)/O(4)$. Then $SO(4)$ acts on $F$ by the left multiplication. We let $F\to P^{u}\to \BSO$ be the associated bundle. Since $F, P^{u}$ and $\BSO$ are all simply-connected, we may take a localization at $\bQ$ and obtain a fibration 
\begin{equation}\label{eq: universal FQ bundle}
F_{\bQ} \to P^{u}_{\bQ}\to \BSO_{\bQ}.    
\end{equation}

Consider the frame bundle $SO(4)\to \Fr^{X}\to X$, classified by a map $f^{X}: X\to \BSO$. Consider 
the associated bundle 
\[
F\to P^{X}\xrightarrow{\pi^{X}} X.
\]
Then a linear structure on $X$ is given by a section of this bundle. By composing with the localization map $\BSO\to \BSO_{\bQ}$, we obtain the map $f^{X}_{\bQ}: X\to \BSO_{\bQ}$. By pulling back (\ref{eq: universal FQ bundle}) via $f^{X}_{\bQ}$, we obtain the fibration 
\[F^{X}_{\bQ}\to P^{X}_{\bQ}\xrightarrow{\pi^{X}_{\bQ}} X.\] 
which fits into the commutative diagram 
\[
\xymatrix{ F\ar@{^{(}->}[r]\ar[d]& P^{X}\ar[r]\ar[d]& X\ar@{=}[d]\\
F_{\bQ}\ar@{^{(}->}[r]& P^{X}_{\bQ}\ar[r]& X
}.
\]
\begin{Lemma}\label{lem: bundle factor though S4} (1)  Suppose $\partial X\neq \emptyset$. Then $f^{X}_{\bQ}$ is null homotopic. So $P^{\bQ}_{X}$ is trivial.

(2) Suppose $X$ is closed, then $f^{X}_{\bQ}$ is homotopic to a composition
\[
X\xrightarrow{q} S^4\xrightarrow{f} \BSO_{\bQ}.
\]
Here 
\begin{equation}\label{eq: pinch}
q: X\to S^4    
\end{equation}
is the map that pinches all points in $X\setminus \Inte(D^4)$ into a single point $\infty\in S^{4}$, and $f$ is a continuous map. If we further assume  $\sigma(X)=0$, then $f$ has a decomposition 
\[
S^4\xrightarrow{q'} S^4\xrightarrow{f^{S^{4}}} \BSO\to \BSO_{\bQ}
\]
Here $q'$ is a map of degree $\frac{\chi(X)}{2}$, and $f^{S^4}$ is the map that classifies $TS^4$.
\end{Lemma}
\begin{proof} Since $\BSO_{\bQ}$ is 3-connected, $f^{X}_{\bQ}$ is null homotopic if $\partial X\neq\emptyset$ and it factors through the map $q: X\to S^4$ when $X$ is closed. Furthermore, we have
\[
\pi_{4}(\BSO_{\bQ})\cong \pi_{4}(\BSO)\otimes \bQ \cong \pi_{3}(SO(4))\otimes \bQ\cong \bQ\oplus \bQ.
\]
The two summands corresponds to the Euler number and signature respectively. This implies the conclusion when $\sigma(X)=0$. 
\end{proof}

Now we let $\eta$ be a nonzero integral class in $\cA^{\even}_{k}$. Recall that we have defined the linear structure $l(\eta,D^{k+1}\times X-\Inte(D^{k+5}))$ on $\cT^{v}(D^{k+1}\times X)|_{D^{k+1}\times X-\Inte(D^{k+5})}$ and its restriction $l(\eta,\partial(D^{k+1}\times X))$ on $\cT^{v}(D^{k+1}\times X)|_{\partial(D^{k+1}\times X)}$. They give rise to sections 
\[
\gamma^{\circ}: D^{k+1}\times X-\Inte(D^{k+5})\to D^{k+1}\times P^{X}
\]
and
 \[
\gamma^{\partial}: \partial{D^{k+1}\times X}\to D^{k+1}\times P^{X}
\]
respectively. We let $\gamma^{\circ}_{\bQ}$ and  $\gamma^{\partial}_{\bQ}$ be their compositions with the rationalization map $D^{k+1}\times P^{X}\to D^{k+1}\times P^{X}_{\bQ}$. Note that linear structures on $\cT^{v}(D^{k+1}\times X)$ one-to-one correspond to sections of the bundle $D^{k+1}\times P^{X}\to D^{k+1}\times X$. So Theorem  \ref{thm: bundle nonextension} asserts that the lifting problem
\[\xymatrix{\partial(D^{k+1}\times X)\ar[r]^-{\gamma^{\partial}}\ar@{^{(}->}[d]& D^{k+1}\times P^{X}\ar[d]^{\id\times\pi^{X}}\\
D^{k+1}\times X\ar@{-->}[ur]\ar@{=}[r] & D^{k+1}\times X}
\]
has no solution. We consider the following lifting problem instead
 \begin{equation}\label{eq: rational lifting}
 \xymatrix{\partial(D^{k+1}\times X)\ar[r]^-{\gamma^{\partial}_{\bQ}}\ar@{^{(}->}[d]& D^{k+1}\times P^{X}_{\bQ}\ar[d]^{\id\times \pi^{X}_{\bQ}}\\
D^{k+1}\times X\ar@{-->}[ur]\ar@{=}[r] & D^{k+1}\times X}    
 \end{equation}
Theorem \ref{thm: bundle nonextension} directly follows from the following propositions.

\begin{Proposition}\label{prop: nonextension nonempty boundary}
Suppose $\partial X\neq \emptyset$ and suppose $\eta\neq 0\in \cA^{\even}_{k}$. Then (\ref{eq: rational lifting}) has no solution.
\end{Proposition}

\begin{Proposition}\label{prop: signature zero}
Suppose $\sigma(X)=0$ and suppose $\eta\neq 0\in \cA^{\even}_{k}$. Then (\ref{eq: rational lifting}) has no solution.
\end{Proposition}

We need some preparations before proving these results. 

By Lemma \ref{lem: bundle factor though S4}, $P^{X}_{\bQ}\to X$ is pulled back from  a fibration $P^{S}_{\bQ}\to S^4$ so is a simple fibration. We consider its Moore-Postnivkov tower 
\begin{equation}\label{eq: tower X}
P^{X}_{\bQ}\simeq Z^X_{\infty}\to \cdots \to Z^X_{n}\to \cdots\to Z^S_{0}=X,     
\end{equation}
which is pulled back from the tower 
\begin{equation}\label{eq: tower S}
P^{S}_{\bQ}\simeq Z^S_{\infty}\to \cdots \to Z^S_{n}\to \cdots \to Z^S_{0}=S^4. 
\end{equation}
We also consider the Postnivkov tower of the fiber $F_{\bQ}$
\begin{equation}\label{eq: tower F}
F_{\bQ}\simeq Z^F_{\infty}\to \cdots\to Z^F_{n}\to \cdots \to Z^F_{0}=*.
\end{equation}

Since 
\[
H^{m}(D^{k+1}\times X, (D^{k+1}\times X)-\Inte(D^{k+5}); G)=0
\]
for any $0\leq m\leq k+5$ and any abelian group $G$. The section $\gamma^{\circ}_{\bQ}$ has a unique lifting $l_{k+3}: D^{k+1}\times X\to D^{k+1}\times Z_{k+3}^X$ that fits into the commutative diagram 
 \[
 \xymatrix{(D^{k+1}\times X)-\Inte(D^{k+5})\ar[r]^-{\gamma^{\circ}_{\bQ}}\ar@{^{(}->}[dd]& D^{k+1}\times P^{X}_{\bQ}\ar[d]\\ & D^{k+1}\times Z^{X}_{k+3}\ar[d]\\
D^{k+1}\times X\ar@{=}[r]\ar[ur]^{l_{k+3}} & D^{k+1}\times X}     \]
Therefore, the obstruction $o_{k+3}(\gamma^{\circ}_{\bQ},\id)$ for the relative lifting problem 
 \begin{equation}\label{eq: rational extension outside a disk}
 \xymatrix{(D^{k+1}\times X)-\Inte(D^{k+5})\ar[r]^-{\gamma^{\circ}_{\bQ}}\ar@{^{(}->}[d]& D^{k+1}\times P^{X}_{\bQ}\ar[d]^{\id\times \pi^{X}_{\bQ}}\\
D^{k+1}\times X\ar@{=}[r]\ar@{-->}[ur] & D^{k+1}\times X}       
 \end{equation}
is nonempty and has no indeterminacy. 
\begin{Lemma}\label{lem: disk obstruction nonvanishing}
The obstruction class     $o_{k+3}(\gamma^{\circ}_{\bQ},\id)\in \pi_{k+4}(F_{\bQ})$ is nonzero.
\end{Lemma}
\begin{proof} By the alternative definition of obstruction class via relative CW structure of $(D^{k+1}\times X, (D^{k+1}\times X)-\Inte(D^{k+5}))$, we see that the element $o_{k+3}(\gamma^{\circ}_{\bQ},\id)$ is mapped to $\alpha_{\eta}$ under the isomorphism $\pi_{k+4}(F_{\bQ})\cong \pi_{k+4}(\Top(4))\otimes \bQ$. . Hence it is nonzero.
\end{proof}

\begin{proof}[Proof of Proposition \ref{prop: nonextension nonempty boundary}]The case $k=2$ has been proved so we may assume $k\geq 3$. Since $\cA^{\text{even}}_{3}=\cA^{\text{even}}_{4}=0$, we have $k\geq 5$.

We consider the tower (\ref{eq: tower F}). For $m=1,2,3$, $Z^F_{m+k}\to Z^F_{m+k-1}$ is a fibration with fiber $K(\pi_{m+k}F_{\bQ},m+k)$. We claim that the induced map 
\[
H^{k+5}(Z^F_{m+k-1};\pi_{m+4}(F_{\bQ}))\to 
H^{k+5}(Z^F_{m+k};\pi_{m+4}(F_{\bQ}))
\]
is surjective. To see this, we note that since $F_{\bQ}$ is 4-connected, $H^{i}(Z^F_{n};G)=0$
for any abelian group $G$, any $n\geq 1$ and any $1\leq i\leq 4$. Note that $H^{*}(K(\pi_{n}F_{\bQ},n);\mathbb{Z})$ is supported in degree $n$ when $n$ is odd and degree $ln$ with $l\geq 1$ when $n$ is even. In particular, we have
\[
H^{k+i}(K(\pi_{m+k}F_{\bQ},m+k);\mathbb{Z})=0,\quad \forall 1\leq m\leq 3\text{ and } m<i\leq 5 
\]
The claim is then proved by a straightforward application of the Serre spectral sequence. From this claim, we see that the $\bfk$-invariant $k_{k+3}(F_{\bQ}/\ast)$ is pulled back from some class in $H^{k+5}(Z^{F}_{k})$.

By Lemma \ref{lem: bundle factor though S4}, the fibration $P^{X}_{\bQ}\to X$ is pulled back from the fiberation $F_{\bQ}\to \ast$. By funtoriality of the $\bfk$-invariant, we see that 
$\bfk_{k+3}(P^{X}_{\bQ},X)$ is pulled back from a class \[\bfk'\in H^{k+5}(Z^{X}_{k};\pi_{k+4}(F_{\bQ})).\]
Note that for any abelian group $G$ and any $i\leq k$, we have 
\[H^{i}(D^{k+1}\times X, \partial (D^{k+1}\times X);G)=H_{k+5-i}(X;G)=0.\]
This allows us to apply Lemma \ref{lem: obstruction independence with lift} (2) to conclude that $o_5(\gamma^{\partial}_{\bQ},\id)$ has no indeterminacy. We may assume $o_5(\gamma^{\partial}_{\bQ},\id)$ is nonempty because otherwise a lifting
$l_{l+3}: D^{k+1}\times X\to Z^{X}_{k+3}$ does not exist and the proof is finished. By naturality of the obstruction classes, we see that $o_{3}(\gamma^{\partial}_{\bQ},\id)$ is mapped to $o_{k+3}(\gamma^{\circ}_{\bQ},\id)$ under the isomorphism 
\[
H^{k+5}(D^{k+1}\times X,\partial (D^{k+1}\times X);\bQ)\cong H^{k+5}(D^{k+1}\times X, D^{k+1}\times X-\Inte(D^{k+5});\bQ).
\]
By Lemma \ref{lem: disk obstruction nonvanishing}, this obstruction class is nonzero. So the proof is finished.
\end{proof}

We are left to prove Proposition \ref{prop: nonextension nonempty boundary}. From now on, we assume $X$ is a closed manifold with $\sigma(X)=0$, and $k\geq 5$. We pick a base point $b\in X^\circ$. We may assume that $b$ is sent to $b_0:=\infty\in S^4$ under the pinch map $p:X\to S^4$ defined in (\ref{eq: pinch}).

We restrict the bundle $P^{X}_{\bQ}$ to get the bundle $F_{\bQ}\to P^{b}_{\bQ}\to D^{k+1}\times \{b\}$. Let $\gamma^{b}$ be any section of this bundle that coincide with $\gamma^{\partial}_{\bQ}$ on $S^{m}\times \{b\}$. We consider the relative lifting problem 
 \begin{equation}\label{eq: rational lifting with a section at infty}
 \xymatrix{(\partial D^{k+1}\times X)\cup (D^{k+1}\times \{b\}) \ar[r]^-{\gamma^{\partial}_{\bQ}\cup \gamma^{b}}\ar@{^{(}->}[d]& D^{k+1}\times P^{X}_{\bQ}\ar[d]^{\id\times \pi^{X}_{\bQ}}\\
D^{k+1}\times X\ar@{-->}[ur]\ar@{=}[r] & D^{k+1}\times X}    
 \end{equation}
We have the following three lemmas.
\begin{Lemma}\label{lem: o-gamma-b no indeterminacy} For any $\gamma^{b}$, the 
 obstruction class  
 \[
o_{k+3}(\gamma^{\partial}_{\bQ}\cup \gamma^{b},\id)\in \pi_{k+4}(F_{\bQ}).
 \]
of (\ref{eq: rational lifting with a section at infty}) has no indeterminacy.
 \end{Lemma}

\begin{Lemma}\label{lem: obstruction class nonvaishing for gamma-b-0}
Let $\gamma^{b}_{0}$ be the restriction of $\gamma^{\circ}_{\bQ}$ to $D^{k+1}\times \{b\}$. Then we have $o_{k+3}(\gamma^{\partial}_{\bQ}\cup \gamma^{b}_{0},\id)\neq 0$.
\end{Lemma}

\begin{Lemma}\label{lem: obstruction class independent with gamma-b}
For any choice of $\gamma_{b}$, one has 
\begin{equation}\label{eq: obstruction independent of section at infinity}
o_{k+3}(\gamma^{\partial}_{\bQ}\cup \gamma^{b}_{0},\id)=o_{k+3}(\gamma^{\partial}_{\bQ}\cup \gamma^{b},\id).  \end{equation}
 \end{Lemma}

\begin{proof}[Proof of Proposition \ref{prop: signature zero}] By Lemma \ref{lem: o-gamma-b no indeterminacy}, Lemma \ref{lem: obstruction class nonvaishing for gamma-b-0} and Lemma \ref{lem: obstruction class independent with gamma-b}, the obstruction class $o_{k+3}(\gamma^{\partial}_{\bQ}\cup \gamma^{b},\id)$ has no indeterminacy and is nonzero for any choice of $\gamma^{b}$. Hence the lifting problem (\ref{eq: rational lifting with a section at infty}) has no solution.
\end{proof}

 \begin{proof}[Proof of Lemma \ref{lem: o-gamma-b no indeterminacy}] 
We consider the tower (\ref{eq: tower S}). Since $S^4$ is 3-connected and $F_{\bQ}$ is 4-connected, we see that $H^{i}(Z^{S}_{n};G)=0$ for any abelian group $G$, any $n\geq 1$ and any $1\leq i\leq 3$. By exactly the same argument as the proof of Proposition \ref{prop: nonextension nonempty boundary}, we see that the $\bfk$-invariant $\bfk_{k+3}(P^{D}_{X}/X)$ is pulled back from a class 
\[
\bfk'\in H^{k+5}(Z^{X}_{k+1};\pi_{k+4}(F_{\bQ})).
\]
Note that for any abelian group $G$ and any $0\leq m\leq k+1$ we have 
\[
H^{m}(D^{k+1}, (\partial D^{k+1}\times X)\cup (D^{k+1}\times \{b\});G)=0.
\]
So we can apply Lemma \ref{lem: obstruction independence with lift} (2) to finish the proof. \end{proof}

\begin{proof}[Proof of Lemma \ref{lem: obstruction class nonvaishing for gamma-b-0}] By functoriality of the obstruction class, $o_{k+3}(\gamma^{\partial}_{\bQ}\cup \gamma^{b}_{0},\id)$ is sent to $o_{k+3}(\gamma^{\circ}_{\bQ},\id)$ under the isomorphism  beween
\[
H^{k+5}(D^{k+1}\times X, (\partial D^{k+3}\times X)\cup (D^{k+3}\times \{b\});\pi_{k+4}(F_{\bQ}))\]
and 
\[H^{k+5}(D^{k+3}\times X, (D^{k+3}\times X)-\Inte(D^{7});\pi_{k+4}(F_{\bQ})).
\]
So $o_{k+3}(\gamma^{\partial}_{\bQ}\cup \gamma^{b}_{0},\id)\neq 0$ by Lemma \ref{lem: disk obstruction nonvanishing}.
\end{proof}

\begin{proof}[Proof of Lemma \ref{lem: obstruction class independent with gamma-b}] We pick the base point $*\in F_{\bQ}$ to be the image of $[e]\in F$ under the localization map $F\to F_{\bQ}$. We fix a trivilization of $T_{b}X$. This gives an identification between $F_{\bQ}$ and the fiber of $P^{X}_{\bQ}\to X$ at $b$. Under this identification, $\gamma^{b}$ is just a map 
\[
\gamma^{b}:(D^{k+1},\partial D^{k+1})\to (F_{\bQ},*).
\]
and $\gamma^{b}_{0}$ is the constant map.

We first assume that $\gamma^{b}$ factors through a map
\[
\beta^b: (D^{k+1},\partial D^{k+1})\to (\Top(4),*).
\]
Then $\beta^{b}$ induces a homeomorphism 
\[
\beta_{1}: D^{k+1}\times \mathbb{R}^{4}\xrightarrow{\cong_{\Top}} D^{k+1}\times \mathbb{R}^{4}
\]
that covers the the identity map on $D^{k+1}$ and fixes $\partial D^{k+1}\times \mathbb{R}^4$ pointwisely. After a fiberwise compactification, we obtain the homeomorphism 
\[
\beta_2: D^{k+1}\times (\mathbb{R}^{4}\cup \infty)\xrightarrow{\cong_{\Top}} D^{k+1}\times (\mathbb{R}^{4}\cup \infty)
\]
We consider the inversion 
\[\iota: S^{4}=\mathbb{R}^{4}\cup \infty\to \mathbb{R}^{4}\cup \infty,\quad x\mapsto \frac{x}{|x|^2}, 
\]
which interchanges $0$ and $\infty$. Then the bundle automorphism 
\[
\beta_3:=(\iota\times \id)\circ\beta_2\circ (\iota\times \id): D^{k+1}\times S^{4}\to D^{k+1}\times S^4
\]
fixes the section at $\infty$. 

Let \[\Top(4)\to \Fr^{S^4}_{\Top}\to S^4\] be the frame bundle of $S^4$ (with structure group extended to $\Top(4)$). Let \[F\to \bar{P}^{S}\to S^4\] be the associated bundle and let \[F_{\bQ}\to \bar{P}_{\bQ}^{S}\to S^4\] be the rationalization. Up to homotopy relative to $\partial D^{k+1}\times \Fr^{S^4}_{\Top}$, the  the bundle automorphism $\beta_3$ canonically induces a bundle automorphism 
\[
D^{k+1}\times \Fr^{S^4}_{\Top}\to D^{k+1}\times \Fr^{S^4}_{\Top}
\]
and hence an automorphism 
\[
\beta_{4}: D^{k+1}\times \bar{P}_{\bQ}^{S}\to D^{k+1}\times \bar{P}_{\bQ}^{S}
\]
By Lemma \ref{lem: bundle factor though S4}, the bundle $P^{X}_{\bQ}\to X$ is pulled back from the bundle $\bar{P}_{\bQ}^{S}\to S^4$. (Note that $\bar{P}_{\bQ}^{S}$ and $P_{\bQ}^{S}$ are different.) Hence $\beta_4$ induces a bundle automorphism 
\[
\beta_{5}: D^{k+1}\times P_{\bQ}^{X}\to D^{k+1}\times P_{\bQ}^{X}.
\]
Unwinding definitions, we see that $\beta_{5}$ sends the section $\gamma^{\partial}\cup \gamma^{b}_{0}$ to $\gamma^{\partial}\cup \gamma^{b}$ (up to homotopy). To this end, (\ref{eq: obstruction independent of section at infinity}) follows from  functorality of obstruction classes.

Now we consider a general $\gamma^{b}$. We pick a map $\rho: D^{k+1}\to D^{k+1}$ of nonzero degree and consider the lifting problem 
\[
 \xymatrix{(\partial D^{k+1}\times X)\cup (D^{k+1}\times \{b\}) \ar[r]^-{\gamma^{\partial}_{\bQ}\cup \gamma^{b}}\ar@{^{(}->}[d]&D^{k+1}\times P^{X}_{\bQ}\ar[r]^{\rho\times \id} &D^{k+1}\times P^{X}_{\bQ}\ar[d]^{\id\times \pi^{X}_{\bQ}}\\
D^{k+1}\times X\ar@{-->}[urr]\ar^{\rho\times \id}[rr] & & D^{k+1}\times X}    
 \]
The obstruction class is just given by $\deg(\rho)\cdot o_{k+3}(\gamma^{\partial}_{\bQ}\cup \gamma^{b},\id)$. Note that the composition 
\[
\Top(4)\to F\to F_{\bQ}
\]
induces an isomorphism on $\pi_{k+1}(-)\otimes\bQ$, so by choosing $\rho$  with a suitable degree, the composition 
\[
\rho\circ \gamma^{b}: (D^{k+1},\partial D^{k+1})\to (F_{\bQ},*)
\]
factors through a map $(D^{k+1},\partial D^{k+1})\to (\Top(4),*)$. By previous discussion, we have 
\[
\deg(\rho)\cdot o_{k+3}(\gamma^{\partial}_{\bQ}\cup \gamma^{b},\id)=\deg(\rho)\cdot o_{k+3}(\gamma^{\partial}_{\bQ}\cup \gamma^{b}_0,\id)\in \pi_{k+4}(F_{\bQ}).
\]
Since $\pi_{k+4}(F_{\bQ})$ is torsion free, (\ref{eq: obstruction independent of section at infinity}) is again proved.
\end{proof}
\begin{remark} For a closed 4-manifold $X$ with nonvanishing signature, the indeterminacy of $o(\gamma^{\partial}_{\bQ},\id)$ is contained the image of the Whitehead product 
\[
[-,-]:\pi_{4}(\BO) \otimes \pi_{k+2}(\BTop(4))\to \pi_{k+5}(\BTop(4))
\]
So if one can show that the image of $[-,-]$ is contained in the kernel of $I^{\Top}(-)$, then our argument would  show that $\pi_{k}(\Sm(X))\otimes \bQ\geq \cA^{\even}_{k}$.  
\end{remark}

\bibliographystyle{plain}
\bibliography{mainref}

\end{document}